\documentclass[11pt]{article}
\usepackage{amsmath,amssymb,amsthm,amsfonts,latexsym}
\usepackage{mathrsfs,bbm,bm}
\usepackage{xcolor}
\usepackage{enumitem}
\allowdisplaybreaks

\usepackage[letterpaper,margin=1in]{geometry}

\DeclareMathOperator*{\esup}{ess\,sup}
\DeclareMathOperator*{\einf}{ess\,inf}

\newtheorem{thm}{Theorem}[section]
\newtheorem{corollary}[thm]{Corollary}
\newtheorem{lemma}[thm]{Lemma}
\newtheorem{proposition}[thm]{Proposition}
\newtheorem{defn}[thm]{Definition}
\theoremstyle{definition}
\newtheorem{remark}[thm]{Remark}
\theoremstyle{plain}

\numberwithin{equation}{section}

\def\E{{\mathbb E}}
\def\P{{\mathbb P}}
\def\R{\mathbb R}

\def\A{\mathcal A}
\def\L{\mathcal L}
\def\S{\mathcal S}
\def\wh{\widehat}
\def\e{\mathrm{e}}
\def\1{{\mathbbm 1}}
\def\eps{\varepsilon}

\begin{document}

\title{\bf Dirichlet heat kernel estimates for fractional Laplacian under non-local perturbation}
\author{{\bf Zhen-Qing Chen}
 \quad and \quad {\bf Ting Yang}\thanks{
The research of this author is supported by NSFC (Grant No. 12371143).}}

\date{}
\maketitle

\begin{abstract} \small
For $d\geq 1$
and $0<\beta<\alpha<2$, consider a family of non-local
operators $\mathcal{L}^{b}=\Delta^{\alpha/2}+\mathcal{S}^{b}$ on $\R^d$, where
\[
\mathcal{S}^{b}f(x):=\lim_{\varepsilon\to
0}\mathcal{A}(d,-\beta)\int_{ \{z\in \R^d: |z|>\varepsilon\}}
(f(x+z)-f(x))\frac{b(x,z)}{|z|^{d+\beta}}\,{\mathrm{d}z},
\]
and $b(x,z)$ is a bounded measurable function on
$\R^{d}\times\R^{d}$ with $b(x,z)=b(x,-z)$ for every
$x,z\in\R^{d}$.
 Here ${\mathcal{A}}(d, -\beta)$
is the normalizing constant for which $\S^b=-(-\Delta)^{\beta/2}$
when $b(x, z)\equiv 1$. Chen and Wang \cite{CW} proved
that if
$b(x, z) \geq -\frac{{\mathcal{A}}(d, -\alpha)}
{{\mathcal{A}}(d, -\beta)}\, |z|^{\beta -\alpha}$ for a.e. $z\in\R^{d}$,
then
$\L^b$ admits a unique fundamental solution
 $p^{b}(t,x,y)$, which determines a conservative Feller process $X^{b}$ with the strong Feller property.
In this paper, we study
the subprocess $X^{b,D}$ of $X^{b}$ killed upon leaving
a  $C^{1,1}$ open set $D\subset \R^d$. We
establish explicit small-time upper and lower bounds
for its transition density function $p^{b}_{D}(t,x,y)$.
These estimates
are uniform in $b$: all constants depend on $b$ only through an upper bound of $\|b\|_{\infty}$. Moreover, $b$ enters the form of the upper bound only through the size of its large positive part, $\esup_{x\in\R^{d},\,|z|>\lambda}b(x,z)^{+}$ for a fixed $\lambda>0$.
When the jump kernel of $\mathcal{L}^b$
satisfies matching upper and lower bounds,
either of $\alpha$-stable type or of mixed $\alpha$- and $\beta$-stable type,
the small-time estimates for $p^{b}_{D}(t,x,y)$ are sharp.
 For bounded $D$, we also obtain large-time estimates
 governed by the principal eigenvalue
    of the generator  of $X^{b,D}$.
Under a non-degeneracy condition on the jump kernel, we further establish sharp two-sided estimates on the Green function of $X^{b, D}$ in bounded $C^{1,1}$ open sets.
\end{abstract}

\medskip

\noindent\textbf{Mathematics Subject Classification (2020).} Primary 60J35, 47G20, 60J76; Secondary 47D07, 35K08.

\medskip

\noindent\textbf{Keywords and Phrases.}  Rotationally symmetric $\alpha$-stable process, fractional Laplacian,
non-local perturbation, non-symmetric operator,
Dirichlet heat kernel, Green function,
derivative estimate.

\medskip

\section{Introduction}

Discontinuous Markov processes and non-local operators are of central importance in probability and analysis.
Non-Gaussian jump processes model a wide range of physical and economic systems, such as anomalous diffusion in physics \cite{BG90, KSZ96, MK00}, the transport of contaminants in heterogeneous media \cite{BCDS06}, and asset prices with jumps in mathematical finance \cite{CT04}.
The infinitesimal generator of a discontinuous Markov process in $\R^{d}$
is not a differential operator. It is a non-local, or
integro-differential, operator. The basic example is the rotationally
 symmetric (also called isotropic) $\alpha$-stable process in $\R^{d}$
 with
$\alpha \in (0, 2)$, whose generator is the fractional Laplacian $c
\Delta^{\alpha /2}:=- c\, (-\Delta)^{\alpha /2}$. Non-local operators also arise in the theory of PDE, for instance in obstacle problems; see
 \cite{CSS, Silvestre} and the references
therein.

The transition density function of a Markov process, also called its heat kernel,
encodes all the information about the process.
Explicit formulas are available only in a few special cases.
Sharp two-sided estimates are therefore the natural substitute.
For diffusions such estimates are classical. For jump processes in $\R^d$, two-sided heat
kernel estimates have been studied
systematically only since around 2000; see, for example, \cite{CK1, Chen-Kumagai, CKK2}.
The transition densities of jump processes killed upon leaving an open set,
called Dirichlet heat kernels, have been studied more recently still.
We refer to \cite{Chen} for a survey of the symmetric theory.

Many of the systems mentioned above evolve in a confined region. A particle in a vessel is absorbed at the walls; a contaminant is removed once it leaves a domain; a financial position triggers default once it leaves a solvency region. For jump processes, the exit may also occur by a jump, without the process touching the boundary. Such situations are described by the subprocess killed upon leaving an open set $D$, and its transition density, the Dirichlet heat kernel, is the central object. For the operators studied in this paper, it is the function $p^{b}_{D}(t,x,y)$ defined in \eqref{defiforpbd} below.

Dirichlet heat kernels carry important information beyond that contained in global heat kernels. The Dirichlet heat kernel $p_{D}(t,x,y)$ of a Markov process $X$ in an open set $D$ is the transition density of the process killed upon leaving $D$, that is,  for $x\in D$ and Borel
 sets $E\subset D$,
$$
\P_{x}(X_{t}\in E,\,t<\tau_{D})=\int_{E}p_{D}(t,x,y)\,\mathrm{d}y,
$$
 where $\tau_{D}$ is the first exit time of $X$ from $D$. In particular, $\P_{x}(\tau_{D}>t)=\int_{D}p_{D}(t,x,y)\,\mathrm{d}y$, so sharp estimates on $p_{D}$ give sharp estimates on the survival probability. Analytically, $p_{D}$ is the fundamental solution of the parabolic problem $\partial_{t}u=\mathcal{L}u$ in $D$ with zero exterior condition on $D^{c}$, where $\mathcal{L}$ is the generator of the process. Its estimates therefore describe solutions of this problem up to the boundary.
Integrating $p_{D}$ in $t$ gives the Green function, and integrating the latter in $y$ gives the expected exit time.
 For bounded $D$,
 under suitable compactness and irreducibility assumptions,
 the large-time behavior of $p_{D}$ is governed by the principal eigenvalue and eigenfunction of the generator, so sharp estimates yield the decay rate of the semigroup and the boundary behavior of the
  principal eigenfunction.
Finally, the boundary factor in a sharp estimate identifies the exact rate at which the process is killed near $\partial D$. For non-local operators this rate depends on the order of the operator, and it cannot be read off from the global heat kernel.

Most existing heat kernel estimates for discontinuous Markov processes
concern symmetric processes.
Non-symmetric perturbations of the fractional Laplacian by gradient operators are studied in \cite{BJ} for the whole space and in \cite{CKS2} for open sets.
Heat kernels of non-symmetric jump diffusions in the whole space are studied in \cite{CZ}.
Chen and Wang \cite{CW} studied a class of non-symmetric non-local operators:
fractional Laplacians under non-local perturbations.
Let $d\geq 1$
and $0<\beta<\alpha<2$. Consider the non-local
operator $\mathcal{L}^{b}=\Delta^{\alpha/2}+\mathcal{S}^{b}$,
 where
\begin{equation}
\mathcal{S}^{b}f(x):=\lim_{\varepsilon\to 0}\mathcal{A}(d,-\beta)
\int_{ \{z\in \R^d: |z|>\varepsilon\}}
(f(x+z)-f(x))\frac{b(x,z)}{|z|^{d+\beta}}\,{\mathrm{d}z},
\nonumber
\end{equation}
$\mathcal{A}(d,-\beta)=\beta
2^{\beta-1}\pi^{-d/2}\Gamma((d+\beta)/2)\Gamma(1-\beta/2)^{-1}$ is the normalizing
constant for which $\S^b=\Delta^{\beta/2}:=-(-\Delta)^{\beta/2}$
when $b(x, z)\equiv 1$, and
$b(x,z)$ is a bounded measurable function on
$\R^{d}\times\R^{d}$ such that
\begin{equation}
b(x,z)=b(x,-z)  \quad \hbox{for }  x,z\in\R^{d}.\label{condi1}
\end{equation}
Equivalently,
\[
\L^bf (x)=\int_{\R^d} \big( f(y)-f(x)-\langle\nabla f(x), y-x \rangle
\1_{\{|y-x|\leq 1 \}} \big)
j^b(x, y){\mathrm{d}y},
\]
where
\begin{equation}
j^{b}(x,y)
= \frac{{\mathcal{A}}(d, -\alpha)}{|y-x|^{d+\alpha}}
\big( 1+ \frac{{\mathcal{A}}(d, -\beta)}{{\mathcal{A}}(d, -\alpha)}\,
b(x, y-x)\, |y-x|^{\alpha-\beta} \big)
\1_{\{y-x\not=0\}}.
\nonumber
\end{equation}
When $b\equiv 0$, we write $j(x, y)$ for $j^0(x, y)$.

We recall the main results of Chen and Wang \cite{CW}. Assume that
\begin{equation}\label{condi2}
\hbox{for every } x\in \R^d, \quad
j^b(x,y)\ge 0 \quad \hbox{for a.e. } y\in \R^{d}
\end{equation}
(that is, for every $x\in \R^d$, $b(x,z)\ge -\frac{\mathcal{A}(d,-\alpha)}{\mathcal{A}(d,-\beta)}\,|z|^{\beta-\alpha}$
for a.e.\ $z\in\R^{d}$). Then $\L^b$ admits a unique fundamental solution
$p^{b}(t,x,y)$ \cite[Theorem 1.1]{CW}. This kernel is jointly continuous
on $(0,\infty)\times\R^{d}\times\R^{d}$, and it is strictly positive by the lower bound in \cite[(1.22)]{CW}.
It determines a conservative Feller process $X^{b}$ with the strong Feller property
on the Skorokhod space
$\mathbb{D}([0,\infty),\R^{d})$ such that
\[\E_x\big[f(X^{b}_{t})\big]=\int_{\R^{d}}f(y)p^{b}(t,x,y){\mathrm{d}y}\]
for every bounded measurable function $f$ on $\R^d$; see \cite[Theorem 1.3]{CW} and the paragraph following it.
The jumps of $X^{b}$ are described by its L\'{e}vy system $(j^{b}(x,y){\mathrm{d}y},t)$ \cite[(1.25)--(1.26)]{CW}: for every non-negative measurable function $f$ on
$\R_{+}\times \R^{d}\times \R^{d}$ vanishing
on $\{(s,y,y):s\ge 0,y\in\R^{d}\}$, every $ x\in
\R^{d}$ and stopping time $T$ (with respect to the minimal
admissible filtration of $X^{b}$),
\begin{equation}
\E_x\Big[\sum_{
 s\le T
}f(s,X^{b}_{s-},X^{b}_{s})\Big]=\E_x\Big[\int_{0}^{T}\int_{\R^{d}}f(s,X^{b}_{s},y)j^{b}(X^{b}_{s},y){\mathrm{d}y}\,{\mathrm{d}s}\Big].\label{levysystemforxb}
\end{equation}
Moreover, $(X^b_t, \P_x, x\in \R^d)$ is the unique solution
to the martingale problem for $(\L^b, \mathcal{S}(\R^d))$ \cite[Theorems 1.3 and 1.6]{CW}.
 Sharp two-sided estimates on $p^b(t, x, y)$ are derived
in several explicit forms in \cite[Theorem 1.3, Corollary 1.4 and Theorem 1.5]{CW}.
 Proposition \ref{P:2.1} below summarizes those that will be needed in this paper.
These results describe $\L^b$ on the whole space.
For the same operator, Chen, Ren and Yang \cite{CRY16} proved a uniform boundary Harnack principle and gradient estimates for non-negative harmonic functions. Shi and Yi \cite{SY22} studied intrinsic ultracontractivity for this and related operators.

This paper studies $\L^b$ in
 $C^{1,1}$ open sets with zero exterior condition.
 We establish small-time upper and lower bounds on its Dirichlet heat kernel, which become sharp two-sided estimates when the jump kernel satisfies matching upper and lower bounds. Under a non-degeneracy condition on the jump kernel, we also obtain sharp two-sided estimates on its Green function in bounded $C^{1,1}$ open sets.
Before stating the results, we fix notation and recall three model cases.

\medskip
\noindent{\bf Notation.} We use ``:=" to denote a
definition. We define $a\wedge b=\min\{a,b\}$, $a\vee
b=\max\{a,b\}$ and $\log^{+}u:=\log(u\vee 1)$ for $u>0$; $\e$ denotes Euler's number. For any two positive functions $f$ and $g$,
$f\stackrel{c}\lesssim g$ means that there is a positive constant
$c$ such that $f\le c g$ on their common domain of definition, and
$f\stackrel{c}\asymp g$ means that $c^{-1}g\le f\le cg$. We also
write ``$\lesssim$" and ``$\asymp$" if $c$ is unimportant or
understood. We use $B(x,r)$ to denote the open ball centered at $x$
with radius $r>0$. Let $\delta_{D}(x)$ denote the Euclidean distance
between $x$ and $\partial D$.
Capital letters $C_{0},\ C_{1},\ C_{2},\cdots$ denote constants in the statements of results.
Lower case constants $c_{0},\ c_{1},\ c_{2},\cdots$ may change from one appearance to another. We write ${\mathrm{d}x}$ for the Lebesgue measure in $\R^{d}$ and
$\mathrm{diam}(D):=\sup\{|x-y|:x,y\in D\}$
for the diameter of a set $D\subset \R^{d}$.
For a complex number $z$,
we use ${\rm Re}\, z$ and ${\rm Im}\, z$ to denote its real and imaginary parts, respectively.
When $f$ is a function of several variables, we sometimes write $\mathcal{S}^{b}_{x}f$ instead of $\mathcal{S}^{b}f(x)$ to indicate the variable on which $\mathcal{S}^{b}$ acts.

\medskip
\noindent{\bf Model cases.} The family of processes $X^{b}$ contains
non-local perturbations
of several important L\'{e}vy processes. Their heat kernels show the range of behavior that any estimate on $p^b(t, x, y)$ must capture.
When $b\equiv 0$,
$X^{0}$
is the rotationally symmetric $\alpha$-stable
process on $\R^d$. We denote its transition density function by $p(t,x,y)$.
 It is well known that
\begin{equation}\label{p0}
p(t,x,y)\stackrel{c_{0}}{\asymp} t^{-d/\alpha}\wedge \frac{t}{|x-y|^{d+\alpha}}
\quad  \hbox{on }  (0, {\infty})\times \R^{d}\times\R^{d},
\end{equation}
for some positive constant $c_{0}=c_{0}(d,\allowbreak \alpha)$.
 When $b\equiv a$ for some constant $a>0$,
$\mathcal{L}^{b}=\Delta^{\alpha/2}+a\Delta^{\beta/2}$ and $X^{b}$ is
the independent sum of a rotationally symmetric $\alpha$-stable process and a
scaled rotationally symmetric $\beta$-stable process. Denote by $p_{a}(t,x,y)$ the corresponding
transition density. It is proved in \cite{Chen-Kumagai} (see also \cite[(1.9)]{CW}) that
\begin{equation}
p_{a}(t,x,y)
 \stackrel{c_{1}}{\asymp}
\big(t^{-\frac{d}{\alpha}}\wedge(at)^{-\frac{d}{\beta}}\big)\wedge
\big(\frac{t}{|x-y|^{d+\alpha}}+\frac{at}{|x-y|^{d+\beta}}\big)
 \quad \hbox{on }  (0,{\infty})\times \R^{d}\times\R^{d},
 \label{pa}
\end{equation}
 for some positive constant $c_{1}=
 c_{1}(d,\allowbreak \alpha, \allowbreak \beta)
 $.
When
$b(x,z)=-\frac{\mathcal{A}(d,-\alpha)}{\mathcal{A}(d,-\beta)}
|z|^{\beta-\alpha}\1_{\{|z|\ge
1\}}$,  $X^{b}$ is a truncated rotationally symmetric $\alpha$-stable process
with L\'{e}vy intensity
$\mathcal{A}(d,-\alpha) |x|^{-d-\alpha}\1_{\{|x|<1\}}{\mathrm{d}x}$.
Denote by
$\bar{p}_{1}(t,x,y)$ its transition
density function. It is proved in
\cite{Chen-Kim-Kumagai} that for $t\in (0,1]$ and $|x-y|\le 1$,
\[\bar{p}_{1}(t,x,y)
 \stackrel{c_{2}}{\asymp}
t^{-\frac{d}{\alpha}}\wedge
\frac{t}{|x-y|^{d+\alpha}},\]
while for $t\in (0,1]$ and $|x-y|>1$,
\[c_{3}\big(\frac{t}{|x-y|}\big)^{c_{4}|x-y|}\le \bar{p}_{1}(t,x,y)\le
c_{5}\big(\frac{t}{|x-y|}\big)^{c_{6}|x-y|}\] for some constants
$c_{i}=c_{i}(d,\allowbreak \alpha)>0,\ i=2,\cdots,6.$
In the first two cases the heat kernel decays polynomially off the diagonal.
In the third case $j^{b}$ vanishes for jumps of size at least $1$, and the off-diagonal decay is much faster.
Condition \eqref{condi2} admits both regimes. This distinction determines the form of the lower bounds in Theorem \ref{T:1.2}.

\medskip
\noindent{\bf Killed process.} For an open set $D \subset \R^{d}$, define $\tau^{b}_{D}:=\inf\{t>0:X^{b}_{t}\not\in D\}$.
Let $X^{b,D}$ denote the subprocess of $X^{b}$ killed upon leaving  $D$;
 that is, $X^{b,D}_{t}(\omega)=X^{b}_{t}(\omega)$ if $t<\tau^{b}_{D}(\omega)$ and $X^{b,D}_{t}(\omega)=\partial $ if $t\ge \tau^{b}_{D}(\omega)$, where $\partial$ is a cemetery state.
 Every function $f$ is extended to $\partial$ by setting $f(\partial)=0$.
For $(t, x, y)\in (0, \infty) \times D\times D$, define
\begin{equation}\label{defiforpbd}
p_{D}^{b}(t,x,y):=p^{b}(t,x,y)-\E_x\big[
p^{b}(t-\tau^{b}_{D},X^{b}_{\tau^{b}_{D}},y);\tau^{b}_{D}<t\big].
\end{equation}
Then $p^{b}_{D}(t,x,y)$ is the transition density of the subprocess $X^{b,D}$, and is called the Dirichlet heat kernel of $X^{b}$ in $D$. It is jointly continuous on $(0,\infty)\times D\times D$ by Theorem \ref{them3.4} below.
Let $L^{2}(D;\mathrm{d}x)$ be the space of all complex-valued measurable functions on $D$ which are square-integrable with respect to the Lebesgue measure.
By the estimates on $p^b(t, x, y)$ (see Proposition \ref{P:2.1}(i) and
Theorem \ref{them3.4} below), the transition semigroup $\{P^{b,D}_t; t\geq 0\}$ of $X^{b, D}$,
defined by $P^{b, D}_t f (x):= \E_x [ f (X^{b, D}_t)]$, is a strongly continuous semigroup in $L^2(D; {\mathrm{d}x})$.
Let $\L^{b, D}$ denote the infinitesimal generator of  $\{P^{b,D}_t; t\geq 0\}$ in $L^2(D; {\mathrm{d}x})$.
Formally, $\L^{b, D}$ is the operator $\L^b$ in $D$ with zero Dirichlet exterior condition
on $D^c$. The (complex) spectrum of $\L^{b, D}$ is denoted by $\sigma (\L^{b, D})$;
see Section \ref{sec8} for its definition.
Set $\lambda^{b,D}_{1}:=-\sup\mathrm{Re}\,\sigma(\mathcal{L}^{b,D})$. When $D$ is bounded, $-\lambda^{b,D}_{1}$ is the principal eigenvalue of $\L^{b,D}$, and it has a non-negative eigenfunction; see Theorem \ref{T:7.1}.

\medskip
\noindent{\bf Boundary behavior.} The central question is how killing at $\partial D$ affects the global estimates above.
For the killed rotationally symmetric $\alpha$-stable process, sharp two-sided Dirichlet heat kernel estimates were obtained in \cite{CKS1}; see also \cite{BGR1}.
Analogous results hold for $\Delta^{\alpha/2}+\Delta^{\beta/2}$ \cite{CKS3} and for the fractional Laplacian with gradient perturbation \cite{CKS2}.
They have been extended to large classes of symmetric L\'{e}vy processes; see \cite{CKSrel, CKSrot, BGR2} and the references therein.
In each case, the Dirichlet heat kernel is comparable to  the global estimate multiplied by one boundary factor for each space variable.
This factor measures the probability that the process survives in $D$ up to time $t$.
These results are obtained for $C^{1,1}$ open sets.
This regularity provides interior and exterior balls of uniform size at every boundary point.
These balls are the basic tool for obtaining the explicit boundary decay rate.
We therefore work with the following class of open sets.

\begin{defn}\label{C11set}
Let $d\geq 2$. An
open set $D$ in $\R^{d}$ is said to be $C^{1,1}$ if there
exist a localization radius $R_{0}>0$ and a constant
$\Lambda_{0}>0$ with the following property. For every $Q\in \partial D$, there exist a
$C^{1,1}$ function $\phi=\phi_{Q}:\R^{d-1}\to\R$
satisfying $\phi(0)=0,\ \nabla\phi(0)=0,\ \|\nabla
\phi\|_{\infty}\le\Lambda_{0},\
|\nabla\phi(x)-\nabla\phi(y)|\le\Lambda_{0}|x-y|$, and an
orthonormal coordinate system $CS_{Q}$ with its origin at $Q$ such
that
\[
B(Q,R_{0})\cap D=\{y=(\tilde{y},y_{d})\mbox{ in } CS_{Q}:|y|<R_{0},y_{d}>\phi(\tilde{y})\}.
\]
The pair $(R_{0},\Lambda_{0})$ is called
the $C^{1,1}$  \textit{characteristics}
of $D$.
An open set $D\subset\R$ is said to be $C^{1,1}$ if it is a union of finitely or countably many disjoint open intervals and there is $R_{0}>0$ with the following property: each interval has length at least $R_{0}$, and any two distinct intervals are at distance at least $R_{0}$.
In this case, $R_0$ is called the $C^{1,1}$ {\it   characteristic} of $D$.
\end{defn}
For notational uniformity, we set $\Lambda_{0}:=1$ when $d=1$ and also call $(R_{0},\Lambda_{0})$ the $C^{1,1}$ characteristics of $D$. In dimension one, $\Lambda_{0}$ plays no geometric role, and constants listed as depending on $\Lambda_{0}$ do not depend on it.

A $C^{1,1}$ open set can be unbounded and disconnected, and a bounded $C^{1,1}$ open set has only finitely many connected components.
For a domain $D\subset \R^{d}$ and a constant $\Lambda_{1}\ge 1$, we say that the path distance in $D$ is comparable to the Euclidean distance with characteristic $\Lambda_{1}$ if for every $x,y\in D$, there is a rectifiable curve $l$ in $D$ connecting $x$ to $y$ so that the length of $l$ is no larger than $\Lambda_{1}|x-y|$. Clearly, such a property holds for all bounded $C^{1,1}$ domains, $C^{1,1}$ domains with compact complements and domains above the graphs of bounded $C^{1,1}$ functions.

\smallskip

For $A\in (0,+\infty)$, let $\mathcal{B}_{A}$ be the set of bounded functions $b(x,z)$ on $\R^{d}\times\R^{d}$ satisfying \eqref{condi1} and \eqref{condi2} with $\|b\|_{\infty}\le A$.

\smallskip

For $t>0$ and $x,y\in D$, we set
\[
\Phi_{D}(t,x,y):=\Big(1\wedge \frac{\delta_{D}(x)^{\alpha}}{t}\Big)^{1/2} \Big(1\wedge \frac{\delta_{D}(y)^{\alpha}}{t}\Big)^{1/2},
\]
and, for $a\ge 0$,
\[
f^{a}_{D}(t,x,y):=\Phi_{D}(t,x,y)\Big(t^{-d/\alpha}\wedge\Big(\frac{t}{|x-y|^{d+\alpha}}+\frac{a t}{|x-y|^{d+\beta}}\Big)\Big).
\]
For $\lambda>0$, let $M_{b^{+},\lambda}:=\esup_{x,z\in \R^{d},\, |z|>\lambda}b^{+}(x,z)$, where $b^{+}:=b\vee 0$.

Our main result shows that the boundary behavior of the killed rotationally symmetric $\alpha$-stable process persists under the perturbation $\mathcal{S}^{b}$. Away from the boundary, the estimates are governed by the global heat kernel.

\begin{thm}\label{T:1.2}
Let $d\ge 1$, let $D\subset\R^{d}$ be a $C^{1,1}$ open set with characteristics $(R_{0},\Lambda_{0})$, and let $A,T>0$.
All constants below depend only on $d,\alpha,\beta,R_{0},\Lambda_{0},A,T$ and on the additional parameters indicated. In {\rm (i)}--{\rm (iii)}, the estimates hold for all $b\in\mathcal{B}_{A}$ satisfying the stated conditions, and for every $t\in(0,T]$ and $x,y\in D$.
\begin{enumerate}[label=\upshape(\roman*), ref=\upshape(\roman*)]
\item \emph{(Upper bounds.)} For every $\lambda>0$, there is $C_{0}=C_{0}(\lambda)>0$ such that
\[
p^{b}_{D}(t,x,y)\le C_{0}f^{a}_{D}(t,x,y)\quad\hbox{with } a:=M_{b^{+},\lambda}.
\]
If, for some $\eps\in(0,1]$,
\begin{equation}\label{condi5}
j^{b}(x,y)\le \eps^{-1} |x-y|^{-d-\alpha}\quad\mbox{ for a.e. }x,y\in \R^{d},
\end{equation}
then $p^{b}_{D}(t,x,y)\le C_{1}f^{0}_{D}(t,x,y)$ with $C_{1}=C_{1}(\eps)>0$.

\item \emph{(Lower bounds.)} There are $R_{1}=R_{1}(d,\alpha,\beta,A)>0$ and $C_{2}>0$ such that $p^{b}_{D}(t,x,y)\ge C_{2} f^{0}_{D}(t,x,y)$ whenever $|x-y|<R_{1}$.
For every $\Lambda_{1}\ge 1$, there are $C_{i}=C_{i}(\Lambda_{1})>0$, $i=3,4,5$, such that
\begin{equation}\label{them1.ii}
p^{b}_{D}(t,x,y)\ge C_{3}\,\Phi_{D}(t,x,y)\Big(\frac{C_{4}\, t}{|x-y|}\Big)^{C_{5}|x-y|}
\end{equation}
whenever $|x-y|\ge R_{1}$ and $x,y$ can be joined by a rectifiable curve in $D$ of length at most $\Lambda_{1}|x-y|$.

\item \emph{(Lower bounds under non-degeneracy.)} If, for some $\eps>0$,
\begin{equation}\label{condi3}
j^{b}(x,y)\ge \eps |x-y|^{-d-\alpha}\quad\mbox{ for a.e. }x,y\in \R^{d},
\end{equation}
then $p^{b}_{D}(t,x,y)\ge C_{6}f^{0}_{D}(t,x,y)$ with $C_{6}=C_{6}(\eps)>0$. If, for some $\eps>0$,
\begin{equation}\label{condi4}
j^{b}(x,y)\ge \eps \big(|x-y|^{-d-\alpha}+|x-y|^{-d-\beta}\big)\quad\mbox{ for a.e. }x,y\in \R^{d},
\end{equation}
then $p^{b}_{D}(t,x,y)\ge C_{7}f^{1}_{D}(t,x,y)$ with $C_{7}=C_{7}(\eps)>0$.

\item \emph{(Large time.)} Let $D$ be bounded. There are $\lambda_{0}=\lambda_{0}(d,\alpha,\beta,\mathrm{diam}(D),A)>0$ and $C_{8}=C_{8}(\mathrm{diam}(D))>0$ such that, for every $b\in\mathcal{B}_{A}$,
\[
\lambda^{b,D}_{1}\ge\lambda_{0}
\quad\hbox{and}\quad
p^{b}_{D}(t,x,y)\le C_{8} \e^{-\lambda_{0}t}\delta_{D}(x)^{\alpha/2}\delta_{D}(y)^{\alpha/2}
\]
for $t\ge T$ and $x,y\in D$. If, moreover, $p^{b}_{D}>0$ on $(0,\infty)\times D\times D$, then this bound holds with $\lambda^{b,D}_{1}$ in place of $\lambda_{0}$ and a constant $C_{9}=C_{9}(D,b)>0$ in place of $C_{8}$.
Finally, suppose that one of the following conditions holds:
\begin{enumerate}[label=\upshape(\alph*), ref=\upshape(\alph*)]
\item $\mathrm{diam}(D)< R_1$, where $R_{1}$ is the constant in {\rm(ii)};
\item $D$ is connected and the path distance in $D$ is comparable to the Euclidean distance with characteristic $\Lambda_{1}\ge 1$;
\item condition \eqref{condi3} holds for some $\eps>0$.
\end{enumerate}
Then $p^{b}_{D}>0$ on $(0,\infty)\times D\times D$, and there is $C_{10}\ge 1$, depending also on $\mathrm{diam}(D)$, on $\Lambda_{1}$ under {\rm(b)} and on $\eps$ under {\rm(c)}, such that
\[
C_{10}^{-1} \e^{-\lambda^{b,D}_{1}t}\delta_{D}(x)^{\alpha/2}\delta_{D}(y)^{\alpha/2}\le p^{b}_{D}(t,x,y)
\le C_{10} \e^{-\lambda^{b,D}_{1}t}\delta_{D}(x)^{\alpha/2}\delta_{D}(y)^{\alpha/2}
\]
for $t\ge T$ and $x,y\in D$.
\end{enumerate}
\end{thm}

Theorem \ref{T:1.2} separates three regimes.
Part (i) gives the small-time upper bound. Only the large jumps of $b^{+}$ enter, through $M_{b^{+},\lambda}$.
Parts (ii) and (iii) give the small-time lower bound.
Under \eqref{condi2} alone, $j^{b}$ may vanish for large jumps, as in the truncated case above.
The lower bound by $f^{0}_{D}$ then holds only for $|x-y|<R_{1}$.
For larger distances it is replaced by a bound of the faster-decaying type satisfied by $\bar{p}_{1}(t,x,y)$.
Part (iv) treats large time on bounded sets. There the decay is exponential, with rate given by $\lambda^{b,D}_{1}$.

When the jumping kernel $j^{b}$ satisfies matching upper and lower bounds, parts (i), (iii) and (iv) combine into the following sharp estimates.

 \begin{corollary}
Suppose $d\ge 1$ and
$D$ is a $C^{1,1}$ open set in $\R^{d}$
with characteristics $(R_{0},\Lambda_{0})$. For every
$A,T\in (0,+\infty)$ and $\eps\in (0,1]$,
there is a constant $C_{11}=C_{11}(d,\allowbreak \alpha,\allowbreak \beta,\allowbreak R_{0},\allowbreak \Lambda_{0},\allowbreak A,\allowbreak T,\allowbreak \eps)\ge 1$ such that
for every $b\in \mathcal{B}_{A}$ satisfying
\[\eps |x-y|^{-d-\alpha}\le j^{b}(x,y)\le \eps^{-1}|x-y|^{-d-\alpha}\quad\mbox{ for a.e. }x,y\in \R^{d},\]
we have
\[C_{11}^{-1}f^{0}_{D}(t,x,y)\le p^{b}_{D}(t,x,y)\le C_{11} f^{0}_{D}(t,x,y)\quad\mbox{ for }t\in (0,T]\mbox{ and }x,y\in D.\]
If, in addition, $D$ is bounded, then there is a constant
$C_{12}=C_{12}(d,\alpha,\beta,R_{0},\Lambda_{0},\mathrm{diam}(D),A,T,\eps)\ge 1$
such that for every $b\in\mathcal{B}_{A}$ satisfying the above conditions,
\[C_{12}^{-1}\e^{-t \lambda^{b,D}_{1}}\delta_{D}(x)^{\alpha/2}\delta_{D}(y)^{\alpha/2}\le p^{b}_{D}(t,x,y)\le C_{12} \e^{-t \lambda^{b,D}_{1}}\delta_{D}(x)^{\alpha/2}\delta_{D}(y)^{\alpha/2}\]
for $t\ge T$ and $x,y\in D$.
\end{corollary}

\begin{corollary}
Suppose $d\ge 1$ and
$D$ is a $C^{1,1}$ open set in $\R^{d}$
with characteristics $(R_{0},\Lambda_{0})$. For every $A,T,\eps\in (0,{\infty})$, there is a constant $C_{13}=C_{13}(d,\allowbreak \alpha,\allowbreak \beta,\allowbreak R_{0},\allowbreak \Lambda_{0},\allowbreak A,\allowbreak T,\allowbreak \eps)\ge 1$ such that
for all $b\in\mathcal{B}_{A}$ satisfying \eqref{condi4},
\[C_{13}^{-1}f^{1}_{D}(t,x,y)\le p^{b}_{D}(t,x,y)\le C_{13} f^{1}_{D}(t,x,y)\quad\mbox{ for }t\in (0,T]\mbox{ and }x,y\in D.\]
If, in addition, $D$ is bounded, then there is a constant
$C_{14}=C_{14}(d,\alpha,\beta,R_{0},\Lambda_{0},\mathrm{diam}(D),A,T,\eps)\ge 1$
such that
for all $b\in\mathcal{B}_{A}$ satisfying \eqref{condi4},
\[C_{14}^{-1}\e^{-t \lambda^{b,D}_{1}}\delta_{D}(x)^{\alpha/2}\delta_{D}(y)^{\alpha/2}\le p^{b}_{D}(t,x,y)\le C_{14} \e^{-t \lambda^{b,D}_{1}}\delta_{D}(x)^{\alpha/2}\delta_{D}(y)^{\alpha/2}\]
for $t\ge T$ and $x,y\in D$.
\end{corollary}

Since $G^b_{D}(x, y) = \int_0^\infty p^b_{D}(t, x, y)\, {\mathrm{d}t}$, integrating the above heat kernel estimates in $t$ yields the following sharp two-sided estimates on the Green function $G^b_D(x, y)$ of $\L^b$ in $D$.
See the proof of \cite[Corollary 1.2]{CKS1} for the details about such integration.

\begin{corollary}
 Suppose $d\ge 1$
 and $D$ is a bounded $C^{1,1}$ open set in $\R^{d}$ with characteristics $(R_{0},\Lambda_{0})$.
 For every $A, \eps \in (0, \infty)$,
there exists a constant
$C_{15}=C_{15}(d,\alpha,\beta,R_{0},\Lambda_{0},\mathrm{diam}(D),A,\varepsilon)\geq 1 $
 such that
for all $b\in\mathcal{B}_{A}$ satisfying \eqref{condi3},
 \begin{equation} \nonumber
G^{b}_{D}(x,y)\stackrel{C_{15}}{\asymp}\begin{cases}
         \big(1\wedge\frac{\delta_{D}(x)^{\alpha/2}\delta_{D}(y)^{\alpha/2}}{|x-y|^{\alpha}}\big)|x-y|^{\alpha-d}
         \quad &\hbox{if } d>\alpha, \smallskip \\
         \log\big(1+\frac{\delta_{D}(x)^{\alpha/2}\delta_{D}(y)^{\alpha/2}}{|x-y|^{\alpha}}\big)
         \quad &\hbox{if } d=1=\alpha, \smallskip \\
         \big(\delta_{D}(x)\delta_{D}(y)\big)^{(\alpha-1)/2}\wedge \frac{\delta_{D}(x)^{\alpha/2}\delta_{D}(y)^{\alpha/2}}{|x-y|}
         &\hbox{if }   d=1<\alpha,
        \end{cases}
\end{equation}
\end{corollary}

\begin{remark}\label{R:dim}
Sections \ref{sec5}--\ref{sec7} treat the case $d>\alpha$. The condition $d>\alpha$ enters the proofs in two places.
First, for $d>\alpha$ the Green function $G_{D}$ of the killed rotationally symmetric $\alpha$-stable process is comparable to $|x-y|^{\alpha-d}$ times boundary factors; see \eqref{greenesti}. The estimates of Section \ref{sec5} on $\mathcal{S}^{b}_{x}G_{D}(x,y)$ and on the Duhamel series are built on this form.
Second, for $d>\alpha$ the Green function of the subprocess $X^{b,q_{0}}$ is infinite exactly on the diagonal, which is condition (A1) in the construction of the dual process in Section \ref{sec6}.
Both properties fail for $d=1\le\alpha$. This case is dealt with in Section \ref{sec9}, where we obtain Theorem \ref{T:1.2} from the case $d=2$ by viewing $X^{b}$ as the first coordinate of a suitable process on $\R^{2}$.
\end{remark}

\medskip
\noindent{\bf Novelty.} To the best of our knowledge, Theorem \ref{T:1.2} gives the first sharp two-sided Dirichlet heat kernel estimates for a
 space-dependent
non-symmetric non-local perturbation of the fractional Laplacian. Previous Dirichlet estimates concern symmetric processes \cite{CKS1, CKS3, CKSrel, CKSrot, BGR1, BGR2}, or non-symmetric perturbations by gradient operators \cite{CKS2, KS}, which are local.
 We also mention two further lines of work. Chen, Kim and Song \cite{CKS15} established stability of Dirichlet
heat kernel estimates for symmetric non-local operators under Feynman--Kac perturbations. Their perturbation acts
through a multiplicative Feynman--Kac weight and is different from the non-local perturbation $\mathcal{S}^{b}$ studied
here. Recently, Cho and Song \cite[Section 5.2]{CS25} obtained sharp Dirichlet heat kernel and Green function estimates for
non-symmetric strictly $\alpha$-stable L\'{e}vy processes.
Their non-symmetry arises from directionally non-symmetric jumps, unlike our space-dependent perturbation $\mathcal{S}^{b}$, and can lead to boundary exponents different from $\alpha/2$.
Our setting differs from all of these in four respects.
First, the perturbation $\mathcal{S}^{b}$ is itself a jump operator. It may be negative and may even remove all large jumps, as in the truncated example above, so $X^{b}$ is in general not comparable to the rotationally symmetric $\alpha$-stable process.
Second, the estimates are uniform in $b$: the small-time constants depend on $b$ only through an upper bound of $\|b\|_{\infty}$, and the upper bounds see $b$ only through its large jumps via $M_{b^{+},\lambda}$.
Third, $\mathcal{L}^{b}$ is not symmetric, so the boundary decay in the second space variable requires an argument using
a suitable dual process.
Fourth, the results cover every dimension $d\ge 1$, including $d=1\le\alpha$, through a lifting construction that appears to be new in this setting.

\medskip
\noindent{\bf Strategy.} Since $\mathcal{L}^{b}=\Delta^{\alpha/2}+\mathcal{S}^{b}$ is a lower order perturbation of $\Delta^{\alpha/2}$, one expects $p^{b}_D$ to be related to the Dirichlet heat kernel $p_{D}$ of the killed rotationally symmetric $\alpha$-stable process $X^{0,D}$ by
\begin{equation}
p^{b}_{D}(t, x,y)=p_{D}(t, x,y)+\int_0^t \int_{D} p^{b}_{D}(s, x,z)\mathcal{S}^{b}_{z} p_{D}(t-s, z,y){\mathrm{d}z} {\mathrm{d}s}\quad \mbox{ for }x,y\in D .
\nonumber
\end{equation}
Pointwise estimates on $\mathcal{S}^{b}_{z}p_{D}(t-s,z,y)$ are, however, difficult to obtain. Following \cite{CKS2}, we work instead with Green functions. On a bounded open set $D$, the Green function $G^{b}_{D}$ of $X^{b,D}$ satisfies Duhamel's formula
\begin{equation}
G^{b}_{D}(x,y)=G_{D}(x,y)+\int_{D}G^{b}_{D}(x,z)\mathcal{S}^{b}_{z}G_{D}(z,y){\mathrm{d}z}\quad\mbox{ for }
x,y\in D, x\not=y,
\nonumber
\end{equation}
where $G_{D}$ is the Green function of $X^{0,D}$. Iterating this formula and bounding $\mathcal{S}^{b}_{z}G_{D}$ (see (2) below) gives $G^{b}_{D}\asymp G_{D}$ on $C^{1,1}$ sets of small diameter.
   This comparison is the starting point for the boundary decay estimates. Together with the duality and chaining arguments, it leads to the small-time heat kernel estimates.

\medskip
\noindent{\bf Main difficulties.} The main difficulties stem from two facts: no explicit formula is available for $p^{b}_{D}$, and several standard tools fail for $\mathcal{L}^{b}$. This accounts for the length of the paper. We list the main obstacles and how they are overcome.
\begin{enumerate}[label=\upshape(\arabic*), ref=\upshape(\arabic*)]
\item \emph{Lack of symmetry.} Upper bounds in the second space variable cannot be obtained from those in the first variable. We construct a dual process of the subprocess $X^{b,q_{0}}$ with respect to an excessive reference measure (Section \ref{sec6}), and we show that the duality constant $M$ is uniform over $\mathcal{B}_{A}$.
\item \emph{The perturbation is non-local and of non-fixed sign.} Unlike a gradient perturbation, $\mathcal{S}^{b}$ acts on the Green function $G_{D}$ through all of $\R^{d}$, and $\mathcal{S}^{b}_{x}G_{D}(x,y)$ may be singular as $x$ approaches $\partial D$. Sharp bounds on $\mathcal{S}^{b}_{x}G_{D}(x,y)$, needed in the Duhamel argument above, require three different kernels according to whether $\alpha>2\beta$, $\alpha=2\beta$ or $\alpha<2\beta$ (Section \ref{sec5}). Together with 3G-type inequalities, they yield the convergence of the Duhamel series and the comparability $G^{b}_{D}\asymp G_{D}$ on small $C^{1,1}$ sets.

\item \emph{Large jumps may be absent.} Since $j^{b}$ may vanish for large jumps, $p^{b}$ is in general not comparable to the $\alpha$-stable heat kernel. The lower bounds therefore split into a short-range regime, handled near the boundary by two-point estimates, and a long-range regime, handled by chaining along curves (Section \ref{sec7}).

\item \emph{Low dimension.} For $d=1\le\alpha$, the methods of Sections \ref{sec5}--\ref{sec7} break down; see Remark \ref{R:dim}. We lift $X^{b}$ to a process on $\R^{2}$ whose first coordinate is $X^{b}$, and transfer the two-dimensional estimates back to dimension one (Section \ref{sec9}).

\item \emph{Large time.} The generator $\mathcal{L}^{b,D}$ is not self-adjoint, so the variational characterization of the principal eigenvalue is unavailable. We use the spectral theory of positive compact operators, and we give sufficient conditions under which a positive principal eigenfunction exists (Section \ref{sec8}).
 \end{enumerate}

\medskip
The paper is organized as follows. Section \ref{sec2} collects preliminaries, and Section \ref{sec3} establishes basic properties of killed processes. Sections \ref{sec5}--\ref{sec7} treat the case $d>\alpha$: Green function estimates on small sets (Section \ref{sec5}), duality (Section \ref{sec6}) and small-time heat kernel estimates (Section \ref{sec7}). Section \ref{sec9} treats the case $d=1\le\alpha$. Section \ref{sec8} treats large time and completes the proof of Theorem \ref{T:1.2}.

\section{Preliminaries}\label{sec2}

This section collects the known estimates and geometric facts used throughout the paper.
We first recall estimates on the heat kernel $p^b(t, x, y)$ of
$\L^b$ from \cite[Theorems 1.3 and 1.5]{CW}.
Throughout Sections \ref{sec2} and \ref{sec3}, $d\geq 1$ is arbitrary.

\begin{proposition}\label{P:2.1}
Suppose $A,\lambda,\varepsilon\in (0,+\infty)$. The following hold.
\begin{enumerate}[label=\upshape(\roman*), ref=\upshape(\roman*)]
\item There are positive constants
$C_{i}=C_{i}(d,\allowbreak \alpha,\allowbreak \beta,\allowbreak A,\allowbreak \lambda)$, $i=16,\cdots,20$ such that for
all $b\in\mathcal{B}_{A}$, with $a:=M_{b^{+},\lambda}$,
\begin{equation}\label{stdensity}
C_{16}^{-1}\bar{p}_{1}(t,C_{17}x,C_{17}y)\le p^{b}(t,x,y)\le
C_{16}p_{a}(t,x,y)
\end{equation}
for $t\in (0,1]$ and $x,y\in\R^{d}$, and
\begin{equation}
C_{18}^{-1}\e^{-C_{19}t}\bar{p}_{1}(t,C_{20}x,C_{20}y)\le p^{b}(t,x,y)\le
C_{18}\e^{C_{19}t}p_{a}(t,x,y)\label{ltdensity}
\end{equation}
 for $t\in (0,{\infty})$ and $x,y\in \R^{d}$.

  \item There are
positive constants $C_{i}=C_{i}(d,\allowbreak \alpha,\allowbreak \beta,\allowbreak A,\allowbreak \lambda,\allowbreak \eps)$, $i=21,22,23$ such that for
all $b\in\mathcal{B}_{A}$ satisfying \eqref{condi3}, with $a:=m_{b^{+},\lambda}$,
 \begin{equation}\nonumber
p^{b}(t,x,y)\ge C_{21}^{-1}p_{a}(t,x,y)\quad\mbox{ for }t\in (0,1]\mbox{ and }x,y\in\R^{d},
\end{equation}
and
\begin{equation}\label{2.4}
p^{b}(t,x,y)\ge C_{22}\e^{-C_{23}t}p_{a}(t,x,y)\quad\mbox{ for }t\in (0,\infty)\mbox{ and }x,y\in \R^{d}.
\end{equation}
Here
$m_{b^{+},\lambda}:=\einf_{x,z\in\R^{d},|z|>\lambda}b^{+}(x,z)$.
\item There
are positive constants $C_{i}=C_{i}(d,\allowbreak \alpha,\allowbreak \beta,\allowbreak A,\allowbreak \eps)$, $i=24,25,26$,
such that for
all $b\in\mathcal{B}_{A}$ satisfying \eqref{condi5},
we have
\begin{equation}\nonumber
 p^{b}(t,x,y)\le C_{24}p(t,x,y)\quad\mbox{ for }t\in (0,1]\mbox{ and }x,y\in\R^{d},
\end{equation}
and
\begin{equation}\label{2.6}
p^{b}(t,x,y)\le C_{25}\e^{C_{26}t}p(t,x,y)\quad\mbox{ for }t\in (0,{\infty})\mbox{ and }x,y\in \R^{d}.
\end{equation}
\end{enumerate}
 \end{proposition}

\begin{proof}
For $t\in (0,1]$, the estimates in (i) and (ii) are \cite[Theorem 1.3]{CW}.
(A gap in the proof of (1.16) in \cite{CW} is corrected in \cite{CW2},
and all main results of \cite{CW} remain valid; see \cite[p.~760]{CW2}.)
The estimates for $t\in (0,\infty)$ follow from those for $t\in (0,1]$ by the Chapman--Kolmogorov equation, since $\{p_{a}(t,\cdot,\cdot)\}_{t>0}$ satisfies the Chapman--Kolmogorov equation and
$\int_{\R^{d}}\bar{p}_{1}(t,cx,cz)\bar{p}_{1}(s,cz,cy)\,\mathrm{d}z=c^{-d}\bar{p}_{1}(t+s,cx,cy)$ for every $c>0$.

For (iii), let $b\in\mathcal{B}_{A}$ satisfy \eqref{condi5}, and define
\[
b_1(x,z):=b(x,z)\vee\Big(-\frac{\mathcal{A}(d,-\alpha)}{2\mathcal{A}(d,-\beta)}|z|^{\beta-\alpha}\Big).
\]
Then $b_1(x,z)=b_1(x,-z)$ and  $\|b_1\|_{\infty} \leq \|b\|_{\infty}\le A$.
  Moreover,
$j^{b_1}=j^{b}\vee \big(\frac{1}{2}\mathcal{A}(d,-\alpha)|x-y|^{-d-\alpha}\big)$, and hence
\[
\frac{1}{2}\mathcal{A}(d,-\alpha)|x-y|^{-d-\alpha}\le j^{b_1}(x,y)\le \big(\eps^{-1}\vee
\frac{1}{2}
\mathcal{A}(d,-\alpha)\big)|x-y|^{-d-\alpha}\quad\mbox{for a.e. }x,y\in\R^{d}.
\]
By \cite[Theorem 1.5]{CW}, $p^{b_1}(t,x,y)\le c_{1}p(t,x,y)$ for $t\in (0,1]$ and $x,y\in\R^{d}$, with $c_{1}=c_{1}(d,\alpha,\beta,A,\eps)>0$.
Note that  $j^{b_1}(x,y)=j^{b}(x,y)$ for $|x-y|<\eta(A)  :=\big(\frac{1}{2A} \frac{\mathcal{A}(d,-\alpha)}{\mathcal{A}(d,-\beta)}\big)^{1/(\alpha-\beta)}$.
Hence
\[
\kappa:=\sup_{x\in\R^{d}}\int_{\R^{d}}\big(j^{b_1}(x,y)-j^{b}(x,y)\big)\mathrm{d}y\le \frac{1}{2}\mathcal{A}(d,-\alpha)\int_{|z|\ge\eta(A)}|z|^{-d-\alpha}\mathrm{d}z<\infty .
\]
By Meyer's construction, the process $X^{b_1}$ started at $x$ can be  obtained from $X^{b}$
by adding jumps with intensity $
j^{b_{1}}
-j^{b}$  (see \cite{Meyer} and \cite[Remark 3.4]{BBCK}).
The identification in distribution of this construction is justified by
 the uniqueness of the martingale problem for $(\mathcal{L}^{b_1},\mathcal{S}(\R^{d}))$ in \cite[Theorem 1.3]{CW}.
On the event that no jump is added before time $t$, which has conditional probability at least $\e^{-\kappa t}$, the two processes coincide on $[0,t]$. Therefore $p^{b}(t,x,y)\le \e^{\kappa t}p^{b_1}(t,x,y)\le c_{1}\e^{\kappa}p(t,x,y)$ for $t\in (0,1]$; see \cite[Lemma 3.6]{BBCK}.
The estimate for $t\in (0,\infty)$ follows by the Chapman--Kolmogorov equation.
\end{proof}

We also need the following known
geometric properties of
 a $C^{1,1}$ open set $D$ with $C^{1,1}$ characteristics
 $(R_0, \Lambda_0)$.

\begin{proposition}\label{P:2.2}
The following hold for a $C^{1,1}$ open set $D$ with $C^{1,1}$ characteristics
 $(R_0, \Lambda_0)$.
\begin{enumerate}[label=\upshape(\roman*), ref=\upshape(\roman*)]
 \item {(Outer and inner ball property) There is
 a positive constant $r_{0}=r_{0}(R_{0},\Lambda_{0})$
 such that for every $Q\in\partial D$
and $0<r<r_{0}$, there are balls $B(x',r)\subset D,\ B(x'',r)\subset
D^{c}$ tangent at $Q$. In this case, we say $D$ is a $C^{1,1}$ open set
at scale $r_{0}$.}

\item {There exists a constant $L=L(d,R_{0},\Lambda_{0})>0$ such that
for every $
Q\in\partial D
$ and $0<r\le R_{0}$, there is a $C^{1,1}$ domain $V=V(Q,r)$ with characteristics $(rR_{0}/L,\Lambda_{0}L/r)$ such that $D\cap
 B(Q,r) \subset V \subset D\cap
 B(Q,2r) $.
}

\item {There exists a constant $\kappa=\kappa(\Lambda_{0})\in (0,1/2)$ such that for every $r\in (0,R_{0})$ and $Q\in \partial D$, there
is a point $\xi$ in $D\cap B(Q,r)$, denoted by $\xi_{r}(Q)$, such that
$B(\xi,\kappa r)\subset D\cap B(Q,r)$. The pair $(R_{0},\kappa)$ is called the
$\kappa$-fat characteristics of $D$. }
\end{enumerate}
\end{proposition}
For (i), (ii) and (iii) we refer to
 \cite[Lemma 2.2]{AKSZ}, \cite[Lemma 2.2]{Song} and \cite[p. 47, property P2]{Bogdan1},
respectively.

\begin{proposition}\label{prop1}
Suppose $D$ is a $C^{1,1}$ open set with characteristics
$(R_{0},\Lambda_{0})$. Then there is a constant $\theta_{0}=\theta_{0}(\Lambda_{0})\in
(0,1)$ such that for every $x\in D$ and $r\in (0,R_{0})$, there
exists a ball $B(\xi,\theta_{0} r)\subset D\cap B(x,r)$.
\end{proposition}

\begin{proof} Let $\kappa=\kappa(\Lambda_{0})\in (0,1/2)$ be the constant in Proposition \ref{P:2.2}(iii), and set $\theta_{0}:=\kappa/2$.
Fix $x\in D$ and $r\in (0,R_{0})$.
If $\delta_{D}(x)\ge r/2$, then $B(x,\theta_{0}r)\subset B(x,r/2)\subset D\cap B(x,r)$.
If $\delta_{D}(x)<r/2$, choose $Q\in\partial D$ with $|x-Q|=\delta_{D}(x)$. Since $r/2<R_{0}$, Proposition \ref{P:2.2}(iii) gives a point $\xi$ with $B(\xi,\kappa r/2)\subset D\cap B(Q,r/2)$.
Since $B(Q,r/2)\subset B(x,r/2+\delta_{D}(x))\subset B(x,r)$, we get $B(\xi,\theta_{0}r)\subset D\cap B(x,r)$.
\end{proof}

In dimension one, Proposition \ref{P:2.2}(i), (iii) and Proposition \ref{prop1} follow directly from Definition \ref{C11set}, with $r_{0}$ comparable to $R_{0}$, and universal constants $\kappa$ and $\theta_{0}$. In Proposition \ref{P:2.2}(ii), $V(Q,r)$ can be chosen to be an open interval whose $C^{1,1}$ characteristics are comparable to $r$.

\section{Process killed upon leaving an open set}\label{sec3}

This section establishes the basic properties of the killed process: joint continuity of its transition density,
lower bounds in small balls, and strict positivity. At the end of the section these properties are transferred to the process killed at an independent exponential time. These results are used repeatedly in Sections \ref{sec5}--\ref{sec8}.
Throughout,
$A\in (0,+\infty)$, $b\in \mathcal{B}_{A}$
and $X^{b}$ is the corresponding Feller process.
Let  $D\subset \R^{d}$ be an open subset, and let $X^{b,D}$ be the subprocess of $X^b$ killed upon leaving $D$.
Recall that $\partial $ is a cemetery added to $D$. Let $D_{\partial}:=D\cup\{\partial\}$.
Define for every $x,y\in D$,
\[
N^{D}(x,{\mathrm{d}y}):=j^{b}(x,y){\mathrm{d}y},\quad N^{D}(x,\{\partial\})=\int_{D^{c}}j^{b}(x,y){\mathrm{d}y}.
\]
It follows from \eqref{levysystemforxb} that
$(N^{D}, t\wedge \tau^{b}_{D})$
is a
L\'{e}vy system for $X^{b,D}$, \textit{i.e.}, for any $x\in D$, any
non-negative measurable function $f$ on $[0,\infty)\times D\times
D_{\partial}$ vanishing on
$\{(s,y,y): s\ge 0,\ y\in\R^{d}\}$
and stopping time $T$ (with respect to the
filtration of $X^{b,D}$),
\begin{equation}\label{levysysforxbd}
\E_x\Big[\sum_{s\le T}f(s,X^{b,D}_{s-},X^{b,D}_{s})\Big]
=\E_x\Big[\int_{0}^{T}\int_{D_{\partial}}f(s,X^{b,D}_{s},y)N^{D}(X^{b,D}_{s},{\mathrm{d}y})
{\mathrm{d}s} \Big].
\end{equation}

 For every $A\in (0,{\infty})$, define
\begin{equation}\label{defiforetaA}
\eta(A)
:=\big(\frac{1}{2A} \frac{\mathcal{A}(d,-\alpha)}{\mathcal{A}(d,-\beta)}\big)^{1/(\alpha-\beta)}.
\end{equation}
Since $b\in\mathcal{B}_{A}$,
we have  for every
$x\in\R^{d}$ and a.e. $y\in\R^{d}$
\begin{eqnarray}
j^{b}(x,y)&=&\frac{\mathcal{A}(d,-\alpha)}{|x-y|^{d+\alpha}}+\mathcal{A}(d,-\beta)\frac{b(x,y-x)}{|x-y|^{d+\beta}}\1_{\{|x-y|<\eta(A)\}}+\mathcal{A}(d,-\beta)\frac{b(x,y-x)}{|x-y|^{d+\beta}}\1_{\{|x-y|\ge\eta(A)\}}\nonumber\\
&\ge&\frac{\mathcal{A}(d,-\alpha)}{|x-y|^{d+\alpha}}-\mathcal{A}(d,-\beta)\frac{A}{|x-y|^{d+\beta}}\1_{\{|x-y|<\eta(A)\}}
-\mathcal{A}(d,-\alpha)\frac{1}{|x-y|^{d+\alpha}}\1_{\{|x-y|\ge\eta(A)\}}\nonumber\\
&=&\frac{\mathcal{A}(d,-\alpha)}{|x-y|^{d+\alpha}}\1_{\{|x-y|<\eta(A)\}}-\mathcal{A}(d,-\beta)\frac{A}{|x-y|^{d+\beta}}\1_{\{|x-y|<\eta(A)\}}\nonumber\\
&=&\frac{\mathcal{A}(d,-\alpha)}{|x-y|^{d+\alpha}}\1_{\{|x-y|<\eta(A)\}}\big(1-\frac{\mathcal{A}(d,-\beta)}{\mathcal{A}(d,-\alpha)}\,A|x-y|^{\alpha-\beta}\big)\nonumber\\
&\ge&\frac{1}{2}\bar{j}_{\eta(A)}(x,y),
\nonumber
\end{eqnarray}
where
$\bar{j}_{\eta(A)}(x,y):=\mathcal{A}(d,-\alpha)|x-y|^{-d-\alpha}\1_{\{|x-y|<\eta(A)\}}$.
In
other words, we have for every $x\in\R^{d}$,
\begin{equation} \label{jb&j}
j^{b}(x,y)\ge \frac{1}{2}j(x,y)\quad \mbox{a.e. on
}\{y\in\R^{d}:\ |x-y|<\eta(A)\}.
\end{equation}

\begin{lemma}\label{lem3.1}
There is a constant $C_{27}=C_{27}(d,\allowbreak \alpha,\allowbreak \beta,\allowbreak A)>0$ such that for every $\delta>0$ and $s>0$,
\begin{equation}\label{lem3.1.1}
\sup_{x\in\R^{d}}\P_x(\tau^{b}_{B(x,\delta)}\le s)\le C_{27}\,\e^{C_{27}s}\,s\big(\delta^{-\alpha}+\delta^{-\beta}\big).
\end{equation}
In particular, $\lim_{s\downarrow 0}\sup_{x\in\R^{d}}\P_x(\tau^{b}_{B(x,\delta)}\le s)=0$ for every $\delta>0$.
\end{lemma}

\begin{proof} For every $x\in\R^{d}$, we have
\begin{eqnarray}
&&\P_x(\tau^{b}_{B(x,\delta)}\le s)\nonumber\\
&\le&\P_x(\tau^{b}_{B(x,\delta)}\le s,X^{b}_{s}\in B(x,\delta/2))+\P_x(X^{b}_{s}\in B(x,\delta/2)^{c})\nonumber\\
&\le&\E_x\big[\P_{X^{b}_{\tau^{b}_{B(x,\delta)}}}\big(|X^{b}_{s-\tau^{b}_{B(x,\delta)}}-X^{b}_{0}|\ge \delta/2\big);\,\tau^{b}_{B(x,\delta)}\le s\big]
+\P_x\big(|X^{b}_{s}-X^{b}_{0}|\ge \delta/2\big)\nonumber\\
&\le& 2\sup_{t\in[0,s]}\sup_{z\in\R^{d}}\P_z\big(|X^{b}_{t}-X^{b}_{0}|\ge \delta/2\big).\nonumber
\end{eqnarray}
By \eqref{ltdensity},
\begin{eqnarray}
&&\sup_{t\in[0,s]}\sup_{z\in\R^{d}}\P_z\big(|X^{b}_{t}-X^{b}_{0}|\ge \delta/2\big)\nonumber\\
&=&\sup_{t\in[0,s]}\sup_{z\in\R^{d}}\P_z\big(|X^{b}_{t}-z|\ge \delta/2\big)\nonumber\\
&=&\sup_{t\in[0,s]}\sup_{z\in\R^{d}}\int_{|y-z|\ge \delta/2}p^{b}(t,z,y){\mathrm{d}y}\nonumber\\
&\le&\sup_{t\in[0,s]}\sup_{z\in\R^{d}}c_{1}\e^{c_{2}t}\int_{|y-z|\ge \delta/2}t^{-d/\alpha}\wedge t^{-d/\beta}\wedge\big(\frac{t}{|z-y|^{d+\alpha}}+\frac{t}{|z-y|^{d+\beta}}\big){\mathrm{d}y}\nonumber\\
&\le&c_{3}\e^{c_{2}s}s\int_{\delta/2}^{\infty}(r^{-\alpha-1}+r^{-\beta-1}){\mathrm{d}r}\nonumber
\end{eqnarray}
for some constants $c_{i}=c_{i}(d,\allowbreak \alpha,\allowbreak \beta,\allowbreak A)>0$, $i=1,2,3$.
Since $\int_{\delta/2}^{\infty}(r^{-\alpha-1}+r^{-\beta-1}){\mathrm{d}r}\le 2^{\alpha}\alpha^{-1}\delta^{-\alpha}+2^{\beta}\beta^{-1}\delta^{-\beta}$, this proves \eqref{lem3.1.1}. \end{proof}

\begin{thm}\label{them3.4}
Let $D$ be an open set in $\R^{d}$. The density function $p^{b}_{D}(t,x,y)$ is jointly continuous in $(0, \infty)\times D\times D$ and satisfies
\begin{equation}\label{e:3.5}
p^{b}_{D}(t+s,x,y)=\int_{D}p^{b}_{D}(t,x,z)p^{b}_{D}(s,z,y){\mathrm{d}z},\quad
 \mbox{ for } t,s>0\mbox{ and }x,y\in D.
\end{equation}
\end{thm}

\begin{proof}
Define $k^{b}_{D}(t,x,y):=\E_x\big[p^{b}(t-\tau^{b}_{D},X^{b}_{\tau^{b}_{D}},y);\,\tau^{b}_{D}<t\big]$ for $(t, x, y)\in (0, \infty)\times \R^d\times \R^d$.   By \eqref{defiforpbd}, it suffices to show that $k^{b}_{D}(t,x,y) $ is jointly continuous on $(0, \infty)\times D\times D$. By \eqref{ltdensity}, there are positive constants $c_{i}=c_{i}(d,\allowbreak \alpha,\allowbreak \beta,\allowbreak A)$, $i=1,\cdots,4$ such that
\begin{eqnarray}
p^{b}(t,x,y)&\le&c_{1}\e^{c_{2}t}p_{\|b\|_{\infty}}(t,x,y)\nonumber\\
&\le&c_{3}\e^{c_{2}t}\big[t^{-d/\alpha}\wedge(\|b\|_{\infty}t)^{-d/\beta}\wedge\big(\frac{t}{|x-y|^{d+\alpha}}+\frac{\|b\|_{\infty}t}{|x-y|^{d+\beta}}\big)\big]\nonumber\\
&\le&c_{4}\e^{c_{2}t}\big[t^{-d/\alpha}\wedge
t^{-d/\beta}\wedge\big(\frac{t}{|x-y|^{d+\alpha}}+\frac{t}{|x-y|^{d+\beta}}\big)\big].\nonumber
\end{eqnarray}
Thus for any $t_{0}>0$ and $\delta>0$, we have
\begin{eqnarray}
&&\sup_{t\le t_{0}}\sup_{|x-y|\ge \delta}p^{b}(t,x,y)\nonumber\\
&\le&c_{4}\e^{c_{2}t_{0}}\sup_{t\le t_{0},|x-y|\ge\delta}\big[t^{-d/\alpha}\wedge t^{-d/\beta}\wedge\big(\frac{t}{|x-y|^{d+\alpha}}+\frac{t}{|x-y|^{d+\beta}}\big)\big]\nonumber\\
&\le&c_{4}\e^{c_{2}t_{0}}\big(\frac{t_{0}}{\delta^{d+\alpha}}+\frac{t_{0}}
{\delta^{d+\beta}}\big)
=:
c_{5}(d,\allowbreak \alpha,\allowbreak \beta,\allowbreak A,\allowbreak t_{0},\allowbreak \delta)< \infty.\label{them3.4-eq1}
\end{eqnarray}
Fix $t_{0}>0$, $y_{0}\in D$ and $\delta>0$ with $B(y_{0},2\delta)\subset D$, and let $y\in B(y_{0},\delta)$. Since $X^{b}_{\tau^{b}_{D}}\notin D$, we have $|X^{b}_{\tau^{b}_{D}}-y|\ge\delta$. Define  $q(s,z,y):=p^{b}(s,z,y)$ for $s>0$ and $q(s,z,y):=0$ for $s\le 0$.  Since $p^{b}(s,z,y)\to 0$ uniformly on $\{|z-y|\ge\delta\}$ as $s\downarrow 0$, we have by  \eqref{ltdensity} and \eqref{them3.4-eq1} that
$q$ is bounded by $c_{5}=c_{5}(d,\alpha,\beta,A,2t_{0},\delta)$ and continuous on $(-\infty,2t_{0}]\times\{(z,y):|z-y|\ge\delta\}$.
 Hence $0\le k^{b}_{D}(t,x,y) =\mathbb{E}_{x}[q(t-\tau^{b}_{D},X^{b}_{\tau^{b}_{D}},y)] \le c_{5}$ for $t\le 2t_{0}$, $x\in\R^{d}$ and $y\in B(y_{0},\delta)$, and, by dominated convergence, $(t,y)\mapsto k^{b}_{D}(t, x,y)$ is continuous for each $x\in\R^{d}$.

Let $s\in(0,t_{0}/2)$ and $t\in(t_{0}/2,2t_{0}]$. By the Markov property of $X^b$ at time $s$,
\begin{align*}
k^{b}_{D}(t,x,y)=&\int_{\R^{d}}p^{b}(s,x,z)k^{b}_{D}(t-s,z,y)\,\mathrm{d}z\\
&+\E_x\big[q(t-\tau^{b}_{D},X^{b}_{\tau^{b}_{D}},y);\tau^{b}_{D}\le s\big]-\E_x\big[k^{b}_{D}(t-s,X^{b}_{s},y);\tau^{b}_{D}\le s\big].
\end{align*}
The integral is jointly continuous in $(t,x,y)$. Indeed, $x\mapsto p^{b}(s,x,\cdot)$ is continuous in $L^{1}(\R^{d})$ by Scheff\'e's lemma, since $p^{b}(s,x,\cdot)$ is a probability density, and $k^{b}_{D}(t-s,z,y)$ is bounded by $c_{5}$ and continuous in $(t,y)$. The two expectations are bounded by $2c_{5}\P_x(\tau^{b}_{D}\le s)\le 2c_{5}\P_x(\tau^{b}_{B(x,\varrho)}\le s)$ for $x$ with $\delta_{D}(x)\ge\varrho>0$. By Lemma \ref{lem3.1}, this tends to $0$ as $s\downarrow 0$, uniformly in such $x$. Thus $k^{b}_{D}$ is a locally uniform limit of continuous functions on $(t_{0}/2,2t_{0}]\times\{x\in D:\delta_{D}(x)\ge\varrho\}\times B(y_{0},\delta)$, hence continuous there. Since $t_{0}$, $y_{0}$, $\delta$ and $\varrho$ are arbitrary, $k^{b}_{D}$, and therefore $p^{b}_{D}$, is jointly continuous on $(0,\infty)\times D\times D$.

Finally, we prove \eqref{e:3.5}. Fix $t,s>0$ and $x\in D$. By the Markov property of $X^{b,D}$,
\begin{equation}\label{e:3.6}
 p^{b}_{D} (t+s,x,y)=\int_{D }p^{b}_{D} (t,x,z)p^{b}_{D} (s,z,y) {\mathrm{d} z  }, \quad
 \mbox{ for a.e. } y\in D.
\end{equation}
By what we have just proved, $p^{b}_{D}(t+s,x,y)$ and $p^{b}_{D}(s,z,y)$ are continuous in $y\in D$ for each $z\in D$.
Since $p^{b}_{D}(s,z,y)\le\sup_{z,w}p^{b}(s,z,w)<\infty$ by \eqref{ltdensity}, and
$\int_{D} p^{b}_{D}(t,x,z)\,\mathrm{d}z\le 1$, the dominated convergence theorem shows that the right-hand side of \eqref{e:3.6}
is also continuous in $y\in D$. Consequently, \eqref{e:3.6} holds for every $y\in D$.
\end{proof}

\begin{lemma}\label{prop3.5'}
 Let $a_{1},\kappa_{1}\in(0,1)$, $A>0$ and $R\in (0,1/2]$. There are constants $l=l(d,\allowbreak \alpha,\allowbreak \beta,\allowbreak a_{1},\allowbreak \kappa_{1},\allowbreak R,\allowbreak A)\in (0,1)$ and $C_{28}=C_{28}(d,\allowbreak \alpha,\allowbreak \beta,\allowbreak a_{1},\allowbreak \kappa_{1},\allowbreak R,\allowbreak A)>0$ such that
for all $b\in\mathcal{B}_{A}$,
$x_{0}\in\R^{d}$ and $r\in (0,R]$,
 \begin{equation}\label{e:3.7}
p^{b}_{B(x_{0},r)}(t,x,y)\ge C_{28} r^{-d}
\quad \hbox{for } (t,x,y)\in[\kappa_{1} l r^{\alpha},lr^{\alpha}]\times B(x_{0},a_{1}r)\times B(x_{0},a_{1}r).
\end{equation}
Moreover, let $\eps>0$. Then the same holds for every $R>0$ and all $b\in \mathcal{B}_{A}$ satisfying \eqref{condi3}, with $l$ and $C_{28}$ depending also on $\eps$.
\end{lemma}

\begin{proof} Fix $x_{0}\in\R^{d}$. We use $B_{r}$ to denote $B(x_{0},r)$. By \eqref{ltdensity},
\begin{eqnarray}
&&p^{b}_{B_{r}}(t,x,y)\nonumber\\
&=&p^{b}(t,x,y)-\E_x\big[p^{b}(t-\tau^{b}_{B_{r}},X^{b}_{\tau^{b}_{B_{r}}},y);\,\tau^{b}_{B_{r}}<t\big]\nonumber\\
&\ge&c_{1}\e^{-c_{2}t}\bar{p}_{1}(t,x,y)-\E_x\Big[c_{3}\e^{c_{4}(t-\tau^{b}_{B_{r}})}p_{\|b\|_{\infty}}(t-\tau^{b}_{B_{r}},X^{b}_{\tau^{b}_{B_{r}}},y);\,\tau^{b}_{B_{r}}<t\Big].\nonumber
\end{eqnarray}
For every $0<r\le R\le 1/2$ and $x,y\in B(x_{0},a_{1}r)$,
we have $|x-y|\le 2a_{1}r<1$, and thus for any $t\in
[\kappa_{1} l r^{\alpha},lr^{\alpha}]\subset (0,1]$,
\begin{eqnarray}
\bar{p}_{1}(t,x,y)&\ge&c_{5}(d,\allowbreak \alpha)\big(t^{-d/\alpha}\wedge\frac{t}{|x-y|^{d+\alpha}}\big)
\, \geq \,
c_{6}(d,\allowbreak \alpha,\allowbreak a_{1})t^{-d/\alpha}\big(1\wedge \frac{t^{1/\alpha}}{r}\big)^{d+\alpha}\nonumber\\
&\ge&c_{6}(lr^{\alpha})^{-d/\alpha}\big(1\wedge \frac{(\kappa_{1} l)^{1/\alpha}r}{r}\big)^{d+\alpha}
\, = \, c_{6}\kappa_{1}^{1+d/\alpha}lr^{-d}.\label{prop3.5-eq1}
\end{eqnarray}
On the other hand, since $|X^{b}_{\tau^{b}_{B_{r}}}-y|\ge (1-a_{1})r$ for every $y\in B(x_{0},a_{1}r)$, we have
\begin{eqnarray}
&&p_{\|b\|_{\infty}}(t-\tau^{b}_{B_{r}},X^{b}_{\tau^{b}_{B_{r}}},y)\nonumber\\
&\le&c_{7}(d,\allowbreak \alpha,\allowbreak \beta)\big[(t-\tau^{b}_{B_{r}})^{-d/\alpha}\wedge
(\|b\|_{\infty}(t-\tau^{b}_{B_{r}}))^{-d/\beta}\wedge\big(\frac{t}{|X^{b}_{\tau^{b}_{B_{r}}}-y|^{d+\alpha}}+\frac{\|b\|_{\infty}t}
{|X^{b}_{\tau^{b}_{B_{r}}}-y|^{d+\beta}}\big)\big]\nonumber\\
&\le&c_{8}(d,\allowbreak \alpha,\allowbreak \beta,\allowbreak A,\allowbreak a_{1})\big(\frac{t}{r^{d+\alpha}}+\frac{t}{r^{d+\beta}}\big)\nonumber\\
&\le&c_{9}(d,\allowbreak \alpha,\allowbreak \beta,\allowbreak A,\allowbreak a_{1},\allowbreak R)\frac{t}{r^{d+\alpha}}.\label{prop3.5-eq2}
\end{eqnarray}
Enlarging $c_{4}$ if necessary, we may assume $c_{4}\ge C_{27}$.
Lemma \ref{lem3.1} gives, for every $0<r\le R<+\infty$ and $x\in B(x_{0},a_{1}r)$,
\begin{equation}
\P_x(\tau^{b}_{B_{r}}<t)\le \P_x(\tau^{b}_{B(x,(1-a_{1})r)}\le t)\le c_{10}(d,\allowbreak \alpha,\allowbreak \beta,\allowbreak A,\allowbreak a_{1},\allowbreak R)
\,\e^{c_{4}t}\,\frac{t}{r^{\alpha}}.\label{prop3.5-eq4}
\end{equation}
Thus by \eqref{prop3.5-eq2} and \eqref{prop3.5-eq4}, we get for every $t\in [\kappa_{1} l r^{\alpha},l r^{\alpha}]$ and $x,y\in B(x_{0},a_{1}r)$,
\begin{eqnarray}
&&\E_x\Big[c_{3}\e^{c_{4}(t-\tau^{b}_{B_{r}})}p_{\|b\|_{\infty}}(t-\tau^{b}_{B_{r}},X^{b}_{\tau^{b}_{B_{r}}},y);\,\tau^{b}_{B_{r}}<t\Big]\nonumber\\
&\le& c_{13}(d,\allowbreak \alpha,\allowbreak \beta,\allowbreak A,\allowbreak a_{1},\allowbreak R)\e^{2c_{4}t}\frac{t^{2}}{r^{d+2\alpha}}
 \, \le \, c_{13}\e^{2c_{4}lR^{\alpha}}l^{2}r^{-d}.\label{prop3.5-eq5}
\end{eqnarray}
It follows by \eqref{prop3.5-eq1} and \eqref{prop3.5-eq5} that for every $0<r\le R\le 1/2$ and $(t,x,y)\in [\kappa_{1}l r^{\alpha},l r^{\alpha}]\times B(x_{0},a_{1}r)\times B(x_{0},a_{1}r)$,
\[p^{b}_{B_{r}}(t,x,y)\ge l
\e^{-c_{2}lR^{\alpha}}(c_{1}c_{6}\kappa_{1}^{1+d/\alpha}-c_{13}l
\e^{(2c_{4}+c_{2})lR^{\alpha}})r^{-d}.\] The first assertion of Lemma
\ref{prop3.5'} follows by setting
$l=l(d,\allowbreak \alpha,\allowbreak \beta,\allowbreak a_{1},\allowbreak \kappa_{1},\allowbreak R,\allowbreak A)$ sufficiently small such
that $(c_{1}c_{6}\kappa_{1}^{1+d/\alpha}-c_{13}l
\e^{(2c_{4}+c_{2})lR^{\alpha}})>0$.
On the other hand, by \eqref{2.4}, we have for all $b\in\mathcal{B}_{A}$ satisfying \eqref{condi3}, $t\in (0,+\infty)$, $0<r\le R<+\infty$ and $x,y\in B_{r}$,
\[
p^{b}_{B_{r}}(t,x,y)\ge c_{14}\e^{-c_{15}t}p(t,x,y)-\E_x\Big[c_{16}\e^{c_{16}(t-\tau^{b}_{B_{r}})}p_{\|b\|_{\infty}}(t-\tau^{b}_{B_{r}},X^{b}_{\tau^{b}_{B_{r}}},y);\,\tau^{b}_{B_{r}}<t\Big].
\]
Then a similar argument (with the first inequality of \eqref{prop3.5-eq1} replaced by \eqref{p0}) gives the second assertion.
\end{proof}

\begin{proposition}\label{prop3.5}
Let $a_{1}\in (0,1)$, $a_{3}>a_{2}>0$, $A>0$ and $R\in (0,1/2]$. There is a positive constant $C_{29}=C_{29}(d,\allowbreak \alpha,\allowbreak \beta,\allowbreak a_{1},\allowbreak a_{2},\allowbreak a_{3},\allowbreak R,\allowbreak A)$ such that
for all $b\in\mathcal{B}_{A}$, $x_{0}\in\R^{d}$ and $r\in (0,R]$,
\begin{equation*}
p^{b}_{B(x_{0},r)}(t,x,y)\ge C_{29}r^{-d} \quad \hbox{for }
 (t,x,y)\in [a_{2}r^{\alpha},a_{3}r^{\alpha}]\times B(x_{0},a_{1}r)\times B(x_{0},a_{1}r).
\end{equation*}
Moreover, let $\eps>0$. Then the same holds for every $R>0$ and all $b\in\mathcal{B}_{A}$ satisfying \eqref{condi3}, with $C_{29}$ depending also on $\eps$.
\end{proposition}

\begin{proof} Let $\kappa_{1}:=a_{2}/(1+a_{3})\in (0,1)$ and let
$l$ and $C_{28}$ be the constants given in \eqref{e:3.7} that depend on $d,\alpha,\beta,a_{1},R,A$ and (if \eqref{condi3} holds) $\varepsilon$.
Since $a_{2}/(\kappa_{1}l)-a_{3}/l=1/l>1$, there is an integer $k\ge 1$
such that $a_{3}/l\le k\le a_{2}/(\kappa_{1}l)$.
In this case $t/k\in
[\kappa_{1} lr^{\alpha},lr^{\alpha}]$ for every $t\in
[a_{2}r^{\alpha},a_{3}r^{\alpha}]$. Thus by the Chapman--Kolmogorov equation and
 \eqref{e:3.7},
we have
\begin{align}
p^{b}_{B(x_{0},r)}(t,x,y)
\ge&\int_{B(x_{0},a_{1}r)}\cdots\int_{B(x_{0},a_{1}r)}p^{b}_{B(x_{0},r)}(t/k,x,z_{1})\cdots p^{b}_{B(x_{0},r)}(t/k,z_{k-1},y){\mathrm{d}z}_{1}\cdots {\mathrm{d}z}_{k-1}\nonumber\\
\ge&(C_{28}r^{-d})^{k}
 |B(x_{0},a_{1}r)|^{k-1}
=:c_{1}r^{-d}.
\nonumber
\end{align}
This completes the proof.\end{proof}

\begin{lemma}\label{cor3.6}
Suppose $D$ is an open set in $\R^{d}$. For every $x\in D$,
we use $D_{x}$ to denote the connected component of $D$ that
contains $x$. Then
$p^{b}_{D}(t,x,y)>0$ for every $t>0$ and  $x,y\in D$ with $\mathrm{dist}(D_{x},D_{y})<\eta(A)$
, where $\eta(A)$ is given in \eqref{defiforetaA}.
\end{lemma}

\begin{proof}
Fix $x,y\in D$ with $\mathrm{dist}(D_{x},D_{y})<\eta(A)$.
If $y\in D_{x}$, choose points $x=x_{0},x_{1},\cdots,x_{k}=y$ on a polygonal path in $D_{x}$ and $r\in(0,1/2]$ such that $B(x_{i},r)\subset D_{x}$ for all $i$ and $|x_{i}-x_{i-1}|<r/2$ for $1\le i\le k$. This is possible because the path is compact.
Put $B_{i}:=B(x_{i},r)$, $z_{0}:=x$ and $z_{k}:=y$. By the Chapman--Kolmogorov equation and the domain monotonicity $p^{b}_{D}\ge p^{b}_{B_{i}}$ on $B_{i}\times B_{i}$,
\[
p^{b}_{D}(t,x,y)\ge\int_{B(x_{1},r/4)}\cdots\int_{B(x_{k-1},r/4)}\prod_{i=1}^{k}p^{b}_{B_{i}}(t/k,z_{i-1},z_{i})\,{\mathrm{d}z}_{1}\cdots {\mathrm{d}z}_{k-1}.
\]
For each $i$, both $z_{i-1}$ and $z_{i}$ lie in $B(x_{i},3r/4)$. Proposition \ref{prop3.5} with $a_{1}=3/4$, $a_{2}=t/(kr^{\alpha})$ and $a_{3}=2a_{2}$ shows that each factor is bounded below by a positive constant. Hence $p^{b}_{D}(t,x,y)>0$.
If $y\not\in D_{x}$, then by the
strong Markov property, \eqref{levysystemforxb} and \eqref{jb&j}, we
have
\begin{eqnarray}
&&p^{b}_{D}(t,x,y)\nonumber\\
&=&\E_x\big[p^{b}_{D}(t-\tau^{b}_{D_{x}},X^{b}_{\tau^{b}_{D_{x}}},y);\,\tau^{b}_{D_{x}}<t\big]\nonumber\\
&\ge&\E_x\big[p^{b}_{D}(t-\tau^{b}_{D_{x}},X^{b}_{\tau^{b}_{D_{x}}},y);\,X^{b}_{\tau^{b}_{D_{x}}}\in D_{y},\ \tau^{b}_{D_{x}}<t\big]\nonumber\\
&=&\int_{0}^{t}\int_{D_{y}}p^{b}_{D}(t-s,z,y)\Big[\int_{D_{x}}p^{b}_{D_{x}}(s,x,w)j^{b}(w,z){\mathrm{d}w}\Big]{\mathrm{d}z\mathrm{d}s}\nonumber\\
&\ge&\int_{0}^{t}\int_{D_{y}}\int_{D_{x}}p^{b}_{D_{y}}(t-s,z,y)p^{b}_{D_{x}}(s,x,w)j^{b}(w,z){\mathrm{d}w\mathrm{d}z\mathrm{d}s}\nonumber\\
&\ge&\frac{1}{2}\int_{0}^{t}\int_{D_{y}}\int_{\{w\in D_{x}:\ \mathrm{dist}(w,D_{y})<\eta(A)\}}p^{b}_{D_{y}}(t-s,z,y)p^{b}_{D_{x}}(s,x,w)\bar{j}_{\eta(A)}(w,z){\mathrm{d}w\mathrm{d}z\mathrm{d}s}\nonumber\\
&>&0.\nonumber
\end{eqnarray}
Here the first equality follows from the strong Markov property at $\tau^{b}_{D_{x}}$, together with $p^{b}_{D_{x}}(t,x,y)=0$ for $y\notin D_{x}$.
The second equality follows from the L\'{e}vy system \eqref{levysystemforxb} with $T=\tau^{b}_{D_{x}}$ and $f(s,w,z)=\1_{\{s<t\}}\1_{D_{x}}(w)\1_{D_{y}}(z)\,p^{b}_{D}(t-s,z,y)$: since $D_{x}$ and $D_{y}$ are disjoint open sets, $X^{b}$ can enter $D_{y}$ from $D_{x}$ only by a jump, and
$p^{b}_{D_{x}}(s,x,\cdot)$ is the occupation density of $X^{b}_{s}$ on the event $\{s<\tau^{b}_{D_{x}}\}$.
The last inequality holds because $\{w\in D_{x}:\mathrm{dist}(w,D_{y})<\eta(A)\}$ is a non-empty open set, $p^{b}_{D_{x}}(s,x,\cdot)>0$ on $D_{x}$ and $p^{b}_{D_{y}}(t-s,\cdot,y)>0$ on $D_{y}$ by the first part of the proof, and $\bar{j}_{\eta(A)}(w,z)>0$ for $|w-z|<\eta(A)$.
\end{proof}

\begin{lemma}\label{lem:exit}
Let $D$ be a bounded open set in $\R^{d}$. There exist positive constants
$\lambda_0 = \lambda_0 (d,\allowbreak\alpha,\allowbreak\beta,\allowbreak\mathrm{diam}(D),\allowbreak A)$
and
$ C_{30} = C_{30} (d,\allowbreak\alpha,\allowbreak\beta,\allowbreak\mathrm{diam}(D),\allowbreak A) $
such that $\P_x(\tau^{b}_{D}>t)\le C_{30}\e^{-\lambda_{0}t}$ for every $x\in D$ and $t>0$.
\end{lemma}

\begin{proof}
For $x\in D$, let $A_{x}:=\{y\in\R^{d}:\ \mathrm{diam}(D)<|y-x|<2\mathrm{diam}(D)\}\subset D^{c}$. By Proposition \ref{P:2.1}(i),
\begin{equation*}
\P_x (\tau^{b}_{D} \leq 1)\ge\int_{D^c} p^b(1, x, y)\mathrm{d}y\ge c_{0}\int_{A_{x}}\bar{p}_{1}(1,c_{1}x,c_{1}y)\mathrm{d}y\ge\eps_{0}>0,
\end{equation*}
where $c_{0},c_{1}$ depend on $d,\alpha,\beta,A$, and $\eps_{0}$ depends on $d,\alpha,\beta,A$ and $\mathrm{diam}(D)$.
Hence $\sup_{x\in D} \P_x (\tau^{b}_{D} > 1)\leq 1-\eps_0$, and the Markov property gives $\sup_{x\in D} \P_x (\tau^{b}_{D} > n)\leq (1-\eps_0)^n$ for $n\ge 1$.
The assertion follows with $\lambda_0 := - \log (1-\eps_0)$ and $C_{30}:= \e^{\lambda_0}$.
\end{proof}

\medskip
\noindent{\bf Killing at an independent exponential time.}
Killing $X^{b}$ at an independent exponential time makes its Green function on the whole space finite.
This is needed in Section \ref{sec6}, where a dual process is constructed with respect to a reference measure defined through that Green function.
For each $q\ge 0$, consider the subprocess $X^{b,q}_{t}$ of $X^{b}_{t}$ killed at an independent exponential time $\zeta$ of parameter $q$: $X^{b,q}_{t}:=X^{b}_{t}$ when $t<\zeta$ and $X^{b,q}_{t}:=\partial$ when $t\ge \zeta$, where $\partial$ is the cemetery point. Set $\R^{d}_{\partial}:=\R^{d}\cup\{\partial\}$. $X^{b,0}_{t}$ is simply $X^{b}_{t}$. Let $\{P^{b,q}_{t}:\ t\ge 0\}$ be the transition semigroup of the process $X^{b,q}$. The transition density of $X^{b,q}$ is continuous and given by
\[p^{b,q}(t,x,y)=\e^{-q t}p^{b}(t,x,y)\quad\mbox{ for all }t\ge 0,\mbox{ and }  x,y\in \R^{d}.\]

\begin{proposition}
For each $q\ge 0$, the L\'{e}vy system of $X^{b,q}$ is given by $(N^{q},H^{q})$ where $H^{q}_{t}:=t$ and for $x\in\R^{d}$,
\[N^{q}(x,{\mathrm{d}y}):=j^{b}(x,y){\mathrm{d}y}\mbox{ on }\R^{d},\quad N^{q}(x,\{\partial\}):=q.\]
 This implies that for any $x\in \R^{d}$, any non-negative measurable function $f$ on $\R_{+}\times \R^{d}\times
 \R^{d}_{\partial}
 $ vanishing on $\{(s,y,y):\ s\ge 0,y\in \R^{d}\}$ and stopping time $T$ (with respect to the filtration of $X^{b,q}$),
 \begin{equation}\label{levysysforxbq}
\E_{x}\Big[\sum_{s\le T}f(s,X^{b,q}_{s-},X^{b,q}_{s})\Big]=\E_{x}\Big[\int_{0}^{T}\int_{
\R^{d}_{\partial}
}f(s,X^{b,q}_{s},y)N^{q}(X^{b,q}_{s},{\mathrm{d}y}){\mathrm{d}s}\Big].
\end{equation}
\end{proposition}

\begin{proof} This proposition is trivial for $q=0$. For $q>0$, its proof is essentially the same as that of \cite[Theorem 3.2]{KS}. We omit the details. \end{proof}

For any Borel set $G\subset \R^{d}$, we define $\tau^{b,q}_{G}:=\inf\{t>0:\ X^{b,q}_{t}\not\in G\}$. For an open set $D\subset\R^{d}$, let $X^{b,q,D}$ denote the subprocess of $X^{b,q}$ killed upon leaving $D$. By definition, $X^{b,q,D}$ is also obtained by killing $X^{b,D}$ at an independent exponential time of parameter $q$. Thus the transition density $p^{b,q}_{D}(t,x,y)$ of $X^{b,q,D}$ is related to the transition density $p^{b}_{D}(t,x,y)$ of $X^{b,D}$ as follows:
\begin{equation}\label{n1.2}
p^{b,q}_{D}(t,x,y)=\e^{-q t}p^{b}_{D}(t,x,y)\quad \mbox{ for }t\ge 0\mbox{ and }x,y\in D.
\end{equation}
Let $P^{b,q,D}_{t}$ denote the transition semigroup of $X^{b,q,D}$. By Theorem \ref{them3.4} and Lemma \ref{cor3.6}, we have the following corollaries.

\begin{corollary}\label{ncor1.1}
Let $D$ be an open set in $\R^{d}$. Then $p^{b,q}_{D}(t,x,y)$ is jointly continuous on $(0,{\infty})\times D\times D$ and satisfies
\[p^{b,q}_{D}(t+s,x,y)=\int_{D}p^{b,q}_{D}(t,x,z)p^{b,q}_{D}(s,z,y){\mathrm{d}z}\quad\mbox{ for all }t,s>0\mbox{ and }x,y\in D.\]
\end{corollary}

\begin{corollary}\label{ncor1.3}
Let $D$ be an open set in $\R^{d}$. Then $p^{b,q}_{D}(t,x,y)>0$ for every $t>0$ and $x,y\in D$ with $\mathrm{dist}(D_{x},D_{y})<\eta(A)$.
\end{corollary}

\section{Green function estimates on small sets}\label{sec5}
This section proves that $G^{b}_{D}$ is comparable to $G_{D}$ on balls and on $C^{1,1}$ open sets of small diameter.
The proof bounds $\mathcal{S}^{b}_{x}G_{D}(x,y)$ and sums the Duhamel series.
In Sections \ref{sec5}--\ref{sec7} we assume $d>\alpha$; see Remark \ref{R:dim}.

For any $\lambda,A>0$,
define $A_{\lambda}:=\lambda^{\beta-\alpha}A$ and
\begin{equation}
b_{\lambda}(x,z):=\lambda^{\beta-\alpha}b(\lambda^{-1}x,\lambda^{-1}z),\quad  x,z\in\R^{d}.
\nonumber
\end{equation}
If $\|b\|_{\infty}\le A$, then $\|b_{\lambda}\|_{\infty}\le A_{\lambda}$. A constant
$c$ depending on $D$, $b$ and $A$ (or part of them) is
\textit{scale-invariant} if it satisfies that $c(\lambda
D,b_{\lambda},A_{\lambda})=c(D,b,A)$ for every $\lambda>0$.

One checks that $\lambda X^{b}_{\lambda^{-\alpha}t}$ has the same distribution as $X^{b_{\lambda}}_{t}$, while for any open set $D$, $\lambda X^{b,D}_{\lambda^{-\alpha}t}$ has the same distribution as $X^{b_{\lambda},\lambda D}_{t}$. So for any $\lambda>0$, we have the following scaling properties:
\begin{equation*}
p^{b}(t,x,y)=\lambda^{d}p^{b_{\lambda}}(\lambda^{\alpha}t,\lambda x,\lambda y),\quad \forall x,y\in \R^{d},\ t>0.
 \end{equation*}
\begin{equation*}
p^{b}_{D}(t,x,y)=\lambda^{d}p^{b_{\lambda}}_{\lambda D}(\lambda^{\alpha}t,\lambda x,\lambda y),\quad \forall x,y\in D,\ t>0.
\end{equation*}
 Let $G^{b}_{D}(x,y):=\int_{0}^{{\infty}}p^{b}_{D}(t,x,y){\mathrm{d}t}$ be the Green function of $X^{b,D}$. If $G^{b}_{D}(x,y)$ is finite off the diagonal of $D\times D$, then for any $\lambda>0$,
\begin{equation}
G^{b}_{D}(x,y)=\lambda^{d-\alpha}G^{b_{\lambda}}_{\lambda D}(\lambda
x,\lambda y),\quad\forall  x,y\in D.\label{scalingforgb}
\end{equation}

Suppose $X$ is a rotationally symmetric $\alpha$-stable process and $X^{D}$ is the subprocess of $X$ killed upon leaving $D$.
Let $\tau_{D}$ be the first exit time of $X$ from $D$, and let
$G(x,y)$, $G_{D}(x,y)$ and $K_{D}(x,y)$ denote
the Green function of $X$, the Green function of $X^{D}$ and the Poisson kernel of $X^{D}$, respectively.
Let $B(x_{0},r)$ be a ball in $\R^{d}$. The following explicit formulas are known. For every $x,y\in\R^{d}$,
\begin{equation}
G(x,y)=2^{-\alpha}\pi^{-d/2}\Gamma(\frac{d-\alpha}{2})\Gamma(\frac{\alpha}{2})^{-1}|x-y|^{\alpha-d}.
\nonumber
\end{equation}
For every $x,y\in B(x_{0},r)$,
\begin{equation}
G_{B(x_{0},r)}(x,y)=2^{-\alpha}\pi^{-d/2}\Gamma(\frac{d}{2})\Gamma(\frac{\alpha}{2})^{-2}
\Big(\int_{0}^{z}(u+1)^{-d/2}u^{\alpha/2-1}{\mathrm{d}u}\Big)|x-y|^{\alpha-d},
\nonumber
\end{equation}
where $z=(r^{2}-|x-x_{0}|^{2})(r^{2}-|y-x_{0}|^{2})|x-y|^{-2}$. For
every $x\in B(x_{0},r)$ and $y\in \overline{B(x_{0},r)}^{c}$,
\begin{equation}
K_{B(x_{0},r)}(x,y)=c(d,\allowbreak \alpha)(r^{2}-|x-x_{0}|^{2})^{\alpha/2}(|y-x_{0}|^{2}-r^{2})^{-\alpha/2}|x-y|^{-d},\label{ballpoisson}
\end{equation}
where
$c(d,\allowbreak \alpha)=\Gamma(\frac{d}{2})\pi^{-d/2-1}\sin\frac{\pi\alpha}{2}$.
It is known (see \cite[Corollary 1.2]{CKS1}; see also \cite[Theorems 1.1 and 1.2]{CS} for $d\geq 2$ and \cite[Theorem 4.5]{Ku} for $d\geq 3$)
 that for any bounded
$C^{1,1}$ open set $D$
 with characteristics $(R_{0},\Lambda_{0})$,
there is a constant
$c_{0}=c_{0}(d,\allowbreak \alpha,\allowbreak R_{0},\allowbreak\Lambda_{0},\allowbreak\mathrm{diam}(D))>1$
such that
\begin{equation}\label{greenesti}
G_{D}(x,y)\stackrel{c_{0}}{\asymp}
\Big(1\wedge\frac{\delta_{D}(x)}{|x-y| }\Big)^{\alpha/2} \Big(1\wedge\frac{\delta_{D}(y)}{|x-y| }\Big)^{\alpha/2}
|x-y|^{\alpha-d} ,
\quad  x,y\in D.
\end{equation}
Consequently, the following $3G$-inequality holds: there is a constant
$c_1>0$ depending only on $d$, $\alpha$, $R_{0}$, $\Lambda_{0}$ and $\mathrm{diam}(D)$
such that
\begin{equation}\label{e:3G}
\frac{G_{D}(x,y)G_{D}(y,z)}{G_{D}(x,z)}\le c_1 \big(|x-y|^{\alpha-d}+|y-z|^{\alpha-d}\big)\quad
\mbox{for pairwise distinct }x,y,z\in D.
\end{equation}
Moreover, by the scaling property
 \eqref{scalingforgb},
 the constants $c_{0}$ and $c_1$ can be chosen to be scale-invariant in domain $D$.

\begin{defn}\label{harmon}
We say that a function $u$ defined on $\R^{d}$ is
$\mathcal{L}^{b}$-harmonic on an open set $D$ if for every bounded open set $U$ with
closure $\bar{U} \subset D$,
\begin{equation}
\E_x \big[ | u(X^{b}_{\tau^{b}_{U}}) | \big] < \infty \quad
\hbox{and} \quad   u(x)=\E_x\big[u(X^{b}_{\tau^{b}_{U}})\big]\quad\mbox{ for }x\in U.
\label{harmonicity}
\end{equation}
 The function
is called regular $\mathcal{L}^{b}$-harmonic if
\eqref{harmonicity} holds for $U=D$.
Here we take the usual convention that
$\E_x\big[u(X^{b}_{\tau^{b}_{D}})\big] :=\E_x\Big[u(X^{b}_{\tau^{b}_{D}});\tau^{b}_{D}
< \infty\Big]$, when $D$ is unbounded.
\end{defn}

In particular,
$G^{b}_{D}(\cdot,y)$ is $\mathcal{L}^{b}$-harmonic on $D\setminus
\{y\}$
, since
$G^{b}_{D}(x,y)=G^{b}_{U}(x,y)+\E_xG^{b}_{D}(X^{b}_{\tau^{b}_{U}},y)$
for every open set $U\subset D$. We
point out
that in general
$G^{b}_{D}(x,y)\not=G^{b}_{D}(y,x)$, and $G^{b}_{D}(x,\cdot)$ is not
$\mathcal{L}^{b}$-harmonic. The definition of $\alpha$-harmonicity
for $\Delta^{\alpha/2}$ is analogous to that of
$\mathcal{L}^{b}$-harmonicity.

\begin{lemma}\label{lemma5}
Suppose $D$ is a bounded open set in $\R^{d}$ and $A\in (0, \infty)$. For every $x\in
D$, $y\mapsto G^{b}_{D}(x,y)$ is continuous in $D\setminus \{x\}$. Moreover,
there exists a scale-invariant constant
$C_{31}=C_{31}(d,\allowbreak \alpha,\allowbreak \beta,\allowbreak D,\allowbreak A)>0$ such that
for all $b\in\mathcal{B}_{A}$,
\begin{equation} \label{3.12}
G^{b}_{D}(x,y)\le C_{31}|x-y|^{-d+\alpha},
\quad  \ x,y\in  D.
\end{equation}
\end{lemma}

\begin{proof} First we claim that there exist positive constants $c_{1}$
and $c_{2}$ depending on $d,\alpha,\beta, \mathrm{diam}(D)$ and $A$,
such that for any $1\le t<\infty$, $x,y\in D$,
 and
 $b\in\mathcal{B}_{A}$,
\begin{equation}
p^{b}_{D}(t,x,y)\le c_{1} \e^{-c_{2} t}.\label{lem3.7}
\end{equation}
Indeed, for $t\ge 1$, the Markov property, \eqref{stdensity} and Lemma \ref{lem:exit} give
$p^{b}_{D}(t,x,y)=\int_{D}p^{b}_{D}(t-1,x,z)p^{b}_{D}(1,z,y)\,\mathrm{d}z\le \sup_{z,w\in\R^{d}}p^{b}(1,z,w)\,\P_{x}(\tau^{b}_{D}>t-1)\le c_{1}\e^{-c_{2}t}$.
 By \eqref{stdensity} and \eqref{lem3.7}, we have
\begin{eqnarray}
G^{b}_{D}(x,y)&=&\int_{0}^{\infty}p^{b}_{D}(t,x,y){\mathrm{d}t}\nonumber\\
&\le&\int_{0}^{1}p^{b}(t,x,y){\mathrm{d}t}+\int_{1}^{\infty}p^{b}_{D}(t,x,y){\mathrm{d}t}\nonumber\\
&\lesssim&\int_{0}^{1}t^{-d/\alpha}\wedge\big(\frac{t}{|x-y|^{d+\alpha}}+\frac{A
t}{|x-y|^{d+\beta}}\big){\mathrm{d}t}+\int_{1}^{\infty}c_{1}
\e^{-c_{2}t}{\mathrm{d}t}\nonumber\\
&\le&(1+A|x-y|^{\alpha-\beta})\int_{0}^{1}t^{-d/\alpha}\wedge\frac{t}{|x-y|^{d+\alpha}}\,{\mathrm{d}t}+c_{1}/c_{2}\nonumber\\
&\le&(1+A\mathrm{diam}(D)^{\alpha-\beta})|x-y|^{\alpha-d}+c_{1}/c_{2}\nonumber\\
&\le&\big[(1+A\mathrm{diam}(D)^{\alpha-\beta})+c_{1}
\mathrm{diam}(D)^{d-\alpha}/c_{2}\big]|x-y|^{\alpha-d}.\nonumber
\end{eqnarray}
This proves \eqref{3.12}.
The scale-invariance of $C_{31}$ is implied by \eqref{scalingforgb}.
Moreover, by \eqref{lem3.7}, \eqref{stdensity} and the dominated convergence
theorem, one verifies that $y\mapsto G^{b}_{D}(x,y)$ is continuous in $D\setminus \{x\}$.\end{proof}

The first inequalities in the following two lemmas are proved in Lemma 3.1 and Lemma 3.2 of \cite{BKN}, respectively,
while the second ones can be proved by similar arguments.
We omit their proofs. Lemma \ref{partialf} follows from Lemma \ref{partialK}.
Indeed, for $x\in D$ let
$r\in (0,\delta_{D}(x))$ and $B:=B(x,r/2)$.
Since $\bar{B}\subset D$, $f(y)=\int_{\bar{B}^{c}}K_{B}(y,z)f(z)\,{\mathrm{d}z}$ for $y\in B$. By Lemma \ref{partialK}, differentiation under the integral sign at $y=x$ is justified, and gives the two bounds in Lemma \ref{partialf} with
$\delta_{B}(x)=r/2$. Lemma \ref{partialf} therefore follows by letting $r\to\delta_{D}(x)$.

\begin{lemma}\label{partialK}
There is a positive constant $C_{32}=C_{32}(d,\allowbreak \alpha)$ such that for
every $r>0$ and every ball $B:=B(0,r)$,
\[|\nabla_{x} K_{B}(x,z)|\le C_{32}\frac{K_{B}(x,z)}{\delta_{B}(x)},\quad |\partial_{ij}K_{B}(x,z)|\le C_{32}\frac{K_{B}(x,z)}{\delta_{B}(x)^{2}},  \quad (x,z)\in B\times \bar{B}^{c}.\]
Here $\nabla_{x}:=(\frac{\partial}{\partial x_{1}},\cdots,\frac{\partial}{\partial
x_{d}})$ and $\partial_{ij}:=\frac{\partial^{2}}{\partial x_{i}\partial x_{j}}$.
\end{lemma}

\begin{lemma}\label{partialf}
There is a positive constant $C_{33}=C_{33}(d,\allowbreak \alpha)$ such that for
an arbitrary open set $D$ in $\R^{d}$, and every
non-negative function $f$ which is $\alpha$-harmonic in $D$,
 \[
|\nabla_{x} f(x)|\le C_{33}\frac{f(x)}{\delta_{D}(x)},\quad |\partial_{ij}f(x)|\le C_{33}\frac{f(x)}{\delta_{D}(x)^{2}},
\quad  x\in D,\ i,j\in \{1,\cdots,d\}.
\]
\end{lemma}

\begin{lemma}\label{partialG}
Let $D$ be a
bounded $C^{1,1}$ open set
in $\R^{d}$
with $C^{1,1}$ characteristics $(R_0, \Lambda_0)$.
There exists a
positive constant
$C_{34}=C_{34}(d,\alpha,D):=C_{34}(d,\alpha,R_{0},\Lambda_{0},\mathrm{diam}(D))$
which is scale-invariant in $D$,
such that
 for every $x,y\in D$,
 \begin{equation*}
|\nabla_{x}G_{D}(x,y)|\le C_{34}|x-y|^{\alpha-d-1}\big(1\wedge
\frac{\delta_{D}(y)}{|x-y| }\big)^{\alpha/2} \big(1\vee  \frac{|x-y|}{\delta_{D}(x) }\big)^{1-\alpha/2},
\end{equation*}
\begin{equation}\label{partialG.2}
|\partial_{ij}G_{D}(x,y)|\le C_{34}|x-y|^{\alpha-d-2}\big(1\wedge
\frac{\delta_{D}(y)}{|x-y| }\big)^{\alpha/2}
\big( 1\vee \frac{|x-y|}{\delta_{D}(x) }\big)^{2-\alpha/2}.
\end{equation}
\end{lemma}

\begin{proof}   For each $y\in D$ and $1\leq i, j\leq d$,
 by applying Lemma \ref{partialf} to the function $x\mapsto G_{D}(x,y)$ on the set $D\setminus \{y\}$, we get
\[
|\nabla_{x} G_D(x, y) |\le c_1 \frac{G_D(x, y)}{|x-y|\wedge \delta_{D}(x)},\quad
|\partial_{ij} G_D(x, y)|\le c_1 \frac{G_D(x, y)}{(|x-y|\wedge \delta_{D}(x))^2},\quad
 x\in D\setminus\{y\}.
\]
So it follows from \eqref{greenesti} that
there is
$c_2=c_2 (d, \alpha, R_{0},\Lambda_{0},\mathrm{diam}(D))>0$
so that
\begin{eqnarray*}
|\nabla_{x} G_D(x, y) |
&\le&c_2
 |x-y|^{\alpha-d-1}
\big(1\wedge\frac{\delta_{D}(y)}{|x-y|}\big)^{\alpha/2}
\big(1\wedge\frac{\delta_{D}(x)}{|x-y|}\big)^{\alpha/2}
\big(1\vee \frac{|x-y|}{\delta_D (x)} \big) \nonumber \\
&\leq &
c_2
 |x-y|^{\alpha-d-1}
\big(1\wedge\frac{\delta_{D}(y)}{|x-y|}\big)^{\alpha/2}
 \big(1\vee \frac{|x-y|}{\delta_D (x)} \big)^{1-\alpha/2}.
\end{eqnarray*}
 The proof of \eqref{partialG.2} is similar.\end{proof}

 For $x,y\in D$ and $x\not=y$,
define
\begin{equation}
\nonumber
h_{D}(x,y):=\begin{cases}
         |x-y|^{\alpha-\beta-d}\big(1\wedge\frac{\delta_{D}(y)}{|x-y|}\big)^{\alpha/2}
         \quad &\hbox{if } \alpha >2 \beta , \smallskip \\
         |x-y|^{\beta-d}\big(1\wedge\frac{\delta_{D}(y)}{|x-y|}\big)^{\beta}
         \big(1\vee\log\frac{|x-y|}{\delta_{D}(x)}\big)\quad &\hbox{if } \alpha =2\beta ,
         \smallskip \\
         |x-y|^{\alpha-\beta-d}\big(1\wedge\frac{\delta_{D}(y)}{|x-y|}\big)^{\alpha/2}
         \big(1\vee\frac{|x-y|}{\delta_{D}(x)}\big)^{\beta-\alpha/2}
         &\hbox{if }   \alpha<2 \beta ,
        \end{cases}
\end{equation}
and
\begin{align}
 |\mathcal{S}_{x}^{b}|G_{D}(x,y)
:=&\mathcal{A}(d,-\beta) \Big( \int_{|z|\le
\lambda_{D}(x,y)}\big| G_{D}(x+z,y)-G_{D}(x,y)-\nabla_{x}G_{D}(x,y)\cdot
z\big|\,\frac{|b(x,z)|}{|z|^{d+\beta}}\,{\mathrm{d}z}\nonumber\\
&+\int_{|z|>\lambda_{D}(x,y)}\big(G_{D}(x+z,y)+G_{D}(x,y)\big)\frac{|b(x,z)|}{|z|^{d+\beta}}\,{\mathrm{d}z}\Big),\label{e:4.16}
\end{align}
where
$ \lambda_{D}(x,y) :=(\delta_{D}(x)\wedge |x-y|)/2>0$.

\begin{lemma}\label{lemma7}
Let $D$ be a bounded $C^{1,1}$ open set
with characteristics $(R_{0},\Lambda_{0})$.
Then there is
a positive constant
$C_{35}=C_{35}(d,\alpha,\beta,D):=C_{35}(d,\alpha,\beta,R_{0},\Lambda_{0},\mathrm{diam}(D))$
which is scale-invariant in $D$,
such that for every bounded
function $b$ on $\R^{d}\times\R^{d}$ and all $x,y\in D$ with $x\neq y$,
\begin{equation}
|\mathcal{S}_{x}^{b}|G_{D}(x,y)\le
 C_{35}\|b\|_{\infty} h_D(x, y).
\label{lem7.main}
\end{equation}
\end{lemma}

\begin{proof} By definition,
\begin{eqnarray}
&&|\mathcal{S}_{x}^{b}|G_{D}(x,y)\nonumber\\
&\le&\mathcal{A}(d,-\beta)\|b\|_{\infty}\big(\int_{|z|\le
\lambda_{D}(x,y)}|G_{D}(x+z,y)-G_{D}(x,y)-\nabla_{x}G_{D}(x,y)\cdot
z|\,|z|^{-d-\beta}\,{\mathrm{d}z}\nonumber\\
&&+\int_{|z|>\lambda_{D}(x,y),x+z\in
D}G_{D}(x+z,y)|z|^{-d-\beta}{\mathrm{d}z}+\int_{|z|>\lambda_{D}(x,y)}G_{D}(x,y)|z|^{-d-\beta}\,{\mathrm{d}z}\big)\nonumber\\
&=:&\mathcal{A}(d,-\beta)\|b\|_{\infty}(I(x,y)+II(x,y)+III(x,y)).\nonumber
\end{eqnarray}
Define $r_{D}(x,y):=\delta_{D}(x)+\delta_{D}(y)+|x-y|$. Since
$\delta_{D}(y)\le \delta_{D}(x)+|x-y|$, we have
$\delta_{D}(x)+|x-y|\le r_{D}(x,y)\le 2(\delta_{D}(x)+|x-y|)$, in
other words, we have $r_{D}(x,y)\asymp \delta_{D}(x)+|x-y|\asymp
\delta_{D}(y)+|x-y|$. It is known that for every $a,b,p>0$,
\begin{equation}
a\wedge b\asymp ab/(a+b), \quad a\vee b\asymp a+b,\quad
a^{p}+b^{p}\asymp (a+b)^{p}.\label{fact}
\end{equation}
Thus $
\lambda_{D}(x,y)
=(\delta_{D}(x)\wedge |x-y|)/2\asymp \delta_{D}(x)|x-y|/r_{D}(x,y)$.
 Using \eqref{fact} repeatedly, we have
\begin{eqnarray}
III(x,y)&\le&c_{0}\beta^{-1}|x-y|^{\alpha-d}\big(1\wedge\frac{\delta_{D}(y)^{\alpha/2}}{|x-y|^{\alpha/2}}\big)
\big(1\wedge\frac{\delta_{D}(x)^{\alpha/2}}{|x-y|^{\alpha/2}}\big)\int_{|z|>\lambda_{D}(x,y)}|z|^{-d-\beta}{\mathrm{d}z}\nonumber\\
&=&c_{0}\beta^{-1}|x-y|^{\alpha-d}\big(1\wedge\frac{\delta_{D}(y)^{\alpha/2}}{|x-y|^{\alpha/2}}\big)
\big(1\wedge\frac{\delta_{D}(x)^{\alpha/2}}{|x-y|^{\alpha/2}}\big)\lambda_{D}(x,y)^{-\beta}\nonumber\\
&\asymp&c_{0}\beta^{-1}|x-y|^{\alpha-d}\big(1\wedge\frac{\delta_{D}(y)^{\alpha/2}}{|x-y|^{\alpha/2}}\big)\frac{\delta_{D}(x)^{\alpha/2}}{r_{D}(x,y)^{\alpha/2}}\,\frac{r_{D}(x,y)^{\beta}}{\delta_{D}(x)^{\beta}|x-y|^{\beta}}\nonumber\\
&\asymp&c_{0}\beta^{-1}|x-y|^{\alpha-d-\beta}\big(1\wedge\frac{\delta_{D}(y)^{\alpha/2}}{|x-y|^{\alpha/2}}\big)\big(\frac{r_{D}(x,y)}{\delta_{D}(x)}\big)^{\beta-\alpha/2}.\label{III}
\end{eqnarray}

Next we deal with $I$. For $|z|\le\lambda_{D}(x,y)$, by
\eqref{partialG.2},
\begin{align}
&|G_{D}(x+z,y)-G_{D}(x,y)-\nabla_{x}G_{D}(x,y)\cdot
z|\nonumber\\
\le&\frac{1}{2}|z|^{2}\sup_{|u|\le\lambda_{D}(x,y)}\sum_{1\le
i,j\le
d}|\partial_{ij}G_{D}(x+u,y)|\nonumber\\
\le&\frac{1}{2}|z|^{2}d^{2}C_{34}\sup_{|u|\le\lambda_{D}(x,y)}|x+u-y|^{\alpha-d-2}\big(1\wedge\frac{\delta_{D}(y)^{\alpha/2}}{|x+u-y|^{\alpha/2}}\big)
\big(1\vee\frac{|x+u-y|^{2-\alpha/2}}{\delta_{D}(x+u)^{2-\alpha/2}}\big).\label{lem7.eq1}
\end{align}
For every $|u|\le
\lambda_{D}(x,y)
=\frac{1}{2}(\delta_{D}(x)\wedge|x-y|)$, we have $|x-y|/2\le
|x+u-y|\le 3|x-y|/2$ and $\delta_{D}(x+u)\ge \delta_{D}(x)-|u|\ge
\delta_{D}(x)/2$. Thus
\[
 \mbox{RHS of }\eqref{lem7.eq1}
\lesssim\frac{1}{2}|z|^{2}d^{2}C_{34}|x-y|^{\alpha-d-2}\big(1\wedge\frac{\delta_{D}(y)^{\alpha/2}}{|x-y|^{\alpha/2}}\big)
\big(1\vee\frac{|x-y|^{2-\alpha/2}}{\delta_{D}(x)^{2-\alpha/2}}\big),\]
and consequently,
\begin{align}
I(x,y)\lesssim&\frac{1}{2}d^{2}C_{34}(2-\beta)^{-1}|x-y|^{\alpha-d-2}\big(1\wedge\frac{\delta_{D}(y)^{\alpha/2}}{|x-y|^{\alpha/2}}\big)
\big(1\vee\frac{|x-y|^{2-\alpha/2}}{\delta_{D}(x)^{2-\alpha/2}}\big)\int_{|z|\le\lambda_{D}(x,y)}|z|^{2-d-\beta}{\mathrm{d}z}\nonumber\\
=&\frac{1}{2}d^{2}C_{34}(2-\beta)^{-1}|x-y|^{\alpha-d-2}\big(1\wedge\frac{\delta_{D}(y)^{\alpha/2}}{|x-y|^{\alpha/2}}\big)
\big(1\vee\frac{|x-y|^{2-\alpha/2}}{\delta_{D}(x)^{2-\alpha/2}}\big)\lambda_{D}(x,y)^{2-\beta}\nonumber\\
\asymp&\frac{1}{2}d^{2}C_{34}(2-\beta)^{-1}|x-y|^{\alpha-d-\beta}\big(1\wedge\frac{\delta_{D}(y)^{\alpha/2}}{|x-y|^{\alpha/2}}\big)\big(\frac{r_{D}(x,y)}{\delta_{D}(x)}\big)^{\beta-\alpha/2}.\label{I}
\end{align}
Now we deal with $II$.
 By \eqref{greenesti}, we have
 \begin{align}
II(x,y)
\le&c_{0}
\int_{\substack{x+z\in D\\ |z|>\lambda_{D}(x,y)}}
|x+z-y|^{-d+\alpha}\big(1\wedge
\frac{\delta_{D}(y)^{\alpha/2}}{|x+z-y|^{\alpha/2}}\big)\big(1\wedge
\frac{\delta_{D}(x+z)^{\alpha/2}}{|x+z-y|^{\alpha/2}}\big)|z|^{-d-\beta}{\mathrm{d}z}\nonumber\\
=&c_{0}\big(
\int_{\substack{x+z\in D\\
\lambda_{D}(x,y)<|z|< 3|x-y|/4}}
+\int_{\substack{x+z\in D\\ |z|\ge 3|x-y|/4}}
\big)\cdots {\mathrm{d}z}\nonumber\\
=:&c_{0}(IV(x,y)+V(x,y)).\nonumber
\end{align}
As for $IV$, we observe that if $\lambda_{D}(x,y)<|z|<3|x-y|/4$, we
have $\frac{1}{4}|x-y|\le |x+z-y|\le 7|x-y|/4$. By this and
\eqref{fact}, we get
\begin{align}
&IV(x,y)\nonumber\\
\asymp&|x-y|^{-d+\alpha}\big(1\wedge
\frac{\delta_{D}(y)^{\alpha/2}}{|x-y|^{\alpha/2}}\big)
\int_{\substack{x+z\in D\\ \lambda_{D}(x,y)<|z|<3|x-y|/4}}
\big(1\wedge \frac{\delta_{D}(x+z)^{\alpha/2}}
{|x+z-y|^{\alpha/2}}\big)|z|^{-d-\beta}{\mathrm{d}z}\nonumber\\
\asymp&|x-y|^{\alpha-d-\beta}\big(1\wedge
\frac{\delta_{D}(y)^{\alpha/2}}{|x-y|^{\alpha/2}}\big)
\int_{\substack{x+z\in
D\\ \lambda_{D}(x,y)<|z|<3|x-y|/4}}
|x-y|^{\beta}\,\frac{\delta_{D}(x+z)^{\alpha/2}}
{(\delta_{D}(x+z)+|x-y|)^{\alpha/2}}|z|^{-d-\beta}{\mathrm{d}z}\label{3.22}
\end{align}
We continue to estimate the integral in \eqref{3.22}. If
$\delta_{D}(x)\ge |x-y|$, then $\lambda_{D}(x,y)=|x-y|/2$, and
\begin{eqnarray}
&&\int_{\substack{x+z\in
D\\ |x-y|/2<|z|< 3|x-y|/4}}  |x-y|^{\beta}\,\frac{\delta_{D}(x+z)^{\alpha/2}}{(\delta_{D}(x+z)+|x-y|)^{\alpha/2}}|z|^{-d-\beta}{\mathrm{d}z}\nonumber\\
&\le&\int_{|x-y|/2<|z|<3|x-y|/4}|x-y|^{\beta}|z|^{-d-\beta}{\mathrm{d}z}\nonumber\\
&=&\int_{1/2}^{3/4}r^{-\beta-1}{\mathrm{d}r} < \infty.\nonumber
\end{eqnarray}
 Thus by \eqref{3.22} we have
for $\delta_{D}(x)\ge |x-y|$,
\begin{equation}
IV(x,y)\lesssim |x-y|^{\alpha-d-\beta}\big(1\wedge
\frac{\delta_{D}(y)^{\alpha/2}}{|x-y|^{\alpha/2}}\big).\label{IVcase1}
\end{equation}
Otherwise if $\delta_{D}(x)<|x-y|$, then $\lambda_{D}(x,y)=\delta_{D}(x)/2$.
Moreover, $\delta_{D}(x+z)\le \delta_{D}(x)+|z|\le 3|z|$ for any $z$
satisfying $\delta_{D}(x)/2<|z|<3|x-y|/4$.

When $\alpha> 2\beta$,
\begin{eqnarray}
&&\int_{\substack{x+z\in D\\ \delta_{D}(x)/2 <|z|<3|x-y|/4}}
|x-y|^{\beta}\,\frac{\delta_{D}(x+z)^{\alpha/2}}{(\delta_{D}(x+z)+|x-y|)^{\alpha/2}}|z|^{-d-\beta}{\mathrm{d}z}\nonumber\\
&=&\int_{\substack{x+z\in D\\ \delta_{D}(x)/2 <|z|<3|x-y|/4}}
|x-y|^{\beta-\alpha/2}|z|^{-d+\alpha/2-\beta}\frac{|x-y|^{\alpha/2}}{(\delta_{D}(x+z)+|x-y|)^{\alpha/2}}\,\frac{\delta_{D}(x+z)^{\alpha/2}}{|z|^{\alpha/2}}\,{\mathrm{d}z}\nonumber\\
&\lesssim&\int_{ \delta_{D}(x)/2 <|z|< 3|x-y|/4}|x-y|^{\beta-\alpha/2}|z|^{-d+(\alpha/2-\beta)}{\mathrm{d}z}\nonumber\\
&=&\int_{ \delta_{D}(x)/ (2|x-y|)\le|u|\le 3/4}
|u|^{-d+(\alpha/2-\beta)}{\mathrm{d}u}\nonumber\\
&\le&\int_{0}^{3/4}r^{\alpha/2-\beta-1}{\mathrm{d}r} < \infty.\label{IVcase2.1}
\end{eqnarray}

When $\alpha=2\beta$,
\begin{eqnarray}
&&\int_{\substack{x+z\in D\\ \delta_{D}(x)/2<|z|<3|x-y|/4}}
 |x-y|^{\beta}\,\frac{\delta_{D}(x+z)^{\beta}}{(\delta_{D}(x+z)+|x-y|)^{\beta}}|z|^{-d-\beta}{\mathrm{d}z}\nonumber\\
&\le&
\int_{\substack{x+z\in D\\ \delta_{D}(x)/2<|z|<3|x-y|/4}}
(3|z|)^{\beta}|z|^{-d-\beta}{\mathrm{d}z}\nonumber\\
&\asymp&\int_{1/2<|u|<3|x-y| /(4\delta_{D}(x))}|u|^{-d}{\mathrm{d}u}\nonumber\\
&\lesssim&\log\frac{|x-y|}{\delta_{D}(x)}.\label{IVcase2.2}
\end{eqnarray}

When
$\alpha<2 \beta$,
\begin{eqnarray}
&&\int_{\substack{x+z\in D\\ \delta_{D}(x)/2<|z|<3|x-y|/4}}
|x-y|^{\beta}\,\frac{\delta_{D}(x+z)^{\alpha/2}}{(\delta_{D}(x+z)+|x-y|)^{\alpha/2}}|z|^{-d-\beta}{\mathrm{d}z}\nonumber\\
&\lesssim&
\int_{\delta_{D}(x)/2<|z|<3|x-y|/4}
|x-y|^{\beta-\alpha/2}|z|^{-d-(\beta-\alpha/2)}{\mathrm{d}z}\nonumber\\
&\le&\big(\frac{|x-y|}{\delta_{D}(x)}\big)^{\beta-\alpha/2}
\int_{1/2<|u|< 3|x-y| / (4\delta_{D}(x))}
|u|^{-d-(\beta-\alpha/2)}{\mathrm{d}u}\nonumber\\
&\le&\big(\frac{|x-y|}{\delta_{D}(x)}\big)^{\beta-\alpha/2}\int_{1/2}^{\infty}r^{-(\beta-\alpha/2)-1}{\mathrm{d}r}\nonumber\\
&\lesssim&\big(\frac{|x-y|}{\delta_{D}(x)}\big)^{\beta-\alpha/2}.\label{IVcase2.3}
\end{eqnarray}
Thus we conclude from \eqref{3.22}--\eqref{IVcase2.3} that
\begin{equation}\label{IV.1}
IV(x,y)\le c_{1}|x-y|^{\alpha-d-\beta}\big(1\wedge
\frac{\delta_{D}(y)^{\alpha/2}}{|x-y|^{\alpha/2}}\big)
\quad \hbox{when } \alpha >2\beta.
\end{equation}

\begin{equation} \label{IV.2}
IV(x,y)\le
c_{1}|x-y|^{-d+\beta}\big(1\wedge\frac{\delta_{D}(y)^{\beta}}{|x-y|^{\beta}}\big)\big(1\vee
\log\frac{|x-y|}{\delta_{D}(x)}\big)\quad \hbox{when } \alpha =2 \beta,
\end{equation}
and
\begin{equation}\label{IV.3}
IV(x,y)\le c_{1}|x-y|^{\alpha-d-\beta}\big(1\wedge
\frac{\delta_{D}(y)^{\alpha/2}}{|x-y|^{\alpha/2}}\big)\big(1\vee\frac{|x-y|^{\beta-\alpha/2}}{\delta_{D}(x)^{\beta-\alpha/2}}\big)
\quad \hbox{when  }\alpha <2 \beta ,
\end{equation}
where $c_{1}=c_{1}(d,\allowbreak \alpha,\allowbreak \beta)>0$. As for $V$, note that
\begin{eqnarray}
V(x,y)&=&\int_{\substack{u\in D\\ |u-x|\ge \frac{3}{4}|x-y|}}|u-y|^{-d+\alpha}|u-x|^{-d-\beta}\big(1\wedge\frac{\delta_{D}(y)^{\alpha/2}}{|u-y|^{\alpha/2}}\big)\big(1\wedge\frac{\delta_{D}(u)^{\alpha/2}}{|u-y|^{\alpha/2}}\big){\mathrm{d}u}\nonumber\end{eqnarray}
Let $x':=x/|x-y|$ and $y':=y/|x-y|$. On one hand,
\begin{eqnarray}
V(x,y)&\le&\int_{|u-x|\ge
\frac{3}{4}|x-y|}|u-y|^{-d+\alpha}|u-x|^{-d-\beta}{\mathrm{d}u}\nonumber\\
&=&|x-y|^{\alpha-d-\beta}\int_{|u'-x'|\ge
3/4}|u'-y'|^{-d+\alpha}|u'-x'|^{-d-\beta}{\mathrm{d}u}'\nonumber\\
&=&|x-y|^{\alpha-d-\beta}\int_{|v-(x'-y')|\ge
3/4}|v|^{-d+\alpha}(|v-(x'-y')|\vee
\frac{3}{4})^{-d-\beta}{\mathrm{d}v}\nonumber\\
&\lesssim&|x-y|^{\alpha-d-\beta}\int_{|v-(x'-y')|\ge
3/4}|v|^{-d+\alpha}(|v-(x'-y')|+
\frac{3}{4})^{-d-\beta}{\mathrm{d}v}\nonumber\\
&\le&|x-y|^{\alpha-d-\beta}\int_{\R^{d}}|v|^{-d+\alpha}(||v|-1|+\frac{3}{4})^{-d-\beta}{\mathrm{d}v}\nonumber\\
 &\lesssim&|x-y|^{\alpha-d-\beta}\label{3.31}
\end{eqnarray}
On the other hand,
 \begin{align}
V(x,y)\le&\int_{\substack{u\in D\\ |u-x|\ge \frac{3}{4}|x-y|}}|u-y|^{-d+\alpha}|u-x|^{-d-\beta}\big(1\wedge\frac{\delta_{D}(y)^{\alpha/2}}{|u-y|^{\alpha/2}}\big){\mathrm{d}u}\nonumber\\
\asymp&|x-y|^{\alpha-d-\beta}\frac{\delta_{D}(y)^{\alpha/2}}{|x-y|^{\alpha/2}}\int_{\substack{u\in D\\ |u-x|\ge
\frac{3}{4}|x-y|}}|x-y|^{d+\beta-\alpha/2}|u-y|^{-d+\alpha}|u-x|^{-d-\beta}r_{D}(u,y)^{-\alpha/2}{\mathrm{d}u}\nonumber\\
\le&|x-y|^{\alpha-d-\beta}\frac{\delta_{D}(y)^{\alpha/2}}{|x-y|^{\alpha/2}}\int_{\substack{u\in D\\ |u-x|\ge
\frac{3}{4}|x-y|}}|x-y|^{d+\beta-\alpha/2}|u-y|^{-d+\alpha/2}|u-x|^{-d-\beta}{\mathrm{d}u}\nonumber\\
\le&|x-y|^{\alpha-d-\beta}\frac{\delta_{D}(y)^{\alpha/2}}{|x-y|^{\alpha/2}}\int_{|u'-x'|\ge
\frac{3}{4}}|u'-y'|^{-d+\alpha/2}|u'-x'|^{-d-\beta}{\mathrm{d}u}'\nonumber\\
\asymp&|x-y|^{\alpha-d-\beta}\frac{\delta_{D}(y)^{\alpha/2}}{|x-y|^{\alpha/2}}\int_{|v-(x'-y')|\ge
\frac{3}{4}}|v|^{-d+\alpha/2}(|v-(x'-y')|+\frac{3}{4})^{-d-\beta}{\mathrm{d}v}\nonumber\\
\le&|x-y|^{\alpha-d-\beta}\frac{\delta_{D}(y)^{\alpha/2}}{|x-y|^{\alpha/2}}\int_{\R^{d}}|v|^{-d+\alpha/2}(||v|-1|+\frac{3}{4})^{-d-\beta}{\mathrm{d}v}\nonumber\\
\lesssim&|x-y|^{\alpha-d-\beta}\frac{\delta_{D}(y)^{\alpha/2}}{|x-y|^{\alpha/2}}.\label{3.32}
\end{align}
Therefore by \eqref{3.31} and \eqref{3.32} we have
 \begin{equation}
 V(x,y)\le
 c_{2}|x-y|^{\alpha-d-\beta}\big(1\wedge\frac{\delta_{D}(y)^{\alpha/2}}{|x-y|^{\alpha/2}}\big)\label{V}
 \end{equation}
for some positive constant $c_{2}=c_{2}(d,\allowbreak \alpha,\allowbreak \beta)$. The proof is completed by combining
\eqref{III}, \eqref{I}, \eqref{IV.1}, \eqref{IV.2}, \eqref{IV.3} and
 \eqref{V}, and using the fact that for $\alpha \ge 2 \beta$,
 \[\frac{r_{D}(x,y)^{\beta-\alpha/2}}{\delta_{D}(x)^{\beta-\alpha/2}}\asymp 1\wedge
 \frac{\delta_{D}(x)^{\alpha/2-\beta}}{|x-y|^{\alpha/2-\beta}},\]
while for $\alpha <2 \beta$,
\[\frac{r_{D}(x,y)^{\beta-\alpha/2}}{\delta_{D}(x)^{\beta-\alpha/2}}\asymp
1\vee
 \frac{|x-y|^{\beta-\alpha/2}}{\delta_{D}(x)^{\beta-\alpha/2}}.\] \end{proof}

By \eqref{e:4.16} and Lemma \ref{lemma7},
we have for any bounded function $b$ satisfying
\eqref{condi1},
\begin{equation}
|\mathcal{S}_{x}^{b}G_{D}(x,y)|\le
|\mathcal{S}_{x}^{b}|G_{D}(x,y)\le
 C_{35}\|b\|_{\infty}h_{D}(x,y)
 \quad\mbox{ for }x,y\in D,\ x\not=y.
 \label{star}
\end{equation}
When $\alpha=2\beta$, the kernel $h_{D}$ contains a logarithmic factor. For the estimates below only upper bounds of $h_{D}$ are needed, and it is simpler to dominate this factor by a small power. For $x,y\in D$ with $x\not=y$, define
\begin{equation}\label{e:hhat}
\wh{h}_{D}(x,y):=\begin{cases}
h_{D}(x,y)&\hbox{if }\alpha\not=2\beta,\smallskip\\
|x-y|^{\beta-d}\big(1\wedge\frac{\delta_{D}(y)}{|x-y|}\big)^{\beta}\big(1\vee\frac{|x-y|}{\delta_{D}(x)}\big)^{\beta/2}&\hbox{if }\alpha=2\beta.
\end{cases}
\end{equation}
Since $1\vee\log u\le (2/\beta)(1\vee u)^{\beta/2}$ for $u>0$, we have $h_{D}\le (2/\beta)\wh{h}_{D}$. Hence, by \eqref{star},
\begin{equation}\label{e:starhat}
|\mathcal{S}_{x}^{b}|G_{D}(x,y)\le 2\beta^{-1}C_{35}\|b\|_{\infty}\,\wh{h}_{D}(x,y)\quad\mbox{for }x,y\in D,\ x\not=y.
\end{equation}
Note that $\wh{h}_{\lambda D}(\lambda x,\lambda y)=\lambda^{\alpha-\beta-d}\wh{h}_{D}(x,y)$ for every $\lambda>0$.

For every $x,y\in D$ with $x\not=y$, define
\[
g_{D}(x,y):=|x-y|^{\alpha-d}
\big(1\wedge\frac{\delta_{D}(x)}{|x-y| }\big)^{\alpha/2} \big(1\wedge\frac{\delta_{D}(y)}{|x-y| }\big)^{\alpha/2}.
\]

\begin{lemma}\label{lemma8}
Let $D$ be a bounded $C^{1,1}$ open set.
There exists a constant
$ C_{36}= C_{36}(d,\allowbreak \alpha,\allowbreak \beta)>0$ such that the following hold for every pairwise distinct $x,y,z\in
D$.
\begin{enumerate}[label=\upshape(\roman*), ref=\upshape(\roman*)]
\item If $\alpha >2 \beta$, then
\begin{equation}
\frac{g_{D}(x,z)\wh{h}_{D}(z,y)}{g_{D}(x,y)}\le
 C_{36}\,\mathrm{diam}(D)^{\alpha/2-\beta}\big(\frac{1}{|x-z|^{d-\alpha/2}}+\frac{1}{|y-z|^{d-\alpha/2}}\big);\label{7>}
\end{equation}
if $\alpha=2\beta$, then
\begin{equation}
\frac{g_{D}(x,z)\wh{h}_{D}(z,y)}{g_{D}(x,y)}\le
 C_{36}\big(1+\mathrm{diam}(D)^{\beta/2}\big)\big(\frac{1}{|x-z|^{d-\beta/2}}+\frac{1}{|y-z|^{d-\beta/2}}\big);\label{7=}
\end{equation}
 if $\alpha < 2\beta$, then
\begin{equation}
\frac{g_{D}(x,z)\wh{h}_{D}(z,y)}{g_{D}(x,y)}\le
 C_{36}\,\big(\frac{1}{|x-z|^{d-\alpha+\beta}}+\frac{1}{|y-z|^{d-\alpha+\beta}}\big).\label{7<}
\end{equation}
\item
\begin{equation}
\frac{\wh{h}_{D}(x,z)\wh{h}_{D}(z,y)}{\wh{h}_{D}(x,y)}\le
 C_{36}\,\big(\frac{1}{|x-z|^{d-\alpha+\beta}}+\frac{1}{|y-z|^{d-\alpha+\beta}}\big);
 \label{8.1}
\end{equation}
\end{enumerate}

\end{lemma}

\begin{proof}  (i) Recall that $\wh{h}_{D}=h_{D}$ if $\alpha\not=2\beta$. Let
$f_{D}(x,y):=|x-y|^{-d+\alpha-\beta}\big(1\wedge\frac{\delta_{D}(y)}{|x-y|}\big)^{\alpha/2}$
 and $r_{D}(x,y):=\delta_{D}(x)+\delta_{D}(y)+|x-y|$
 for any $x,y\in D$ with $x\not=y$.
 Since $\delta_{D}(y)\le |x-y|+\delta_{D}(x)$, we have $r_{D}(x,y)\asymp |x-y|+\delta_{D}(x)$.
 Using \eqref{fact}, we have
\[1\wedge\frac{\delta_{D}(y)^{\alpha/2}}{|x-y|^{\alpha/2}}\asymp \frac{\delta_{D}(y)^{\alpha/2}}{r_{D}(x,y)^{\alpha/2}},\]
\[\big(1\wedge\frac{\delta_{D}(x)^{\alpha/2}}{|x-y|^{\alpha/2}}\big)\big(1\wedge\frac{\delta_{D}(y)^{\alpha/2}}{|x-y|^{\alpha/2}}\big)\asymp \frac{\delta_{D}(x)^{\alpha/2}\delta_{D}(y)^{\alpha/2}}{r_{D}(x,y)^{\alpha}}.\]
 Thus
\begin{equation}
\frac{g_{D}(x,z)f_{D}(z,y)}{g_{D}(x,y)}\asymp\frac{|x-y|^{d-\alpha}}{|x-z|^{d-\alpha}|y-z|^{d-\alpha+\beta}}
\big(\frac{\delta_{D}(z)r_{D}(x,y)}{r_{D}(y,z)r_{D}(x,z)}\big)^{\alpha/2}
\big(\frac{r_{D}(x,y)}{r_{D}(x,z)}\big)^{\alpha/2}.\label{8'}
\end{equation}
Moreover,
\begin{equation}
\frac{\delta_{D}(z)r_{D}(x,y)}{r_{D}(y,z)r_{D}(x,z)}\le\frac{\delta_{D}(z)(r_{D}(x,z)+r_{D}(y,z))}{r_{D}(y,z)r_{D}(x,z)}=\frac{\delta_{D}(z)}{r_{D}(y,z)}+\frac{\delta_{D}(z)}{r_{D}(x,z)}\le
2,\label{10}
\end{equation}
and
\begin{equation}
\frac{r_{D}(x,y)}{r_{D}(x,z)}\asymp
\frac{|x-y|+\delta_{D}(x)}{|x-z|+\delta_{D}(x)}\le
1+\frac{|x-y|}{|x-z|}. \label{9}
\end{equation}
If $|x-z|>|x-y|/2$, then $\mbox{RHS of }\eqref{9}\le 3$, and consequently by \eqref{8'} and
\eqref{10},
\begin{eqnarray}
\frac{g_{D}(x,z)f_{D}(z,y)}{g_{D}(x,y)}&\lesssim&\frac{|x-y|^{d-\alpha}}{|x-z|^{d-\alpha}|y-z|^{d-\alpha+\beta}}\nonumber\\
&\lesssim&\frac{|x-z|^{d-\alpha}+|y-z|^{d-\alpha}}{|x-z|^{d-\alpha}|y-z|^{d-\alpha+\beta}}\nonumber\\
&\lesssim&\frac{1}{|x-z|^{d-\alpha+\beta}}+\frac{1}{|y-z|^{d-\alpha+\beta}}.\label{2.41}
\end{eqnarray}
Otherwise if $|x-z|\le |x-y|/2$, then $\mbox{RHS of }\eqref{9}\le
\frac{3}{2}|x-y|/|x-z|$, and consequently,
\begin{eqnarray}
\frac{g_{D}(x,z)f_{D}(z,y)}{g_{D}(x,y)}&\lesssim&\frac{|x-y|^{d-\alpha}}{|x-z|^{d-\alpha}|y-z|^{d-\alpha+\beta}}\,\frac{|x-y|^{\alpha/2}}{|x-z|^{\alpha/2}}\nonumber\\
&=&\frac{|x-y|^{d-\alpha/2}}{|x-z|^{d-\alpha/2}|y-z|^{d-\alpha+\beta}}\nonumber\\
&\lesssim&\frac{1}{|y-z|^{d-\alpha+\beta}}+\frac{1}{|x-z|^{d-\alpha/2}|y-z|^{\beta-\alpha/2}}.\label{11}
\end{eqnarray}
If, in addition, $\alpha > 2\beta$, then
\begin{eqnarray}
 \mbox{RHS of }\eqref{11}
&=&\frac{|y-z|^{\alpha/2-\beta}}{|y-z|^{d-\alpha/2}}+\frac{|y-z|^{\alpha/2-\beta}}{|x-z|^{d-\alpha/2}}\nonumber\\
&\le&\mathrm{diam}(D)^{\alpha/2-\beta}\big(\frac{1}{|x-z|^{d-\alpha/2}}+\frac{1}{|y-z|^{d-\alpha/2}}\big).\label{12}
\end{eqnarray}
Since
$
\wh{h}_{D}(z,y)
=f_{D}(z,y)$
for $\alpha>2\beta$,
 \eqref{7>} follows from
\eqref{2.41}, \eqref{11} and \eqref{12}.
On the other hand, if $\alpha\le 2\beta$, then for $|x-z|\le|x-y|/2$,
\begin{equation*}
 \mbox{RHS of }\eqref{11}
\lesssim
\frac{1}{|x-z|^{d-\alpha+\beta}}+\frac{1}{|y-z|^{d-\alpha+\beta}}.
\end{equation*}
This combined with \eqref{2.41} yields that, for $\alpha\le 2\beta$,
\begin{equation}
\frac{g_{D}(x,z)f_{D}(z,y)}{g_{D}(x,y)}\lesssim\frac{1}{|x-z|^{d-\alpha+\beta}}+\frac{1}{|y-z|^{d-\alpha+\beta}}.\label{13}
\end{equation}
If $\alpha=2\beta$, then $\wh{h}_{D}(z,y)=f_{D}(z,y)\big(1\vee(|y-z|/\delta_{D}(z))^{\beta/2}\big)$.
On $\{|y-z|\le\delta_{D}(z)\}$ the last factor equals $1$, and \eqref{13} gives the bound $|x-z|^{-d+\beta}+|y-z|^{-d+\beta}$.
On $\{|y-z|>\delta_{D}(z)\}$ we have $r_{D}(y,z)\asymp|y-z|$. Hence, by \eqref{8'} and $r_{D}(x,y)\le r_{D}(x,z)+r_{D}(y,z)$,
\begin{align*}
\frac{g_{D}(x,z)\wh{h}_{D}(z,y)}{g_{D}(x,y)}
\asymp&\frac{|x-y|^{d-2\beta}}{|x-z|^{d-2\beta}|y-z|^{d-\beta}}\,
\frac{\delta_{D}(z)^{\beta/2}r_{D}(x,y)^{2\beta}}{|y-z|^{\beta/2}r_{D}(x,z)^{2\beta}}\\
\lesssim&\frac{|x-y|^{d-2\beta}}{|x-z|^{d-2\beta}|y-z|^{d-\beta}}\Big(\big(\frac{\delta_{D}(z)}{|y-z|}\big)^{\beta/2}+\big(\frac{\delta_{D}(z)}{r_{D}(x,z)}\big)^{\beta/2}\frac{|y-z|^{3\beta/2}}{r_{D}(x,z)^{3\beta/2}}\Big)\\
\le&\frac{|x-y|^{d-2\beta}}{|x-z|^{d-2\beta}|y-z|^{d-\beta}}\Big(1+\frac{|y-z|^{3\beta/2}}{|x-z|^{3\beta/2}}\Big).
\end{align*}
We bound the right-hand side using $|x-y|^{d-2\beta}\lesssim|x-z|^{d-2\beta}+|y-z|^{d-2\beta}$, which holds since $d>2\beta$, and Young's inequality $a^{-p}b^{-q}\le a^{-p-q}+b^{-p-q}$ for $a,b,p,q>0$.
The first term is at most a constant multiple of $|x-z|^{-d+\beta}+|y-z|^{-d+\beta}$. The second term is at most a constant multiple of
\[
|x-z|^{-3\beta/2}|y-z|^{-(d-5\beta/2)}+|x-z|^{-(d-\beta/2)}|y-z|^{\beta/2}.
\]
The last summand is at most $\mathrm{diam}(D)^{\beta/2}|x-z|^{-(d-\beta/2)}$. If
$d
>
5\beta/2$, the first summand is at most $|x-z|^{-d+\beta}+|y-z|^{-d+\beta}$ by Young's inequality. If
$d
\le
5\beta/2$, then $|y-z|^{5\beta/2-d}\le\mathrm{diam}(D)^{5\beta/2-d}$ and $|x-z|^{-3\beta/2}\le\mathrm{diam}(D)^{d-2\beta}|x-z|^{-(d-\beta/2)}$, so the first summand is at most $\mathrm{diam}(D)^{\beta/2}|x-z|^{-(d-\beta/2)}$.
Finally, $|w|^{-d+\beta}\le\mathrm{diam}(D)^{\beta/2}|w|^{-(d-\beta/2)}$ for $0<|w|\le\mathrm{diam}(D)$. Combining these bounds proves \eqref{7=}.

If $\alpha < 2\beta$,
note that
\begin{align}
&\frac{g_{D}(x,z)
\wh{h}_{D}
(z,y)}{g_{D}(x,y)}\nonumber\\
=&\frac{g_{D}(x,z)f_{D}(z,y)}{g_{D}(x,y)}\big(1\vee
\frac{|y-z|^{\beta-\alpha/2}}{\delta_{D}(z)^{\beta-\alpha/2}}\big)\nonumber\\
=&\frac{g_{D}(x,z)f_{D}(z,y)}{g_{D}(x,y)}\1_{\{\delta_{D}(z)\ge|y-z|\}}+\frac{g_{D}(x,z)f_{D}(z,y)}{g_{D}(x,y)}\,\frac{|y-z|^{\beta-\alpha/2}}{\delta_{D}(z)^{\beta-\alpha/2}}\1_{\{\delta_{D}(z)<|y-z|\}}.\nonumber\\
=:&I_{3}(x,y,z)+I_{4}(x,y,z).\label{III+IV}
\end{align}
By \eqref{13},
\begin{equation}
I_{3}(x,y,z)\lesssim
\frac{1}{|x-z|^{d-\alpha+\beta}}+\frac{1}{|y-z|^{d-\alpha+\beta}}.\label{lem8.1}
\end{equation}
For $I_{4}$, since $r_{D}(y,z)\asymp |y-z|$ for $y,z\in D$ with
$\delta_{D}(z)<|y-z|$, we have
\begin{align}
I_{4}(x,y,z)\asymp&\frac{|x-y|^{d-\alpha}}{|x-z|^{d-\alpha}|y-z|^{d-\alpha+\beta}}\,\frac{\delta_{D}(z)^{\alpha/2}r_{D}(x,y)^{\alpha}}{r_{D}(y,z)^{\alpha/2}r_{D}(x,z)^{\alpha}}\,\frac{|y-z|^{\beta-\alpha/2}}{\delta_{D}(z)^{\beta-\alpha/2}}\1_{\{\delta_{D}(z)<|y-z|\}}\nonumber\\
\asymp&\frac{|x-y|^{d-\alpha}}{|x-z|^{d-\alpha}|y-z|^{d-\alpha+\beta}}\,\frac{\delta_{D}(z)^{\alpha-\beta}r_{D}(x,y)^{\alpha}}{|y-z|^{\alpha-\beta}r_{D}(x,z)^{\alpha}}\1_{\{\delta_{D}(z)<|y-z|\}},\label{114}
\end{align}
Moreover,
\begin{align}
\frac{\delta_{D}(z)^{\alpha-\beta}r_{D}(x,y)^{\alpha}}{|y-z|^{\alpha-\beta}r_{D}(x,z)^{\alpha}}\1_{\{\delta_{D}(z)<|y-z|\}}
\lesssim&\frac{\delta_{D}(z)^{\alpha-\beta}}{|y-z|^{\alpha-\beta}}\big(1+\frac{r_{D}(y,z)^{\alpha}}{r_{D}(x,z)^{\alpha}}\big)\1_{\{\delta_{D}(z)<|y-z|\}}\nonumber\\
\le&1+\frac{\delta_{D}(z)^{\alpha-\beta}}{r_{D}(x,z)^{\alpha-\beta}}\frac{r_{D}(y,z)^{\alpha}}{|y-z|^{\alpha-\beta}r_{D}(x,z)^{\beta}}\1_{\{\delta_{D}(z)<|y-z|\}}\nonumber\\
\le&1+\frac{|y-z|^{\beta}}{|x-z|^{\beta}}.
\nonumber
\end{align}
Thus
\begin{align}
 \mbox{RHS of }\eqref{114}
\lesssim&\frac{|x-y|^{d-\alpha}}{|x-z|^{d-\alpha}|y-z|^{d-\alpha+\beta}}+\frac{|x-y|^{d-\alpha}}{|x-z|^{d-\alpha+\beta}|y-z|^{d-\alpha}}\nonumber\\
\lesssim&\frac{1}{|x-z|^{d-\alpha+\beta}}+\frac{1}{|y-z|^{d-\alpha+\beta}}.\label{115}
\end{align}
 Consequently, \eqref{7<} follows from \eqref{III+IV}--\eqref{115}.

\medskip

\noindent(ii) If $\alpha > 2\beta$, by \eqref{10} we have
\begin{align*}
\frac{\wh{h}_{D}(x,z)\wh{h}_{D}(z,y)}{\wh{h}_{D}(x,y)}\asymp&\frac{|x-y|^{d-\alpha+\beta}}{|x-z|^{d-\alpha+\beta}|y-z|^{d-\alpha+\beta}}\big(\frac{\delta_{D}(z)r_{D}(x,y)}{r_{D}(x,z)r_{D}(y,z)}\big)^{\alpha/2}\nonumber\\
\lesssim&\frac{|x-y|^{d-\alpha+\beta}}{|x-z|^{d-\alpha+\beta}|y-z|^{d-\alpha+\beta}}\nonumber\\
\lesssim&\frac{1}{|x-z|^{d-\alpha+\beta}}+\frac{1}{|y-z|^{d-\alpha+\beta}}.
\end{align*}
If $\alpha < 2\beta$, by \eqref{fact} we have
\[1\vee \frac{|x-y|^{\beta-\alpha/2}}{\delta_{D}(x)^{\beta-\alpha/2}}\asymp 1+\frac{|x-y|^{\beta-\alpha/2}}{\delta_{D}(x)^{\beta-\alpha/2}}\asymp \frac{r_{D}(x,y)^{\beta-\alpha/2}}{\delta_{D}(x)^{\beta-\alpha/2}},\]
and thus
\[\wh{h}_{D}(x,y)\asymp |x-y|^{-d+\alpha-\beta}\frac{\delta_{D}(y)^{\alpha/2}}{r_{D}(x,y)^{\alpha-\beta}\delta_{D}(x)^{\beta-\alpha/2}}.\]
It follows from \eqref{10} that
\begin{align}
\frac{\wh{h}_{D}(x,z)\wh{h}_{D}(z,y)}{\wh{h}_{D}(x,y)}\asymp&\frac{|x-y|^{d-\alpha+\beta}}{|x-z|^{d-\alpha+\beta}|y-z|^{d-\alpha+\beta}}\big(\frac{\delta_{D}(z)r_{D}(x,y)}{r_{D}(x,z)r_{D}(y,z)}\big)^{\alpha-\beta}\nonumber\\
\lesssim&\frac{1}{|x-z|^{d-\alpha+\beta}}+\frac{1}{|y-z|^{d-\alpha+\beta}}.
\nonumber
\end{align}
If $\alpha=2\beta$, then by \eqref{fact},
\[\wh{h}_{D}(x,y)\asymp|x-y|^{\beta-d}\frac{\delta_{D}(y)^{\beta}}{r_{D}(x,y)^{\beta/2}\delta_{D}(x)^{\beta/2}},\]
and the computation above, with the exponent $\alpha-\beta$ of the last factor replaced by $\beta/2$, gives the same bound. This establishes \eqref{8.1} in all cases.
\end{proof}

 \begin{defn}
Suppose $\gamma>0$. For a function $f$ defined on $\R^{d}$,
we define for $r>0$,
\[M^{\gamma}_{f}(r):=\sup_{x\in \R^{d}}\int_{B(x,r)}\frac{|f(y)|}{|x-y|^{d-\gamma}}\,{\mathrm{d}y}.\]
$f$ is said to belong to the Kato class $\mathbb{K}_{d,\gamma}$ if
$\lim_{r\downarrow 0}M^{\gamma}_{f}(r)=0.$
For any bounded open set
$D\subset \R^{d}$, we define
\[M^{\gamma}_{f}(D):=\sup_{x\in\R^{d}}\int_{D}\frac{|f(y)|}{|x-y|^{d-\gamma}}\,{\mathrm{d}y}.\]
\end{defn}

\begin{lemma}\label{lemma11}
Let $D$ be a bounded $C^{1,1}$ open set. Then for any bounded
function $b$ satisfying \eqref{condi1},
\begin{equation}
\mathcal{S}^{b}_{x}\int_{D}G_{D}(x,z)\mathcal{S}^{b}_{z}G_{D}(z,y){\mathrm{d}z}=\int_{D}\mathcal{S}^{b}_{x}G_{D}(x,z)\mathcal{S}^{b}_{z}G_{D}(z,y){\mathrm{d}z}\quad  \mbox{for }x,y\in D,\ x\not=y.\label{23}
\end{equation}
Furthermore, let $\gamma  =
(\alpha-\beta)\wedge (\alpha/2)$ if
$\alpha/2\not=\beta$ and $\gamma\in (0,\beta)$ if $\alpha/2=\beta$.
Then for any measurable function $f\in \mathbb{K}_{d,\gamma}$,
\begin{equation}
\mathcal{S}^{b}_{x}\int_{D}G_{D}(x,z)f(z){\mathrm{d}z}=\int_{D}\mathcal{S}^{b}_{x}G_{D}(x,z)f(z){\mathrm{d}z}\quad
\mbox{for }
 x\in D.\label{doublestar}
\end{equation}
\end{lemma}

\begin{proof} Fix $x,y\in D, x\not=y$.
For any $\varepsilon >0$, we have by \eqref{star}
\begin{eqnarray}
&&\big|(G_{D}(x+u,z)-G_{D}(x,z))b(x,u)|u|^{-d-\beta}\mathcal{S}^{b}_{z}G_{D}(z,y)\big|
\1_{\{z\in D,|u|>\varepsilon\}}
\nonumber\\
&\le&\|b\|_{\infty}|u|^{-d-\beta}|\mathcal{S}^{b}_{z}|G_{D}(z,y)\big(G_{D}(x+u,z)\1_{\{z,x+u\in D,|u|>\varepsilon\}}+G_{D}(x,z)\1_{\{z\in D, |u|>\varepsilon\}}\big)\nonumber\\
&\le&C_{31} C_{35}\|b\|_{\infty}|u|^{-d-\beta}h_{D}(z,y)\big(|x+u-z|^{\alpha-d}\1_{\{z,x+u\in D, |u|>\varepsilon\}}+|x-z|^{\alpha-d}\1_{\{z\in
D, |u|>\varepsilon\}}\big).\nonumber
\end{eqnarray}
Thus
\begin{eqnarray}
&&\int_{|u|>\varepsilon}\int_{z\in
D}\big|(G_{D}(x+u,z)-G_{D}(x,z))b(x,u)|u|^{-d-\beta}\mathcal{S}^{b}_{z}G_{D}(z,y)\big|{\mathrm{d}z\mathrm{d}u}\nonumber\\
&\le&C_{31} C_{35}\|b\|_{\infty}\big(\int_{\substack{v\in
D\\ |v-x|>\varepsilon}}\int_{z\in D}|v-x|^{-d-\beta}|v-z|^{-d+\alpha}h_{D}(z,y){\mathrm{d}z\mathrm{d}v}\nonumber\\
&&+\int_{|u|>\varepsilon}\int_{z\in
D}|u|^{-d-\beta}|z-x|^{-d+\alpha}h_{D}(z,y){\mathrm{d}z\mathrm{d}u}\big).\label{lem11.eq1}
\end{eqnarray}
The integrals in \eqref{lem11.eq1} are
finite. Thus the integral
\[\int_{|u|>\varepsilon}\int_{z\in
D}(G_{D}(x+u,z)-G_{D}(x,z))b(x,u)|u|^{-d-\beta}\mathcal{S}^{b}_{z}G_{D}(z,y){\mathrm{d}z\mathrm{d}u}\]
is absolutely convergent. By Fubini's theorem, we have
\begin{eqnarray}
&&\mathcal{S}^{b}_{x}\int_{D}G_{D}(x,z)\mathcal{S}^{b}_{z}G_{D}(z,y){\mathrm{d}z}\nonumber\\
&=&\lim_{\varepsilon\to 0}\int_{|u|>\varepsilon}\int_{z\in
D}\big(G_{D}(x+u,z)-G_{D}(x,z)\big)b(x,u)|u|^{-d-\beta}\mathcal{S}^{b}_{z}G_{D}(z,y){\mathrm{d}z\mathrm{d}u}\nonumber\\
&=&\lim_{\varepsilon\to 0}\int_{z\in
D}\Big[\int_{|u|>\varepsilon}\big(G_{D}(x+u,z)-G_{D}(x,z)\big)b(x,u)|u|^{-d-\beta}{\mathrm{d}u}\Big]\mathcal{S}^{b}_{z}G_{D}(z,y){\mathrm{d}z}\nonumber\\
&=&\int_{z\in
D}\lim_{\varepsilon\to 0}\Big[\int_{|u|>\varepsilon}\big(G_{D}(x+u,z)-G_{D}(x,z)\big)b(x,u)|u|^{-d-\beta}{\mathrm{d}u}\Big]\mathcal{S}^{b}_{z}G_{D}(z,y){\mathrm{d}z}\nonumber\\
&=&\int_{D}\mathcal{S}^{b}_{x}G_{D}(x,z)\mathcal{S}^{b}_{z}G_{D}(z,y){\mathrm{d}z}.\nonumber
\end{eqnarray}
Here the third equality follows from the dominated convergence theorem,
since for $\lambda_{D}(x,z)=(\delta_{D}(x)\wedge|x-z|)/2$ and $\varepsilon>0$
sufficiently small, we have
\begin{eqnarray}
&&\Big|\int_{|u|>\varepsilon}\big(G_{D}(x+u,z)-G_{D}(x,z)\big)b(x,u)|u|^{-d-\beta}{\mathrm{d}u}\Big|\nonumber\\
&=&\Big|\int_{|u|>\varepsilon}\big(G_{D}(x+u,z)-G_{D}(x,z)-\nabla_{x}G_{D}(x,z)\cdot u1_{|u|<\lambda_{D}(x,z)}\big)b(x,u)|u|^{-d-\beta}{\mathrm{d}u}\Big|\nonumber\\
&\le&\int_{|u|<\lambda_{D}(x,z)}\big|G_{D}(x+u,z)-G_{D}(x,z)-\nabla_{x}G_{D}(x,z)\cdot
u\big||b(x,u)||u|^{-d-\beta}{\mathrm{d}u}\nonumber\\
&&+\int_{|u|\ge\lambda_{D}(x,z)}(G_{D}(x+u,z)+G_{D}(x,z))|b(x,u)||u|^{-d-\beta}{\mathrm{d}u}\nonumber\\
&=&|\mathcal{S}^{b}_{x}|G_{D}(x,z),\nonumber
\end{eqnarray}
and by
\eqref{e:starhat}
and
Lemma \ref{lemma8}(ii)
\begin{equation*}
\int_{z\in
D}|\mathcal{S}^{b}_{x}|G_{D}(x,z)|\mathcal{S}^{b}_{z}|G_{D}(z,y){\mathrm{d}z}
\le
4\beta^{-2}\|b\|^{2}_{\infty}
 C_{35}^{2}\int_{z\in D}\wh{h}_{D}(x,z)\wh{h}_{D}(z,y){\mathrm{d}z}
<\infty.
\end{equation*}
This establishes \eqref{23}.
To prove \eqref{doublestar}, first we note
that $\mathbb{K}_{d,\gamma}\subset
\mathbb{K}_{d,\alpha-\beta}\subset \mathbb{K}_{d,\alpha}$. For any
$\varepsilon>0$ and $x\in D$, we have
\begin{eqnarray}
&&\int_{|u|>\varepsilon}\int_{D}|G_{D}(x+u,z)-G_{D}(x,z)|\frac{|b(x,z)|}{|u|^{d+\beta}}\,|f(z)|{\mathrm{d}z\mathrm{d}u}\nonumber\\
&\le&\|b\|_{\infty}\big(\int_{\substack{u+x\in
D\\ |u|>\varepsilon}}\int_{D}G_{D}(x+u,z)|u|^{-d-\beta}|f(z)|{\mathrm{d}z\mathrm{d}u}\nonumber\\
&&+\int_{|u|>\varepsilon}\int_{D}G_{D}(x,z)|u|^{-d-\beta}|f(z)|{\mathrm{d}z\mathrm{d}u}\big)\nonumber\\
&\le&\|b\|_{\infty}C_{31}\big[\varepsilon^{-d-\beta}\int_{\substack{u+x\in
D\\ |u|>\varepsilon}}\Big(\int_{D}|x+u-z|^{\alpha-d}|f(z)|{\mathrm{d}z}\Big){\mathrm{d}u}\nonumber\\
&&+\int_{|u|>\varepsilon}|u|^{-d-\beta}\Big(\int_{D}|f(z)||x-z|^{\alpha-d}{\mathrm{d}z}\Big){\mathrm{d}u}\big]\nonumber\\
&\lesssim&\|b\|_{\infty}C_{31}\big(\varepsilon^{-d-\beta}(\mathrm{diam}(D)+|x|)^{d}+\varepsilon^{-\beta}\big)M^{\alpha}_{f}(D)<\infty.\nonumber
\end{eqnarray}
In other words,
$\int_{|u|>\varepsilon}\int_{D}(G_{D}(x+u,z)-G_{D}(x,z))\frac{b(x,z)}{|u|^{d+\beta}}\,f(z){\mathrm{d}z\mathrm{d}u}$
is absolutely convergent. We observe that for any $\gamma>0$
satisfying our assumption, $h_{D}(x,y)\le c_{1}|x-y|^{-d+\gamma}$
for some positive constant
$c_{1}=c_{1}(d,\alpha,\beta,\mathrm{diam}(D),\delta_{D}(x))$, thus
\[\int_{D}h_{D}(x,y)|f(y)|{\mathrm{d}y}\le c_{1}M^{\gamma}_{f}(D)<\infty.\]
The argument used for \eqref{23} now gives \eqref{doublestar}. \end{proof}

\begin{lemma}\label{prop4.2}
Let $D$ be a bounded $C^{1,1}$ open set in $\R^{d}$
and $A>0$.
Then  for every $b\in\mathcal{B}_{A}$,
\begin{equation}
G_{D}^{b}(x,y)=G_{D}(x,y)+\int_{D}G_{D}^{b}(x,z)\mathcal{S}^{b}_{z}G_{D}(z,y){\mathrm{d}z},\quad
\mbox{ for }x,y\in D,\ x\not=y.
\label{duhamel}
\end{equation}
\end{lemma}

\begin{proof}
The proof
is  similar to that of \cite[Proposition 4.2]{CKS2}.
By Lemmas \ref{lemma5} and \ref{lemma7}, the integral on the right hand side of
\eqref{duhamel} is absolutely convergent, and is continuous in $y\in
D\setminus \{x\}$ for every $x\in D$.
It is proved in \cite{CW} that for every $f \in C^\infty_c(\R^d)$, (see \cite[Theorem 1.3]{CW}; note that $C^{\infty}_{c}(\R^{d})\subset\mathcal{S}(\R^{d})$)
$M^f_t:=f(X^b_t)-f(X^b_0)-\int^t_0 \L^b f(X^b_s)\, {\mathrm{d}s}$ is a
$\P_{x}$-martingale.
Since $\L^bf $ is bounded for $f\in C^\infty_c(\R^d)$,
we have by Lemma \ref{lemma5} that
$M^{f}_{t\wedge\tau^{b}_{D}}$
is a uniformly integrable martingale.

Define $D_j:=\{x\in D:  \delta_D(x)> 1/j\}$. Let $\phi\in
C^\infty_c(\R^d)$ with
$\phi \le 1$,
${\rm supp} [\phi]\subset
B(0, 1)$ and $\int_{\R^d} \phi (x) {\mathrm{d}x}=1$.
For any $\psi \in
C_c(D)$, define $f:=G_D\psi$ and $f_n:=\phi_n * f$, where $\phi_n
(x):=n^d \phi (nx)$. Clearly $f_n\in C^\infty_c(\R^d)$ and $f_n$
converges uniformly to $f=G_D\psi$. Fix an arbitrary $j\geq 1$. For every $y\in D_j$ and $n>2j$, we have
\[
\Delta^{\alpha/2}f_n(y)=\phi_n * \Delta^{\alpha/2}f (y) = \int_{B(0, 1/n)}
\phi_n(z)\Delta^{\alpha/2}(G_D\psi )(y-z){\mathrm{d}z}
= -\phi_n * \psi (y);
\]
see, for example, \cite[Proposition 3.13]{BB}.
The optional stopping theorem gives $\E_x[M^{f_n}_0]=\E_x[M^{f_n}_{\tau^{b}_{D_j}}]$. Thus we have
\begin{eqnarray}\label{lem4.10.n1}
\E_x\big[f_n\big(X^b_{\tau^{b}_{D_j}}\big)\big]
-f_n(x)
 &=& \int_{D_j} G^b_{D_j}(x,y) \L^b f_n (y)\, {\mathrm{d}y}\nonumber\\
&=& \int_{D_j} G^b_{D_j}(x,y) \big(-\phi_n*\psi(y) + \S^b f_n(y)\big)\, {\mathrm{d}y} .
\end{eqnarray}
The operator $\S^{b}$ does not commute with convolution, because $b(x,z)$ depends on $x$. We argue as follows.
For $y\in\R^{d}$, let $\S^{b,y}$ denote the operator $\S^{b}$ with coefficient frozen at $y$, that is, with $b(\cdot,z)$ replaced by $b(y,z)$.
Since $\S^{b,y}$ commutes with translations, $\S^{b}f_{n}(y)=\S^{b,y}f_{n}(y)=\int_{\R^{d}}\phi_{n}(w)\,\S^{b,y}f(y-w)\,{\mathrm{d}w}$.
Lemma \ref{lemma7} applies to the bounded function $(v,z)\mapsto b(y,z)$ with the same constant $C_{35}$, and \eqref{doublestar} gives $\S^{b,y}f(v)=\int_{D}\S^{b,y}_{v}G_{D}(v,u)\psi(u)\,{\mathrm{d}u}$.
Hence $|\S^{b,y}f(v)|\le C_{35}\|b\|_{\infty}\|\psi\|_{\infty}\sup_{v\in D_{2j}}\int_{D}h_{D}(v,u)\,{\mathrm{d}u}=:c_{j}<\infty$ for all $v\in D_{2j}$ and all $y\in\R^{d}$; see the proof of Lemma \ref{lemma11}.
Moreover, $v\mapsto \S^{b,y}_{v}G_{D}(v,u)$ is continuous on $D_{2j}\setminus\{u\}$ by Lemma \ref{partialG} and dominated convergence, and is bounded by $C_{35}\|b\|_{\infty}h_{D}(v,u)\le c\,|v-u|^{-d+\gamma}$ with $\gamma$ as in Lemma \ref{lemma11}. Therefore $v\mapsto\S^{b,y}f(v)$ is continuous on $D_{2j}$.
Consequently, for $y\in D_{j}$ and $n>2j$ we have $|\S^{b}f_{n}(y)|\le c_{j}$, and $\S^{b}f_{n}(y)\to\S^{b,y}f(y)=\S^{b}f(y)$ as $n\to\infty$.
 Taking $n\to \infty$
in \eqref{lem4.10.n1},
we get, by Lemma \ref{lemma5}, the bounds above, the dominated convergence theorem and \eqref{doublestar},
\begin{equation*}
\E_x\big[f \big(X^b_{\tau^{b}_{D_j}}\big)\big]
-f (x) = \int_{D}
G^b_{D_j}(x,y) \big( -\psi (y)+ (\S^b G_D\psi) (y)
\big){\mathrm{d}y}.
\end{equation*}
Now  taking
$j\to \infty$ yields, by monotone convergence $G^{b}_{D_{j}}\uparrow G^{b}_{D}$, quasi-left continuity of $X^{b}$, and the facts that $f$ is continuous on $\R^{d}$ and vanishes on $D^{c}$,
\[
-G_D \psi (x)= -f (x)= \int_{D } G^b_{D }(x,y)
 \big(-\psi (y)+ (\S^b G_D\psi) (y) \big){\mathrm{d}y}.
\]
This shows that for each $x\in D$, \eqref{duhamel} holds for a.e. $y\in D$. Since $y\mapsto G^{b}_{D}(x,y)$
is continuous in $D\setminus \{x\}$, we establish \eqref{duhamel} for every
$x,y\in D$ with $x\not=y$.
\end{proof}

\begin{thm}\label{them4.7}
Suppose $A\in (0,\infty)$.
 There exist positive constants
$r_{1}=r_{1}(d,\allowbreak \alpha,\allowbreak \beta,\allowbreak A)$ and $ C_{37}= C_{37}(d, \alpha, \beta, A)$
such that for any
$b\in\mathcal{B}_{A}$,
any ball
$B_{r}$ with radius $r\in (0,r_{1}]$
 and any $x,y\in B_{r}$ with $x\neq y$,
\begin{equation}\label{lem13.main}
\frac{1}{2}G_{B_{r}}(x,y)\le G^{b}_{B_{r}}(x,y)\le \frac{3}{2}G_{B_{r}}(x,y)
\qquad  \hbox{and} \qquad |\mathcal{S}^{b}_{x}G_{B_{r}}^b (x,y)| \leq  C_{37} \wh{h}_{B_{r}}(x, y).
\end{equation}
Moreover, we have
\begin{equation}\label{e:4.65}
G_{B_{r}}^{b}(x,y)=\sum_{k=0}^{\infty}G_{k}(x,y),
\end{equation}
where
\begin{equation}\label{e:4.66}
G_{0}(x,y):=G_{B_{r}}(x,y) \quad \hbox{ and } \quad
G_{n}(x,y):=\int_{B_{r}}G_{n-1}(x,z)\mathcal{S}^{b}_{z}G_{B_{r}}(z,y){\mathrm{d}z}
\quad \hbox{for } n\ge 1.
\end{equation}
Furthermore, the constant $r_{1}$ satisfies
the following scaling property:
\begin{equation}
r_{1}(d,\allowbreak \alpha,\allowbreak \beta,\allowbreak A_{\lambda})=\lambda r_{1}(d,\allowbreak \alpha,\allowbreak \beta,\allowbreak A)\quad\mbox{ for all }
 \lambda>0.\label{r*scale}
\end{equation}
\end{thm}

\begin{proof} Fix $A>0$. We first prove \eqref{lem13.main}--\eqref{e:4.66} for all balls of radius $r\in (0,\tilde{r}_{1}]$, for some $\tilde{r}_{1}=\tilde{r}_{1}(d,\allowbreak \alpha,\allowbreak \beta,\allowbreak A)\in (0,1]$ and with some constant $\tilde{C}_{37}(d,\allowbreak \alpha,\allowbreak \beta,\allowbreak A)$ in place of $C_{37}$. The constants
$r_{1}$
and $C_{37}$ are then defined by scaling at the end of the proof. By \eqref{greenesti} and \eqref{e:starhat}, there are positive
constants $c_{1}=c_{1}(d,\allowbreak \alpha)$ and $c_{2}=c_{2}(d,\allowbreak \alpha,\allowbreak \beta)$
such that for every ball $B_{r}$ with radius  $r\in (0,\tilde{r}_{1}]$
and every $x,y\in B_{r}$  with $x\not=y$,
\begin{equation}
G_{B_{r}}(x,y)\stackrel{c_{1}}{\asymp}g_{B_{r}}(x,y),\label{lem13.greenesti}
\end{equation}
and
\begin{equation}
|\mathcal{S}^{b}_{x}G_{B_{r}}(x,y)|\le c_{2}A
\wh{h}_{B_{r}}(x,y).\label{lem13.star}
\end{equation}
For pairwise distinct $x,y,z\in B_{r}$, define
$$I_{B_{r}}(x,y,z):=\frac{g_{B_{r}}(x,z)\wh{h}_{B_{r}}(z,y)}{g_{B_{r}}(x,y)}\mbox{ and }J_{B_{r}}(x,y,z):=\frac{\wh{h}_{B_{r}}(x,z)\wh{h}_{B_{r}}(z,y)}{\wh{h}_{B_{r}}(x,y)}.$$
We claim that, for every $r\in (0,1]$, there exists a positive constant $C(r)=C(r,d,\alpha,\beta)$ such that $C(r)\downarrow 0$ as $r\downarrow 0$, and for every ball $B_{r}$ and every $x,y\in B_{r}$ with $x\not=y$,
\begin{align}
\int_{B_{r}}I_{B_{r}}(x,y,z)\mathrm{d}z\le& C(r),\label{lem13.1}\\
\int_{B_{r}}J_{B_{r}}(x,y,z)\mathrm{d}z\le& C(r).\label{lem13.2}
\end{align}
This claim will be proved at the end of this proof.
Let $G_{k}(x,y)$ be given by
\eqref{e:4.66}. Then by \eqref{lem13.greenesti}--\eqref{lem13.2} and \eqref{23}, for every ball $B_{r}$ with radius $r\in (0,1]$ and every  $x,y\in B_{r}$ with $x\not=y$,
\begin{align*}
|G_{1}(x,y)|\le&\int_{B_{r}}G_{B_{r}}(x,z)|\mathcal{S}^{b}_{z}G_{B_{r}}(z,y)|{\mathrm{d}z}\le c_{1}c_{2}A\int_{B_{r}}g_{B_{r}}(x,z)\wh{h}_{B_{r}}(z,y){\mathrm{d}z}\nonumber\\
=&c_{1}c_{2}A g_{B_{r}}(x,y)\int_{B_{r}}I_{B_{r}}(x,y,z)\mathrm{d}z
\le c_{1}c_{2}AC(r)g_{B_{r}}(x,y),
\end{align*}
and
\begin{align*}
|\mathcal{S}^{b}_{x}G_{1}(x,y)|\le&\int_{B_{r}}\big|\mathcal{S}^{b}_{x}G_{B_{r}}(x,z)\mathcal{S}^{b}_{z}G_{B_{r}}(z,y)\big|{\mathrm{d}z}\le (c_{2}A)^{2}\int_{B_{r}}\wh{h}_{B_{r}}(x,z)\wh{h}_{B_{r}}(z,y){\mathrm{d}z}\nonumber\\
=&(c_{2}A)^{2}\wh{h}_{B_{r}}(x,y)\int_{B_{r}}J_{B_{r}}(x,y,z)\mathrm{d}z
\le(c_{2}A)^{2}C(r)\wh{h}_{B_{r}}(x,y).
\end{align*}
For every $n\ge 1$, we have
\begin{equation}
G_{n}(x,y)=\int_{B_{r}}G_{B_{r}}(x,z)\mathcal{S}^{b}_{z}G_{n-1}(z,y){\mathrm{d}z},\label{16}
\end{equation}
and
\begin{equation*}
\mathcal{S}^{b}_{x}G_{n}(x,y)=\int_{B_{r}}\mathcal{S}^{b}_{x}G_{n-1}(x,z)
\mathcal{S}^{b}_{z}G_{B_{r}}(z,y){\mathrm{d}z}.
\end{equation*}
These equalities can be proved consecutively by induction. Thus by
\eqref{lem13.1}, \eqref{lem13.2} and induction, we have
\begin{equation}
|G_{n}(x,y)|\le c_{1}(c_{2}AC(r))^{n}g_{B_{r}}(x,y)\le
c_{1}^{2}(c_{2}AC(r))^{n}G_{B_{r}}(x,y),\label{Gn}
\end{equation}
and
\begin{equation}\label{DeltaGn}
|\mathcal{S}^{b}_{x}G_{n}(x,y)|\le   c_{2}A(c_{2}AC(r))^{n}
\wh{h}_{B_{r}}(x,y).
\end{equation}
Applying Duhamel's formula \eqref{duhamel} recursively $n$ times, we get for
$n\ge 0$ and $x,y\in B,x\not=y$,
\begin{equation}\label{18}
G_{B_{r}}^{b}(x,y)=\sum_{k=0}^{n}G_{k}(x,y)+\int_{B_{r}}G^{b}_{B_{r}}(x,z)\mathcal{S}^{b}_{z}G_{n}(z,y){\mathrm{d}z}.
\end{equation}
Now let $\tilde{r}_{1}\in (0,1]$ be sufficiently small
so that $\delta :=c_2 A C(\tilde{r}_1) \leq 1/(2c_1^{2}+1)$.
By \eqref{3.12}
and
\eqref{DeltaGn}, we have for any $r\in (0, \tilde{r}_1]$,
\begin{eqnarray}
\lim_{n\to \infty} \Big|\int_{B_{r}}G^{b}_{B_{r}}(x,z)\mathcal{S}^{b}_{z}G_{n}(z,y){\mathrm{d}z}\Big|
\le \lim_{n\to \infty}  c_{2}A \delta^n C_{31}
\int_{B_{r}}|x-z|^{-d+\alpha}\wh{h}_{B_{r}}(z,y){\mathrm{d}z}
=0.   \nonumber
\end{eqnarray}
This together with \eqref{18} establishes \eqref{e:4.65}.
The first assertion in \eqref{lem13.main} then follows from the fact that
\[  \sum_{n=1}^\infty |G_n (x, y)| \leq \sum_{n=1}^\infty c_1^2 \delta^n G_{B_{r}} (x, y)
\leq G_{B_{r}}(x, y) /2
\]
for every $r\in (0, \tilde{r}_1]$.
We next prove the second assertion in \eqref{lem13.main}.
By \eqref{e:4.65},
\begin{align}
\mathcal{S}^{b}_{x}G^{b}_{B_{r}}(x,y)=&\mathcal{A}(d,-\beta)
\lim_{\varepsilon\to 0}\int_{|z|>\varepsilon}\frac{G^{b}_{B_{r}}(x+z,y)-G^{b}_{B_{r}}(x,y)}{|z|^{d+\beta}}\,b(x,z){\mathrm{d}z}\nonumber
\\
=&\mathcal{A}(d,-\beta)\lim_{\varepsilon\to
0}\int_{|z|>\varepsilon}\big(\lim_{n\to\infty}\sum_{k=0}^{n}\frac{G_{k}(x+z,y)
-G_{k}(x,y)}{|z|^{d+\beta}}\big)b(x,z){\mathrm{d}z}.\label{2.7}
\end{align}
By \eqref{Gn} and \eqref{lem13.greenesti}, for any $n\ge 1$,
\begin{eqnarray}
\sum_{k=0}^{n}\big|G_{k}(x+z,y)-G_{k}(x,y)\big|&\le&\sum_{k=0}^{n}c_{1}\delta^{k}\big(g_{B_{r}}(x+z,y)+g_{B_{r}}(x,y)\big)\nonumber\\
&\le&c_{1}(1-\delta)^{-1}\big(|x+z-y|^{\alpha-d}+|x-y|^{\alpha-d}\big).\nonumber
\end{eqnarray}
The last term is absolutely convergent with respect to
$|b(x,z)||z|^{-d-\beta}{\mathrm{d}z}$ on the domain
$\{z\in\R^{d}:|z|>\varepsilon\}$ for any $\varepsilon>0$.
By the dominated convergence theorem, \eqref{2.7} continues as
\begin{eqnarray}
\mathcal{S}^{b}_{x}G^{b}_{B_{r}}(x,y)&=&\mathcal{A}(d,-\beta)\lim_{\varepsilon\to
0}\lim_{n\to\infty}\int_{|z|>\varepsilon}\big( \sum_{k=0}^{n}\frac{G_{k}(x+z,y)-G_{k}(x,y)}{|z|^{d+\beta}}\big)
b(x,z){\mathrm{d}z}\nonumber\\
&=&\mathcal{A}(d,-\beta)\lim_{\varepsilon\to
0}\lim_{n\to\infty}F_{n}(\varepsilon),
\nonumber
\end{eqnarray}
where
\[
F_{n}(\varepsilon):=\sum_{k=0}^{n}\int_{|z|>\varepsilon}\big(G_{k}(x+z,y)
-G_{k}(x,y)\big)b(x,z)|z|^{-d-\beta}{\mathrm{d}z} \quad\mbox{ for }\varepsilon>0.
\]
 For any $n,m\in \mathbb{Z}_{+},\ n>m$
and any $\varepsilon>0$, we have by \eqref{16},
\begin{align}
&\big|F_{n}(\varepsilon)-F_{m}(\varepsilon)\big|\nonumber\\
=& \big| \sum_{k=m+1}^n \int_{|z|>\eps} \frac{G_k(x+z, y)-G_k (x, y)}{|z|^{d+\beta}}
b(x, z) {\mathrm{d}z} \big| \nonumber \\
=&
\big|\sum_{k=m+1}^{n}\int_{|z|>\varepsilon}\int_{B_{r}}\frac{G_{B_{r}}(x+z,u)
-G_{B_{r}}(x,u)}{|z|^{d+\beta}}\, \mathcal{S}^{b}_{u}G_{k-1}(u,y)\,b(x,z)\,{\mathrm{d}u}\,{\mathrm{d}z}\big|\nonumber\\
=&\bigg|\sum_{k=m+1}^{n}\int_{B_{r}}\mathcal{S}^{b}_{u}G_{k-1}(u,y)\mathrm{d}u\int_{|z|>\varepsilon}\frac{G_{B_{r}}(x+z,u)-G_{B_{r}}(x,u)-\nabla_{x}
G_{B_{r}}(x,u)\cdot z \1_{\{|z|<\lambda(u)\}}}{|z|^{d+\beta}}\,b(x,z){\mathrm{d}z}\bigg|.\nonumber
\end{align}
Here $\lambda(u):=\lambda_{B_{r}}(x,u)$ depends on the integration variable $u$, and the final equality holds because the compensator term integrates to zero by the symmetry \eqref{condi1}.
It then follows from \eqref{DeltaGn}, \eqref{lem7.main} and \eqref{lem13.2}
that
\begin{align}
&\big|F_{n}(\varepsilon)-F_{m}(\varepsilon)\big|\nonumber\\
\le&
 \sum_{k=m+1}^{n}  c_{2}A\delta^{k-1}
 \int_{B_{r}}\bigg(\int_{|z|<\lambda(u)}\frac{\big|G_{B_{r}}(x+z,u)-G_{B_{r}}(x,u)-\nabla_{x}
G_{B_{r}}(x,u)\cdot z\big|}{|z|^{d+\beta}}|b(x,z)|{\mathrm{d}z}\nonumber\\
&\qquad+\int_{|z|\ge
\lambda(u)}\frac{G_{B_{r}}(x+z,u)+G_{B_{r}}(x,u)}{|z|^{d+\beta}}|b(x,z)|{\mathrm{d}z}\bigg) \wh{h}_{B_{r}}(u,y){\mathrm{d}u}\nonumber\\
=&\sum_{k=m+1}^{n}c_{2}A\delta^{k-1}\int_{B_{r}}|\mathcal{S}^{b}_{x}|G_{B_{r}}(x,u)\wh{h}_{B_{r}}(u,y)\mathrm{d}u\nonumber\\
\le&
c_{3}\sum_{k=m+1}^{n}\delta^{k-1}\int_{B_{r}}\wh{h}_{B_{r}}(x,u)\wh{h}_{B_{r}}(u,y){\mathrm{d}u}\nonumber\\
\le&c_{3}C(r)\wh{h}_{B_{r}}(x,y)\sum_{k=m+1}^{n}\delta^{k-1}.\nonumber
\end{align}
Here $c_{3}=c_{3}(d,\alpha,\beta,A)>0$.
Therefore,
$\sup_{\varepsilon>0}\big|F_{n}(\varepsilon)-F_{m}(\varepsilon)\big|\to
0$ as $m,n\to\infty$. This implies that $\{F_{n}(\varepsilon):n\ge
1\}$ is
a uniformly convergent sequence of continuous functions. It
follows that
\begin{equation}\label{5.72}
\mathcal{S}^{b}_{x}G^{b}_{B_{r}}(x,z)=\lim_{\varepsilon\to
0}\lim_{n\to\infty}F_{n}(\varepsilon)=\lim_{n\to\infty}\lim_{\varepsilon\to
0}F_{n}(\varepsilon)=\lim_{n\to\infty}\sum_{k=0}^{n}
\mathcal{S}^{b}_{x}G_{k}(x,z)=\sum_{k=0}^{\infty}\mathcal{S}^{b}_{x}G_{k}(x,z).
\end{equation}
The second assertion in \eqref{lem13.main}
 thus follows from
\eqref{DeltaGn}, \eqref{5.72} and the fact that $c_{2}AC(r)\le \delta\le 1/(2c^{2}_{1}+1)$ for all $r\in (0,\tilde{r}_{1}]$.

We next
define $r_{1}$ and $C_{37}$, and verify
\eqref{lem13.main}--\eqref{r*scale}.
Set
\[
r_{1}(d,\allowbreak \alpha,\allowbreak \beta,\allowbreak A):=A^{-1/(\alpha-\beta)}\,\tilde{r}_{1}(d,\allowbreak \alpha,\allowbreak \beta,\allowbreak 1),\qquad
C_{37}(d,\allowbreak \alpha,\allowbreak \beta,\allowbreak A):=A\,\tilde{C}_{37}(d,\allowbreak \alpha,\allowbreak \beta,\allowbreak 1).
\]
Then \eqref{r*scale} holds by definition. Let
$b\in\mathcal{B}_{A}$ and
$B_r$ be a ball of radius $r\le r_{1}(d,\alpha,\beta,A)$, and put $\lambda:=A^{1/(\alpha-\beta)}$.
Then
$b_{\lambda}\in\mathcal{B}_{1}$,
and $\lambda B_r$ has radius at most $\tilde{r}_{1}(d,\alpha,\beta, 1)$.
Hence \eqref{lem13.main}--\eqref{e:4.66} hold for $b_{\lambda}$ on $\lambda B_r$, with $A=1$.
By \eqref{scalingforgb}, applied to $b$ and to $b\equiv 0$, the first assertion in \eqref{lem13.main} holds for $b$ on $B_r$.
A change of variables gives $
\mathcal{S}^{b}
[u(\lambda\,\cdot)](x)=\lambda^{\alpha}(\mathcal{S}^{b_{\lambda}}u)(\lambda x)$, and \eqref{e:hhat} gives $\wh{h}_{\lambda B_r}(\lambda x,\lambda y)=\lambda^{\alpha-\beta-d}\wh{h}_{B_{r}}(x,y)$.
Together with \eqref{scalingforgb}, these identities yield
\begin{align*}
|\mathcal{S}^{b}_{x}G^{b}_{B_{r}}(x,y)|=&\lambda^{d}\big|\big(\mathcal{S}^{b_{\lambda}}G^{b_{\lambda}}_{\lambda B_r}(\cdot,\lambda y)\big)(\lambda x)\big|\\
\le &\lambda^{d}\tilde{C}_{37}(d,\alpha,\beta,1)\,\wh{h}_{\lambda B_r}(\lambda x,\lambda y)=C_{37}(d,\alpha,\beta,A)\,\wh{h}_{B_{r}}(x,y).
\end{align*}
By the same change of variables and induction on $n$, the term $G_{n}(x,y)$ in \eqref{e:4.66} for $b$ on $B_r$ equals $\lambda^{d-\alpha}$ times the corresponding term for $b_{\lambda}$ on $\lambda B_r$ at $(\lambda x,\lambda y)$. Hence \eqref{e:4.65} also holds for $b$ on $B_r$.

It remains to prove \eqref{lem13.1} and \eqref{lem13.2}. Let $\gamma:=(\alpha-\beta)\wedge \alpha/2$ if $\alpha\not= 2\beta$ and $\gamma:=\beta/2$ if $\alpha=2\beta$. Lemma \ref{lemma8}(i) implies that there is a constant $c_{4}=c_{4}(d,\alpha,\beta)>0$ such that for every $r\in (0,1]$, every ball $B_{r}$ with radius $r$, and every pairwise distinct $x,y,z\in B_{r}$, $I_{B_{r}}(x,y,z)\le c_{4}\big(|x-z|^{-d+\gamma}+|y-z|^{-d+\gamma}\big)$, and hence
$$\int_{B_{r}}I_{B_{r}}(x,y,z)\mathrm{d}z\le 2c_{4}\int_{B_{r}}|x-z|^{-d+\gamma}\mathrm{d}z\le 2c_{4}\int_{B(x,2r)}|x-z|^{-d+\gamma}\mathrm{d}z\lesssim \int_{0}^{2r}s^{\gamma-1}\mathrm{d}s\asymp r^{\gamma}.$$
Therefore, \eqref{lem13.1} holds for all $0<\beta<\alpha<2$, with $C(r)\asymp r^{\gamma}$.
Similarly, \eqref{8.1} gives $J_{B_{r}}(x,y,z)\le C_{36}\big(|x-z|^{-d+\alpha-\beta}+|y-z|^{-d+\alpha-\beta}\big)$, and hence \eqref{lem13.2} holds with $C(r)\asymp r^{\alpha-\beta}$.
This completes the proof.
\end{proof}

The proof of Theorem \ref{them4.7} applies to every family of $C^{1,1}$ open sets whose rescalings to unit size have uniformly bounded characteristics. This gives the following analogue for the sets $V(Q,r)$ of Proposition \ref{P:2.2}(ii); compare \cite[Proposition 4.4]{CKS2} for the fractional Laplacian with gradient perturbation.

 \begin{thm}\label{them4.8}
Suppose $D$ is a
 $C^{1,1}$ open set in $\R^{d}$
 with characteristics $(R_{0},\Lambda_{0})$
 and $A\in (0,\infty)$. There
exists a positive constant $
r_{2}=r_{2}(d,\allowbreak \alpha,\allowbreak \beta,\allowbreak R_{0},\allowbreak \Lambda_{0},\allowbreak A)
\in (0,R_{0}\wedge 1)
$ such that for every
$b\in\mathcal{B}_{A}$,
$Q\in\partial D$ and $r\in (0,r_{2}]$,
 \begin{equation}\label{them4.8.main}
\frac{1}{2}\,G_{V(Q,r)}(x,y)\le G^{b}_{V(Q,r)}(x,y)\le \frac{3}{2}\,G_{V(Q,r)}(x,y),
\quad  x,y\in V(Q,r).
\end{equation}
Here  $V(Q, r)$ is
 the $C^{1,1}$ domain introduced in Proposition \ref{P:2.2}(ii).
 \end{thm}

\begin{proof}
Let $L=L(d,R_{0},\Lambda_{0})$ be the constant in Proposition \ref{P:2.2}(ii). For $r\in (0,R_{0}\wedge 1]$, the set $V:=V(Q,r)$ has characteristics $(rR_{0}/L,\Lambda_{0}L/r)$ and $\mathrm{diam}(V)\le 4r$. Hence $r^{-1}V$ has characteristics $(R_{0}/L,
\Lambda_{0}L)$ and diameter at most $4$.
The constants in \eqref{greenesti}, \eqref{e:3G} and Lemma \ref{lemma7} are
scale-invariant in the underlying set
and depend on the set only through its characteristics and its diameter, and the bounds in Lemma \ref{lemma8}(i) involve the set only through its diameter. All these constants are therefore bounded on the family $\{V(Q,r):Q\in\partial D,\ r\in(0,R_{0}\wedge 1]\}$ by constants depending only on $d,\alpha,\beta,R_{0},\Lambda_{0}$.
Moreover, for $x,y\in V$, $\int_{V}|x-z|^{-d+\gamma}\mathrm{d}z\le\int_{B(x,4r)}|x-z|^{-d+\gamma}\mathrm{d}z\lesssim r^{\gamma}$ for every $\gamma>0$.
With these facts, the proof of the first assertion of \eqref{lem13.main} applies verbatim with $B_{r}$ replaced by $V$, and gives \eqref{them4.8.main} for all $r\in(0,r_{2}]$ if $r_{2}\in (0,R_{0}\wedge 1)$ is sufficiently small.
\end{proof}

It follows from \eqref{levysysforxbd} that for every bounded open
set $D$ in $\R^{d}$, every $f\ge 0$, and $x\in D$,
\begin{equation}
\E_x\big[ f(X^{b}_{\tau^{b}_{D}});\,X^{b}_{\tau^{b}_{D}-}\not=X^{b}_{\tau^{b}_{D}}\big]
=\int_{\bar{D}^{c}}f(z)\Big(\int_{D}G^{b}_{D}(x,y)j^{b}(y,z){\mathrm{d}y}\Big){\mathrm{d}z}.\label{(*)}
\end{equation}
Define
\begin{equation}
K^{b}_{D}(x,z):=\int_{D}G^{b}_{D}(x,y)j^{b}(y,z){\mathrm{d}y},\quad
(x,z)\in D\times \bar{D}^{c}.
\nonumber
\end{equation}
Then \eqref{(*)} can be rewritten as
\begin{equation}
\E_x\big[f(X^{b}_{\tau^{b}_{D}});\,X^{b}_{\tau^{b}_{D}-}\not=X^{b}_{\tau^{b}_{D}}\big]=\int_{\bar{D}^{c}}f(z)K^{b}_{D}(x,z){\mathrm{d}z}.\label{(5)}
\end{equation}

\begin{lemma}\label{lem6.1.1}
Suppose $D$ is a $C^{1,1}$ open set with $\mathrm{diam}(D)\le
r_{3}=r_{3}(d,\allowbreak \alpha,\allowbreak \beta,\allowbreak A):= \frac{1}{4}\eta(A)\wedge 3r_{1}$.
Here $\eta(A)$ and $r_{1}$ are constants given in \eqref{defiforetaA} and Theorem \ref{them4.7}, respectively.
Then
\[
\P_x\big(X^{b}_{\tau^{b}_{D}}\in\partial D\big)=0,\quad x\in D.
\]
In this case, for every non-negative measurable function $f$,
\[
\E_x\big[f(X^{b}_{\tau^{b}_{D}})\big]=\int_{\bar{D}^{c}}f(z)K^{b}_{D}(x,z){\mathrm{d}z}\quad  x\in D.
\]
\end{lemma}

\begin{proof} Fix $x\in D$. Set $r:=\frac{1}{2}(\delta_{D}(x)\wedge r_{1})$.
Then $\frac{1}{2}\delta_{D}(x)\ge r\ge
\frac{1}{2}\big(\delta_{D}(x)\wedge\frac{1}{3}\mathrm{diam}(D)\big)\ge
\frac{1}{12}\delta_{D}(x)$. Let $B:=B(x,r)\subset D$. Let $r_{0}=r_{0}(R_{0},\Lambda_{0})$ be the constant introduced in Proposition \ref{P:2.2}(i), and put $\rho:=r\wedge (r_{0}/2)$.
Let $z_{x}\in\partial D$ satisfy $|x-z_{x}|=\delta_{D}(x)$, and let $B(x'',\rho)\subset D^{c}$ be the exterior ball tangent to $\partial D$ at $z_{x}$.
Set $B':=B(x'',\rho/2)$, so that $B'\subset\bar{D}^{c}$.
For $y\in B$ and $z\in B'$ we have $|y-z|\le r+\delta_{D}(x)+2\rho\le \frac{5}{2}\mathrm{diam}(D)<\eta(A)$, since $\rho\le r\le\delta_{D}(x)/2$.
For $z\in B'$ we also have $2r\le\delta_{D}(x)\le|x-z|\le\delta_{D}(x)+2\rho\le 14r$, since $\delta_{D}(x)\le 12r$. It follows from
\eqref{(5)}, Theorem \ref{them4.7}, \eqref{jb&j} and
\eqref{ballpoisson} that
\begin{eqnarray}
\P_x\big(X^{b}_{\tau^{b}_{B}}\in
B'\big)&=&\int_{B'}\int_{B}G^{b}_{B}(x,y)j^{b}(y,z){\mathrm{d}y\mathrm{d}z}\nonumber\\
&\ge&\frac{1}{4}\int_{B'}\int_{B}G_{B}(x,y)j(y,z){\mathrm{d}y\mathrm{d}z}\nonumber\\
&=&\frac{1}{4}\int_{B'}K_{B}(x,z){\mathrm{d}z}\nonumber\\
&\ge&c_{0}
\nonumber
\end{eqnarray}
for some constant $c_{0}=c_{0}(d,\alpha,R_{0},\Lambda_{0},\mathrm{diam}(D))>0$.
Indeed, by \eqref{ballpoisson} and the bounds on $|x-z|$, $K_{B}(x,z)\ge c_{1}(d,\allowbreak \alpha)r^{-d}$ for $z\in B'$.
Hence $\int_{B'}K_{B}(x,z){\mathrm{d}z}\ge c_{2}(d,\allowbreak \alpha)(\rho/r)^{d}\ge c_{2}\big(1\wedge \frac{r_{0}}{2\,\mathrm{diam}(D)}\big)^{d}$, uniformly in $x\in D$.
Let $D_{n}:=\{y\in D:\delta_{D}(y)>1/n\}$ for every
$n\in\mathbb{N}$. For $n$ sufficiently large, we have $B\subset
D_{n}$. In this case
\begin{eqnarray}
\P_x\big(X^{b}_{\tau^{b}_{D_{n}}}\in\bar{D}\big)&=&\P_x\big(X^{b}_{\tau^{b}_{B}}\in\bar{D}\setminus
D_{n}\big)+\P_x\big(X^{b}_{\tau^{b}_{B}}\in D_{n}\setminus
B,X^{b}_{\tau^{b}_{D_{n}}}\in\bar{D}\big)\nonumber\\
&\le&\P_x\big(X^{b}_{\tau^{b}_{B}}\in \bar{D}\setminus
B\big)
 \le 1-\P_x\big(X^{b}_{\tau^{b}_{B}}\in B'\big)\le 1-c_{0}.\nonumber
\end{eqnarray}
Let $u(x)=\P_x\big(X^{b}_{\tau^{b}_{D}}\in\partial D\big)$ and
$C:=\sup\{u(x):x\in D\}$. By the strong Markov property,
\begin{eqnarray}
u(x)&=&\E_x\big[u(X^{b}_{\tau^{b}_{D_{n}}});\,X^{b}_{\tau^{b}_{D_{n}}}\in\bar{D}\big]\nonumber\\
&=&\E_x\big[u(X^{b}_{\tau^{b}_{D_{n}}});\,X^{b}_{\tau^{b}_{D_{n}}}\in\partial
D\big]+\E_x\big[u(X^{b}_{\tau^{b}_{D_{n}}});\,X^{b}_{\tau^{b}_{D_{n}}}\in
D\big].\nonumber
\end{eqnarray}
Since $\P_x\big(X^{b}_{\tau^{b}_{D_{n}}}\in\partial D\big)=0$ by
\eqref{levysysforxbd}, we get
\[u(x)=\E_x\big[u(X^{b}_{\tau^{b}_{D_{n}}});\,X^{b}_{\tau^{b}_{D_{n}}}\in
D\big],\] and consequently $C\le(1-c_{0})C$. Thus $C=0$. \end{proof}

For the remainder of this paper, we assume $A\in (0,{\infty})$ and $b\in\mathcal{B}_{A}$. Let $C_{i}=C_{i}(d,\alpha,\beta,A,1)$, $i=18,19,20$
be the constants given in \eqref{ltdensity} with $\lambda=1$, and $q_{0}:=2C_{19}$.
Here we claim that \eqref{ltdensity} also holds with $p_{A}$ in place of $p_{M_{b^{+},1}}$, after enlarging $C_{18}$ and $C_{19}$ if necessary. Indeed, since $M_{b^{+},1}\le A$ and $(At)^{-d/\beta}\ge (1\wedge A^{-d/\beta})t^{-d/\alpha}$ for $t\in (0,1]$, \eqref{pa} gives $p_{M_{b^{+},1}}(t,x,y)\le c\,p_{A}(t,x,y)$ for $t\in (0,1]$ with $c=c(d,\alpha,\beta,A)$. The claim then follows from \eqref{stdensity} and the Chapman--Kolmogorov equation.
It follows that for every $(t,x,y)\in (0,{\infty})\times \R^{d}\times \R^{d}$,
\begin{equation}
C_{18}^{-1}\e^{-3C_{19}t}\bar{p}_{1}(t,C_{20}x,C_{20}y)\le p^{b,q_{0}}(t,x,y)\le C_{18}\e^{-C_{19}t}p_{A}(t,x,y).\label{n2.1}
\end{equation}
By this and
\eqref{pa} one verifies that for any open set $D$ in $\R^{d}$, the Green function $G^{b,q_{0}}_{D}(x,y):=\int_{0}^{{\infty}}p^{b,q_{0}}_{D}(t,x,y){\mathrm{d}t}$ of $X^{b,q_{0},D}$ is finite and continuous off the diagonal of $D\times D$, and
\begin{eqnarray}
G^{b,q_{0}}_{D}(x,y)\le G^{b,q_{0}}(x,y)\le c_{1}\big(|x-y|^{-d+\alpha}+|x-y|^{-d+\beta}\big)\label{n2.2}
\end{eqnarray}
for some positive constant $c_{1}=c_{1}(d,\allowbreak \alpha,\allowbreak \beta,\allowbreak A)$.
Moreover, if $D$ is a bounded open set in $\R^{d}$, it follows
from \eqref{3.12} that
\begin{equation}\label{n2.5}
G^{b,q_{0}}_{D}(x,y)\le G^{b}_{D}(x,y)\le c_{2}|x-y|^{\alpha -d }\quad\mbox{ for all }x,y\in D.
\end{equation}
Here $c_{2}$ is a positive constant depending on $d,\alpha,\beta,D$ and $A$.

\begin{lemma}\label{nlem2.1}
There is a positive constant $r_{4}=r_{4}(d,\allowbreak \alpha,\allowbreak \beta,\allowbreak A)$ such that
for any $b\in\mathcal{B}_{A}$
and any ball
$B$
with radius $r\in (0,r_{4}]$,
\begin{equation*}
\frac{1}{4}G_{B}(x,y)\le G^{b,q_{0}}_{B}(x,y)\le \frac{3}{2}G_{B}(x,y)\quad\mbox{ for }x,y\in B.
\end{equation*}
\end{lemma}

\begin{proof}
The proof is similar
to that of \cite[Theorem 5.1]{CS2}.
For any $z\in B$, let $(X^{b,B}_{t},\mathbb{P}^{z}_{x})$ be the $G^{b}_{B}(\cdot,z)$-transform of $(X^{b,B}_{t},\P_{x})$.
Indeed,  by \cite[Lemma 3.5(1)]{Chen1},
\begin{equation}\label{n2.6}
 G^{b,q_{0}}_{B}(x,z)
 =G^{b}_{B}(x,z)\mathbb{E}^{z}_{x}\big[ \e^{-q_{0}\tau^{b}_{B}}\big] \quad\mbox{ for }x,z\in B.
 \end{equation}
Let  $r_{4}\le r_{1}$, where $r_{1}$ is the constant given in Theorem \ref{them4.7}.
It follows from \eqref{lem13.main} and \eqref{e:3G} that
there is a positive constant $c_1= c_{1}(d,\allowbreak \alpha) >0$ so that  for any ball $B$ with radius $r\in (0, r_4]$,
\begin{equation*}
\frac{G^{b}_{B}(x,y)G^{b}_{B}(y,z)}{G^{b}_{B}(x,z)}\le c_{1}\big(|x-y|^{\alpha-d}+|y-z|^{\alpha-d}\big)\quad\mbox{ for }x,y,z\in B.
\end{equation*}
Taking $r_4$ smaller if needed, we have for any ball $B$ with radius $r\in (0, r_4]$,
\[
\sup_{x, z\in B}  \mathbb{E}^{z}_{x}\big[ \tau^{b}_{B}\big] =
 \sup_{x, z\in B} \int_{B}\frac{G^{b}_{B}(x,y)G^{b}_{B}(y,z)}{G^{b}_{B}(x,z)}\,{\mathrm{d}y}\le q^{-1}_{0}\ln 2 .
 \]
Thus by  Jensen's inequality,  we have
$ 1/2 \leq \mathbb{E}^{z}_{x}\big[ \e^{-q_{0}\tau^{b}_{B}}\big]  \leq 1$ for   any $ x, z \in B$ with $x\not=z$.
The desired conclusion now follows from \eqref{n2.6} and \eqref{lem13.main}.
 \end{proof}

\begin{lemma}\label{nlem2.2}
Suppose
$D$ is a $C^{1,1}$ open set with the characteristics $(R_{0},\Lambda_{0})$.
There is a positive constant
$r_{5}=r_{5}(d,\allowbreak \alpha,\allowbreak \beta,\allowbreak  R_0,\allowbreak  \Lambda_0,\allowbreak  A)$
so that for
every $b\in\mathcal{B}_{A}$,
every $Q\in \partial D$ and $r\in (0,r_{5}]$,
\begin{equation}\label{n2.7}
\frac{1}{4}G_{V(Q,r)}(x,y)\le G^{b,q_{0}}_{V(Q,r)}(x,y)\le \frac{3}{2}G_{V(Q,r)}(x,y)\quad\mbox{ for }x,y\in V(Q,r).
\end{equation}
Here $V(Q, r)$ is
the $C^{1,1}$ domain introduced
in Proposition \ref{P:2.2}(ii).
\end{lemma}

\begin{proof}
 Using \eqref{them4.8.main} instead of \eqref{lem13.main}, the proof of this lemma is the same as that of Lemma \ref{nlem2.1}.
\end{proof}

\section{Duality}\label{sec6}
The operator $\L^b$ is not symmetric. Upper bounds on $p^{b}_{D}(t,x,y)$ in the second variable therefore cannot be obtained from those in the first variable by symmetry.
This section constructs a dual process that plays the role of the adjoint.
Throughout,
$A\in (0,+\infty)$ and $b\in\mathcal{B}_{A}$.
It follows from Corollary \ref{ncor1.1} and Corollary \ref{ncor1.3} that the process $X^{b,q_{0}}$ has a jointly continuous and strictly positive transition density $p^{b,q_{0}}(t,x,y)$. Define
\[h^{b,q_{0}}(x):=\int_{\R^{d}}G^{b,q_{0}}(y,x){\mathrm{d}y}\quad\mbox{ and }\xi^{b,q_{0}}({\mathrm{d}x}):=h^{b,q_{0}}(x){\mathrm{d}x}.\]

\begin{proposition}\label{nprop3.1}
The function $h^{b,q_{0}}(x)$ is continuous on $\R^{d}$, and there is a positive constant $ C_{38}= C_{38}(d,\allowbreak \alpha,\allowbreak \beta,\allowbreak A)\ge 1$ such that $ C_{38}^{-1}\le h^{b,q_{0}}(x)\le  C_{38}$ for all $x\in \R^{d}$. Moreover, $\xi^{b,q_{0}}({\mathrm{d}x})$ is an excessive measure for $X^{b,q_{0}}$, that is, for any non-negative Borel function $f$ and $t\ge 0$,
\[\int_{\R^{d}}P^{b,q_{0}}_{t}f(x)\xi^{b,q_{0}}({\mathrm{d}x})\le \int_{\R^{d}}f(x)\xi^{b,q_{0}}({\mathrm{d}x}).\]
\end{proposition}
\begin{proof} The proof of the second claim is essentially the same as that of \cite[Proposition 5.2]{CKS2}. We prove only the first claim. By \eqref{n2.1}, we have
\begin{equation}
h^{b,q_{0}}(x)=\int_{\R^{d}}\int_{0}^{{\infty}}p^{b,q_{0}}(t,y,x)\mathrm{d}t\mathrm{d}y
\le\int_{\R^{d}}\int_{0}^{{\infty}}C_{18} \e^{-C_{19}t}p_{A}(t,y,x)\mathrm{d}t\mathrm{d}y
=\int_{0}^{{\infty}}C_{18} \e^{-C_{19}t}{\mathrm{d}t}<{\infty},\nonumber
\end{equation}
\begin{equation}
h^{b,q_{0}}(x)
\ge\int_{\R^{d}}\int_{0}^{{\infty}}C_{18}^{-1} \e^{-3C_{19}t}\bar{p}_{1}(t,C_{20}y,C_{20}x)\mathrm{d}t\mathrm{d}y
=C_{20}^{-d}\int_{0}^{{\infty}}C_{18}^{-1} \e^{-3C_{19}t}{\mathrm{d}t}>0.\nonumber
\end{equation}
The continuity of $h^{b,q_{0}}$ follows from the continuity of $G^{b,q_{0}}$ off the diagonal, the bound \eqref{n2.2}, which makes $y\mapsto G^{b,q_{0}}(y,x)$ uniformly integrable near $y=x$ for $x$ in compact sets, and the decay $\int_{0}^{\infty}\e^{-C_{19}t}p_{A}(t,y,x)\,{\mathrm{d}t}\to 0$ as $|y-x|\to\infty$, uniformly in $x$.\end{proof}

The transition density of the process $X^{b,q_{0}}$ with respect to the reference measure $\xi^{b,q_{0}}$ is defined by
\[\bar{p}^{b,q_{0}}(t,x,y):=\frac{p^{b,q_{0}}(t,x,y)}{h^{b,q_{0}}(y)},\quad
 (t,x,y)\in(0,\infty)\times \R^{d}\times \R^{d}.\]
Then
\[\bar{G}^{b,q_{0}}(x,y):=\int_{0}^{\infty}\bar{p}^{b,q_{0}}(t,x,y){\mathrm{d}t}=\frac{G^{b,q_{0}}(x,y)}{h^{b,q_{0}}(y)},
\quad  x,y\in \R^{d}\]
is the Green function of $X^{b,q_{0}}$ with respect to $\xi^{b,q_{0}}$. One checks that $\bar{G}^{b,q_{0}}(x,y)$ has the following properties:
\begin{enumerate}[label=\upshape(A\arabic*), ref=\upshape(A\arabic*)]
\item $\bar{G}^{b,q_{0}}(x,y)>0$ on $\R^{d}\times \R^{d}$, and $\bar{G}^{b,q_{0}}(x,y)=\infty$ if and only if $x=y$.
\item For every $x\in \R^{d}$, $\bar{G}^{b,q_{0}}(x,\cdot)$ and $\bar{G}^{b,q_{0}}(\cdot,x)$ are
continuous on $\R^d$ in the extended sense.
\item For every compact set $K\subset \R^{d}$ and $x\in \R^{d}$, $\int_{K}\bar{G}^{b,q_{0}}(x,y)\xi^{b,q_{0}}({\mathrm{d}y})<\infty$.

 \item For any compact set $K\subset \R^{d}$ and $y\in \R^{d}$, $\int_{K}\bar{G}^{b,q_{0}}(x,y)\xi^{b,q_{0}}({\mathrm{d}x})<\infty$.

\end{enumerate}
Indeed, (A1) follows from Corollary \ref{ncor1.3}, Proposition \ref{nprop3.1} and \eqref{n2.1}. (A2) follows from the continuity of $\bar{p}^{b,q_{0}}(t,x,y)$ and \eqref{n2.1}. (A3) and (A4) follow from \eqref{n2.2}.
Properties (A1)--(A3) imply that the process $X^{b,q_{0}}$ satisfies the conditions (R) of \cite[p.211]{Chung} and conditions (a) and (b) of \cite[Theorem 5.4]{Chung}. Thus it satisfies Hunt's Hypothesis (B) by \cite[Theorem 5.4]{Chung}. It follows from \cite[Theorem 13.24]{Chung} that $X^{b,q_{0}}$ has a dual process $\wh{X}^{b,q_{0}}$ with respect to the reference measure $\xi^{b,q_{0}}$, and $\wh{X}^{b,q_{0}}$ is a standard process.
Moreover, $\bar{G}^{b,q_{0}}(x,y)$ also satisfies the following property.

\begin{enumerate}[label=\upshape(A\arabic*), ref=\upshape(A\arabic*), start=5]
\item For every $y\in \R^{d}$, $x\mapsto\bar{G}^{b,q_{0}}(x,y)$ is an excessive
function with respect to $X^{b,q_{0}}$, that is, for every $t>0$ and
$x\in \R^{d}$,
\[\E_x(\bar{G}^{b,q_{0}}(X^{b,q_{0}}_{t},y))\le \bar{G}^{b,q_{0}}(x,y),\mbox{ and }\lim_{t\downarrow 0}\E_x(\bar{G}^{b,q_{0}}(X^{b,q_{0}}_{t},y))=\bar{G}^{b,q_{0}}(x,y).\]
For every $y\in \R^{d}$, $x\mapsto\bar{G}^{b,q_{0}}(x,y)$ is harmonic with respect to $X^{b,q_{0}}$ in $\R^{d}\setminus \{y\}$. Furthermore, for every open set $U\subset \R^{d}$, we have
\[
\E_x \big[\bar{G}^{b,q_{0}}(X^{b,q_{0}}_{T^{b}_{U}},y) \big]=\bar{G}^{b,q_{0}}(x,y)
\quad \hbox{for } (x,y)\in \R^{d}\times U,
\]
where $T^{b}_{U}:=\inf\{t>0:X^{b,q_{0}}_{t}\in U\}$. In particular, for
every $y\in \R^{d}$ and $\eps >0$, $x\mapsto \bar{G}^{b,q_{0}}(x,y)$ is
regularly harmonic on $\R^{d}\setminus B(y,\eps )$ with respect to
$X^{b,q_{0}}$.
\end{enumerate}

\begin{proof}[Proof of (A5)]
By a standard argument (see \cite[Proof of
(A4)]{CKS2} and the references therein), it suffices to show
that for every $x\in U^{c}$,
$y\mapsto \E_{x}\big[\bar{G}^{b,q_{0}}\big(X^{b,q_{0}}_{T^{b}_{U}},y\big)\big]$ is continuous in
$U$. Fix $x\in U^{c}$ and $y\in U$. Let $r:=\delta_{U}(y)$.
For any $\wh{y}\in B(y,r/4)$ and $\delta\in (0,r/2)$, by \eqref{levysysforxbq}, \eqref{n2.2} and Proposition \ref{nprop3.1},
we have
\begin{align}
&\E_{x}\big[\bar{G}^{b,q_{0}}\big(X^{b,q_{0}}_{T^{b}_{U}},\hat{y}\big);X^{b,q_{0}}_{T^{b}_{U}}\in B(y,\delta)\big]\nonumber\\
=&\int_{B(y,\delta)}\bar{G}^{b,q_{0}}(z,\hat{y}){\mathrm{d}z}\int_{U^{c}}G^{b,q_{0}}_{U^{c}}(x,w)j^{b}(w,z){\mathrm{d}w}\nonumber\\
=&\int_{B(y,\delta)}\frac{G^{b,q_{0}}(z,\hat{y})}{h^{b,q_{0}}(\hat{y})}{\mathrm{d}z}\int_{U^{c}}G^{b,q_{0}}_{U^{c}}(x,w)j^{b}(w,z){\mathrm{d}w}\nonumber\\
\le&c_{2}\int_{B(y,\delta)}\big(|z-\hat{y}|^{-d+\alpha}+|z-\hat{y}|^{-d+\beta}\big){\mathrm{d}z}
\int_{U^{c}}G^{b,q_{0}}_{U^{c}}(x,w)\big(|z-w|^{-d-\alpha}+|z-w|^{-d-\beta}\big){\mathrm{d}w}\nonumber\\
\le&c_{3}\big(r^{-d-\alpha}+r^{-d-\beta}\big)\int_{B(y,\delta)}\big(|z-\hat{y}|^{-d+\alpha}+|z-\hat{y}|^{-d+\beta}\big){\mathrm{d}z}
\int_{U^{c}}G^{b,q_{0}}_{U^{c}}(x,w){\mathrm{d}w}.\label{n3.1}
\end{align}
By \eqref{n2.1}, $\int_{U^{c}}G^{b,q_{0}}_{U^{c}}(x,w){\mathrm{d}w}=\int_{U^{c}}\int_{0}^{{\infty}}p^{b,q_{0}}(t,x,w){\mathrm{d}t\mathrm{d}w}\le\int_{0}^{{\infty}}C \e^{-C_{19}t}{\mathrm{d}t}<{\infty}$. Thus it follows from \eqref{n3.1} that for any $\eps >0$, there is $\delta\in (0,r/2)$ sufficiently small so that
\[\sup_{\hat{y}\in B(y,r/4)}\E_{x}\Big[\bar{G}^{b,q_{0}}\big(X^{b,q_{0}}_{T^{b}_{U}},\hat{y}\big);X^{b,q_{0}}_{T^{b}_{U}}\in B(y,\delta)\Big]<\frac{\eps }{8}.\]
Fix a sequence $\{y_{n}:\ n\ge 1\}\subset B(y,r/4)$ such that $y_{n}\to y$ as $n\to{\infty}$. Since $\bar{G}^{b,q_{0}}=G^{b,q_{0}}/h^{b,q_{0}}$ is bounded and jointly continuous on $B(y,\delta)^{c}\times B(y,\delta/2)$, by the bounded convergence theorem, we have
\[\lim_{n\to{\infty}}\big|\E_{x}\big[\bar{G}^{b,q_{0}}(X^{b,q_{0}}_{T^{b}_{U}},y_{n})-\bar{G}^{b,q_{0}}(X^{b,q_{0}}_{T^{b}_{U}},y)
;X^{b,q_{0}}_{T^{b}_{U}}\not\in B(y,\delta)\big]\big|=0.\]
Therefore, for $n$ sufficiently large,
\begin{eqnarray}
&&\big|\E_x\big[\bar{G}^{b,q_{0}}(X^{b,q_{0}}_{T^{b}_{U}},y_{n})\big]-\E_x\big[\bar{G}^{b,q_{0}}(X^{b,q_{0}}_{T^{b}_{U}},y)\big]\big|\nonumber\\
&\le&\E_x\big[\bar{G}^{b,q_{0}}(X^{b,q_{0}}_{T^{b}_{U}},y_{n});\,X^{b,q_{0}}_{T^{b}_{U}}\in B(y,\delta)\big]+\E_x\big[\bar{G}^{b,q_{0}}(X^{b,q_{0}}_{T^{b}_{U}},y);\,X^{b,q_{0}}_{T^{b}_{U}}\in B(y,\delta)\big]\nonumber\\
&&+\big|\E_x\big[\bar{G}^{b,q_{0}}(X^{b,q_{0}}_{T^{b}_{U}},y_{n})-\bar{G}^{b,q_{0}}(X^{b,q_{0}}_{T^{b}_{U}},y);\,X^{b,q_{0}}_{T^{b}_{U}}\not\in B(y,\delta)\big]\big|\nonumber\\
&<&\eps /2.\nonumber
\end{eqnarray}
This completes the proof.  \end{proof}

\begin{thm}\label{them5.4}
For every increasing sequence $\{U_{n}:n\ge 1\}$ of open sets with
$\overline{U}_{n}\subset U_{n+1}$ and $U_{n}\uparrow \R^{d}$, $\lim_{n\to
\infty}\E_x \big[ \bar{G}^{b,q_{0}}(X^{b,q_{0}}_{\tau^{b}_{U_{n}}},y) \big]=0$ for
every $x,y\in \R^{d}$ with $x\not=y$. Moreover, for every $x,y\in \R^{d}$,
$\lim_{t\to \infty}\E_x \big[ \bar{G}^{b,q_{0}}(X^{b,q_{0}}_{t},y)\big]=0$.
\end{thm}

 The proof of the above theorem is essentially the same as that of
\cite[Theorem 5.4]{CKS2},
and is omitted.
Using (A1)--(A5), \eqref{n2.5}
and Theorem \ref{them5.4}, we get from
 \cite[Theorem 7 and the remark following it on p. 391]{Liao2}
that the dual process $\wh{X}^{b,q_{0}}$ is a
 Hunt process. Its transience follows from (A2), (A4) and \cite[p. 368, Remark 1]{Liao2} applied to the dual potential kernel.
Let $\wh{P}^{b,q_{0}}_{t}$ denote the semigroup
of $\wh{X}^{b,q_{0}}$. Then for every $f,g\in L^{2}(\R^{d};\xi^{b,q_{0}})$,
\begin{equation}
\int_{\R^{d}}f(x)P^{b,q_{0}}_{t}g(x)\xi^{b,q_{0}}({\mathrm{d}x})=\int_{\R^{d}}\wh{P}^{b,q_{0}}_{t}f(x)g(x)\xi^{b,q_{0}}({\mathrm{d}x}).\label{duality}
\end{equation}
Define $\bar{H}^{q_{0}}_{t}:=t$ and
\[\bar{N}^{q_{0}}(x,{\mathrm{d}y}):=\frac{j^{b}(x,y)}{h^{b,q_{0}}(y)}\,\xi^{b,q_{0}}({\mathrm{d}y})\quad\mbox{ for } (x,y)\in \R^{d}\times \R^{d},\quad\bar{N}^{q_{0}}(x,\{\partial\}):=q_{0}.\]
Then $(\bar{N}^{q_{0}},\bar{H}^{q_{0}})$ is a L\'{e}vy system for $X^{b,q_{0}}$
with respect to the reference measure $\xi^{b,q_{0}}$. It follows from
 \cite[Theorem 3.5]{Getoor}
that a L\'{e}vy system $(\wh{N}^{q_{0}},\wh{H}^{q_{0}})$ for $\wh{X}^{b,q_{0}}$
 can be chosen so that
$\wh{H}^{q_{0}}_{t}=t$ and
$\wh{N}^{q_{0}}(y,{\mathrm{d}x})\xi^{b,q_{0}}({\mathrm{d}y})=\bar{N}^{q_{0}}(x,{\mathrm{d}y})\xi^{b,q_{0}}({\mathrm{d}x})$.
Therefore,
\[\wh{N}^{q_{0}}(x,{\mathrm{d}y})=\frac{j^{b}(y,x)}{h^{b,q_{0}}(x)}\,\xi^{b,q_{0}} ({\mathrm{d}y})=\frac{j^{b}(y,x)h^{b,q_{0}}(y)}{h^{b,q_{0}}(x)}\,{\mathrm{d}y}\quad\mbox{ for } (x,y)\in \R^{d}\times \R^{d}.
\]
Set $\wh{N}^{q_{0}}(x,\{\partial\}):=q_{0}$.
Consequently, for any $x\in \R^{d}$ and any non-negative measurable function $f$ on $\R_{+}\times \R^{d}\times
\R^{d}_{\partial}
$ vanishing on $\{(s,y,y):\ s\ge 0,y\in \R^{d}\}$ and stopping time $T$ (with respect to the filtration of $\wh{X}^{b,q_{0}}$), we have
 \begin{equation}\label{levysysforwhxbq}
\E_{x}\Big[\sum_{s\le T}f(s,\wh{X}^{b,q_{0}}_{s-},\wh{X}^{b,q_{0}}_{s})\Big]
=\E_{x}\Big[\int_{0}^{T}\int_{
\R^{d}_{\partial}
}f(s,\wh{X}^{b,q_{0}}_{s},y)\wh{N}^{q_{0}}(\wh{X}^{b,q_{0}}_{s},{\mathrm{d}y}){\mathrm{d}s}\Big].
\end{equation}

For any open subset $U$ of $\R^{d}$, we use $\wh{X}^{b,q_{0},U}$ to denote the
subprocess of $\wh{X}^{b,q_{0}}$ killed upon leaving $U$. Then $X^{b,q_{0},U}$ and
$\wh{X}^{b,q_{0},U}$ are dual processes with respect to the reference measure $\xi^{b,q_{0}}$.
By Corollary \ref{ncor1.1}, Proposition \ref{nprop3.1} and the duality relation \eqref{duality}, we have the following
theorem.
\begin{thm}\label{them5.6}
For any open subset $U$ in $\R^{d}$,
\begin{equation}\label{n3.4}
\wh{p}^{b,q_{0}}_{U}(t,x,y):=\frac{p^{b,q_{0}}_{U}(t,y,x)h^{b,q_{0}}(y)}{h^{b,q_{0}}(x)}
\end{equation}
is jointly continuous on $(0,\infty)\times U\times U$, and it is
the transition density of $\wh{X}^{b,q_{0},U}$ with respect to the Lebesgue
measure. Moreover,
\begin{equation}
\wh{G}^{b,q_{0}}_{U}(x,y):=\int_{0}^{\infty}\wh{p}^{b,q_{0}}_{U}(t,x,y){\mathrm{d}t}=\frac{G^{b,q_{0}}_{U}(y,x)h^{b,q_{0}}(y)}{h^{b,q_{0}}(x)}\quad\mbox{ for }  (x,y)\in U\times U\label{n3.5}
\end{equation}
is the Green function of $\wh{X}^{b,q_{0},U}$ with respect to the Lebesgue measure.
\end{thm}

\section{Small-time heat kernel estimates: the case $d>\alpha$}\label{sec7}

This section proves the small-time estimates of Theorem \ref{T:1.2}(i)--(iii) for $d>\alpha$.
The upper bounds combine the boundary decay of exit probabilities with the duality of Section \ref{sec6}.
The lower bounds follow from the Chapman--Kolmogorov equation and a chain argument.
Throughout,
$b\in\mathcal{B}_{A}$ for some $A\in (0,+\infty)$. Define
\begin{equation}\label{defiforM}
M=M(A,q_{0})
:=\sup_{b\in\mathcal{B}_{A}}\sup_{x,y\in\R^{d}}\frac{h^{b,q_{0}}(x)}{h^{b,q_{0}}(y)}.
\end{equation}
By Proposition \ref{nprop3.1}, $1\le M<{\infty}$.

\subsection{Small-time upper bound estimates}

For an open subset $U$ of $\R^{d}$, let
$\wh{\tau}^{b,q_{0}}_{U}:=\inf\{t>0:\wh{X}^{b,q_{0}}_{t}\not\in U\}$. The
proof of the following lemma is essentially the same as that of \cite[Lemma
7.3]{CKS2}, and is omitted.

\begin{lemma}\label{lem7.3}
Suppose $U$ is an open subset of $\R^{d}$. $U_{1},U_{3}$ are
open subsets of $U$ with $\mathrm{dist}(U_{1},U_{3})>0$ and
$U_{2}=U\setminus (U_{1}\cup U_{3})$. Then for any $x\in U_{1}$,
$y\in U_{3}$ and $t>0$, we have
\begin{align}
p^{b,q_{0}}_{U}(t,x,y)\le&\P_x\big(X^{b,q_{0}}_{\tau^{b,q_{0}}_{U_{1}}}\in U_{2}\big)\sup_{\substack{s<t
\\ z\in U_{2}}}p^{b,q_{0}}_{U}(s,z,y)+\big(t\wedge
\E_x(\tau^{b,q_{0}}_{U_{1}})\big)\esup_{\substack{u\in U_{1}
\\ z\in
U_{3}}}j^{b}(u,z),\nonumber\\
p^{b,q_{0}}_{U}(t,y,x)\le& M \P_x\big(\wh{X}^{b,q_{0}}_{\wh{\tau}^{b,q_{0}}_{U_{1}}}\in
U_{2}\big)\sup_{\substack{s<t\\ z\in U_{2}}}p^{b,q_{0}}_{U}(s,y,z)+M\big(t\wedge
\E_x(\wh{\tau}^{b,q_{0}}_{U_{1}})\big)\esup_{\substack{u\in U_{3}\\ z\in U_{1}}}j^{b}(u,z),\label{lem7.3.2}\\
p^{b,q_{0}}_{U}(t,x,y)\ge&
\frac{t}{M}\P_x\big(\tau^{b,q_{0}}_{U_{1}}>t\big)\P_y\big(\wh{\tau}^{b,q_{0}}_{U_{3}}>t\big)\einf_{\substack{u\in
U_{1}\\ z\in U_{3}}}j^{b}(u,z). \label{lem7.3.3}
\end{align}
\end{lemma}

\begin{lemma}\label{nlem4.2}
Suppose $D$ is a $C^{1,1}$ open set with $\mathrm{diam}(D)\le
r_{6}
=r_{6}(d,\allowbreak \alpha,\allowbreak \beta,\allowbreak A)
:=\frac{1}{4}\eta(A)\wedge 3r_{4}$, where
$\eta(A)$ and $r_{4}$ are the constants given in \eqref{defiforetaA} and Lemma \ref{nlem2.1}, respectively.
Then for every $x\in D$,
\[\P_{x}\big(X^{b,q_{0}}_{\tau^{b,q_{0}}_{D}}\in\partial D\big)=0 \quad \hbox{and} \quad
\P_{x}\big(\wh{X}^{b,q_{0}}_{\hat{\tau}^{b,q_{0}}_{D}}\in\partial D\big)=0.\]
\end{lemma}

\begin{proof} In view of \eqref{levysysforxbq}, \eqref{levysysforwhxbq}, Lemma \ref{nlem2.1}, \eqref{n3.5} and the fact that $M^{-1}\le
h^{b,q_{0}}(x)/h^{b,q_{0}}(y)\le M$ for all $x,y\in\R^{d}$, the proof of this lemma is essentially the same
as that of Lemma \ref{lem6.1.1}. We omit the details.\end{proof}

For any open set $D$, define the Poisson kernel for $X^{b,q_{0}}$ of $D$ by
\begin{equation}
K^{b,q_{0}}_{D}(x,z):=\int_{D}G^{b,q_{0}}_{D}(x,y)j^{b}(y,z){\mathrm{d}y},\quad
(x,z)\in D\times \bar{D}^{c},\nonumber
\end{equation}
and the Poisson kernel for $\wh{X}^{b,q_{0}}$ of $D$ by
\begin{equation}
\wh{K}^{b,q_{0}}_{D}(x,z):=\int_{D}\wh{G}^{b,q_{0}}_{D}(x,y)\frac{j^{b}(z,y)h^{b,q_{0}}(z)}{h^{b,q_{0}}(y)}\,{\mathrm{d}y},
\quad
(x,z)\in D\times \bar{D}^{c}.\nonumber
\end{equation}
By \eqref{levysysforxbq} and \eqref{levysysforwhxbq}, for every non-negative measurable function $f$ and $x\in D$,
\begin{equation}
\E_x\big[f(X^{b,q_{0}}_{\tau^{b,q_{0}}_{D}});\,X^{b,q_{0}}_{\tau^{b,q_{0}}_{D}-}\not=X^{b,q_{0}}_{\tau^{b,q_{0}}_{D}}\big]
=\int_{\bar{D}^{c}}f(z)K^{b,q_{0}}_{D}(x,z){\mathrm{d}z},\nonumber
\end{equation}
and
\begin{equation}
\E_x\big[f(\wh{X}^{b,q_{0}}_{\wh{\tau}^{b,q_{0}}_{D}});\,\wh{X}^{b,q_{0}}_{\wh{\tau}^{b,q_{0}}_{D}-}\not=\wh{X}^{b,q_{0}}_{\wh{\tau}^{b,q_{0}}_{D}}\big]
=\int_{\bar{D}^{c}}f(z)\wh{K}^{b,q_{0}}_{D}(x,z){\mathrm{d}z}.\nonumber
\end{equation}
In particular, by Lemma \ref{nlem4.2}, if $D$ is a $C^{1,1}$ open set with $\mathrm{diam}(D)\le r_{6}$, then
\begin{equation*}
\E_x\big[f(X^{b,q_{0}}_{\tau^{b,q_{0}}_{D}})\big]
=\int_{\bar{D}^{c}}f(z)K^{b,q_{0}}_{D}(x,z){\mathrm{d}z},
\end{equation*}
and
\begin{equation}
\E_x\big[f(\wh{X}^{b,q_{0}}_{\wh{\tau}^{b,q_{0}}_{D}})\big]
=\int_{\bar{D}^{c}}f(z)\wh{K}^{b,q_{0}}_{D}(x,z){\mathrm{d}z}.\label{poissonforwhxbq}
\end{equation}

For the rest of this section, let $D$ be a $C^{1,1}$ open set with characteristics $(R_{0},\Lambda_{0})$, and let $r_{0}=r_{0}(R_{0},\Lambda_{0})$ be the constant introduced in Proposition \ref{P:2.2}(i). Set
\begin{equation}\label{e:rstar}
r^{*}=r^{*}(d,\alpha,\beta,R_{0},\Lambda_{0},A):=1\wedge \min\{r_{i}:i=0,1,\cdots,6\}\quad\mbox{and}\quad T_{0}:=(r^{*}/16)^{\alpha}.
\end{equation}
Then
$r^{*}<R_{0}$, and $r^{*}\le r_{3}\le\eta(A)/4$.

\begin{lemma}\label{lem7.2}
Suppose $D$ is a $C^{1,1}$ open set
with characteristics $(R_{0},\Lambda_{0})$.
For every $T>0$, there exists a positive constant
$C_{39}= C_{39}(d,\allowbreak \alpha,\allowbreak \beta,\allowbreak R_{0},\allowbreak \Lambda_{0},\allowbreak A,\allowbreak M,\allowbreak T)$
such that for all $x\in D$ and $t\in (0,T]$,
\begin{equation}
\P_x\big(\tau^{b,q_{0}}_{D}>{t/4} \big)\le
 C_{39}\big(1\wedge\frac{\delta_{D}(x)^{\alpha/2}}{\sqrt{t}}\big),\label{lem7.2.1}
\end{equation}
\begin{equation}
\P_x\big(\wh{\tau}^{b,q_{0}}_{D}>{t/4} \big)\le  C_{39}
\big(1\wedge\frac{\delta_{D}(x)^{\alpha/2}}{\sqrt{t}}\big).\label{lem7.2.2}
\end{equation}
\end{lemma}

\begin{proof}
We prove only \eqref{lem7.2.2}. The proof of \eqref{lem7.2.1} is similar.
Fix $x\in D$, $T>0$ and $t\in (0,T]$.
Let $r^{*}$ be the constant given in \eqref{e:rstar} and $r:=\big(\frac{t}{T}\big)^{1/\alpha}r^{*}$. If $\delta_{D}(x)\ge r/16$, then
$ 1\wedge\frac{\delta_{D}(x)^{\alpha/2}}{\sqrt{t}} \ge   1\wedge\frac{ (r^*/16)^{\alpha/2}}{\sqrt{T}}$ and
\begin{equation}
\P_{x}\big(\wh{\tau}^{b,q_{0}}_{D}>{t/4} \big)\le 1\le
\big( 1\vee \frac{\sqrt{T}}{ (r^*/16)^{\alpha/2}} \big)
\big(1\wedge\frac{\delta_{D}(x)^{\alpha/2}}{\sqrt{t}}\big).\nonumber
\end{equation}
Thus \eqref{lem7.2.2} holds in this case.
It remains to consider the case $\delta_{D}(x)<r/16\le r^{*}/16$. Let $z_{x}\in \partial D$ be such that $|x-z_{x}|=\delta_{D}(x)$. Let $U:=V(z_{x},r/16)$ be a $C^{1,1}$ domain such that $D\cap B(z_{x},r/16)\subset U\subset D\cap B(z_{x},r/8)$. Then we have
\begin{eqnarray}
\P_{x}\big(\wh{\tau}^{b,q_{0}}_{D}>{t/4} \big)&\le&\P_{x}\big(\wh{\tau}^{b,q_{0}}_{U}>{t/4} \big)
+\P_{x}\big(\wh{X}^{b,q_{0}}_{\wh{\tau}^{b,q_{0}}_{U}}\in D\big)\nonumber\\
&\le&\frac{4}{t}\E_{x}[\wh{\tau}^{b,q_{0}}_{U}]
+\P_{x}\big(\wh{X}^{b,q_{0}}_{\wh{\tau}^{b,q_{0}}_{U}}\in D\big).\label{n4.7}
\end{eqnarray}
Take $x_{0}\in \R^{d}$ such that
\[B:=B(x_{0},r/32)\subset B(z_{x}+r_{0}\mathbf{n}(z_{x}),r_{0})\cap \big(B(z_{x},r/4)\setminus B(z_{x},3r/16)\big).\]
Here $\mathbf{n}(z_{x})$ denotes the unit inward normal to $\partial D$ at $z_{x}$.
By \eqref{poissonforwhxbq},
we have
\begin{eqnarray}
\P_{x}\big(\wh{X}^{b,q_{0}}_{\wh{\tau}^{b,q_{0}}_{U}}\in D\big)
&=&\int_{D\setminus U}\int_{U}\wh{G}^{b,q_{0}}_{U}(x,y)j^{b}(z,y)
\frac{h^{b,q_{0}}(z)}{h^{b,q_{0}}(y)}{\mathrm{d}y\mathrm{d}z}\nonumber\\
&=&\E_{x}\big[\int_{0}^{\wh{\tau}^{b,q_{0}}_{U}}\int_{D\setminus U}j^{b}(z,\wh{X}^{b,q_{0}}_{s})\frac{h^{b,q_{0}}(z)}{h^{b,q_{0}}(\wh{X}^{b,q_{0}}_{s})}{\mathrm{d}z\mathrm{d}s}\big]\nonumber\\
&\ge&\E_{x}\big[
\int_{0}^{\wh{\tau}^{b,q_{0}}_{U}}\int_{B}
j^{b}(z,\wh{X}^{b,q_{0}}_{s})\frac{h^{b,q_{0}}(z)}{h^{b,q_{0}}(\wh{X}^{b,q_{0}}_{s})}{\mathrm{d}z\mathrm{d}s}\big)\nonumber\\
&\ge&M^{-1}\E_{x}\big[
\int_{0}^{\wh{\tau}^{b,q_{0}}_{U}}\int_{B}
j^{b}(z,\wh{X}^{b,q_{0}}_{s}){\mathrm{d}z\mathrm{d}s}\big].\label{n4.8}
\end{eqnarray}
For $s\in (0,\wh{\tau}^{b,q_{0}}_{U})$ and $z\in B$, $\wh{X}^{b,q_{0}}_{s}\in U\subset B(z_{x},r/8)$ and
\begin{equation}\nonumber
\frac{r}{16}\le |\wh{X}^{b,q_{0}}_{s}-z|\le \frac{3r}{8}<\eta(A).
\end{equation}
By this, \eqref{jb&j} and \eqref{n4.8}, we have
\begin{eqnarray}
\P_{x}\big(\wh{X}^{b,q_{0}}_{\wh{\tau}^{b,q_{0}}_{U}}\in D\big)
&\ge&\frac{1}{2}M^{-1}\E_{x}\Big[
\int_{0}^{\wh{\tau}^{b,q_{0}}_{U}}\int_{B}
j(z,\wh{X}^{b,q_{0}}_{s}){\mathrm{d}z\mathrm{d}s}\Big]\nonumber\\
&\ge&\frac{1}{2}M^{-1}\mathcal{A}(d,-\alpha)\big(\frac{3r}{8}\big)^{-d-\alpha}
|B|\,\E_{x}[\wh{\tau}^{b,q_{0}}_{U}]\nonumber\\
&\ge&
c_{1}(d,\allowbreak \alpha,\allowbreak M)
r^{-\alpha}\E_{x}[\wh{\tau}^{b,q_{0}}_{U}]\ge
c_{2}(d,\allowbreak \alpha,\allowbreak M,\allowbreak r^{*},\allowbreak T)
t^{-1}\E_{x}[\wh{\tau}^{b,q_{0}}_{U}].\label{n4.9}
\end{eqnarray}
By \eqref{n4.7} and \eqref{n4.9}, we have
\begin{equation}\label{n4.10}
\P_{x}\big(\wh{\tau}^{b,q_{0}}_{D}>t/4\big)\le
c_{3}(d,\allowbreak \alpha,\allowbreak M,\allowbreak r^{*},\allowbreak T)
\P_{x}\big(\wh{X}^{b,q_{0}}_{\wh{\tau}^{b,q_{0}}_{U}}\in D\big).
\end{equation}
By \eqref{poissonforwhxbq}, \eqref{n3.5}, \eqref{n2.7} and \eqref{them4.8.main},
\begin{eqnarray}
\P_{x}\big(\wh{X}^{b,q_{0}}_{\wh{\tau}^{b,q_{0}}_{U}}\in D\big)
&=&\int_{D\setminus U}\int_{U}G^{b,q_{0}}_{U}(y,x)j^{b}(z,y)\frac{h^{b,q_{0}}(z)}{h^{b,q_{0}}(x)}{\mathrm{d}y\mathrm{d}z}\nonumber\\
&\le&M\int_{D\setminus U}\int_{U}G^{b,q_{0}}_{U}(y,x)j^{b}(z,y){\mathrm{d}y\mathrm{d}z}\nonumber\\
&\le&3M\int_{D\setminus U}\int_{U}G^{A}_{U}(y,x)j^{A}(z,y){\mathrm{d}y\mathrm{d}z}\nonumber\\
&\le&3M\P_{x}\big(X^{A}_{\tau^{A}_{U}}\in D\big).\label{n4.11}
\end{eqnarray}
Here $X^{A}$ denotes the Feller process generated by the operator $\Delta^{\alpha/2}+A\Delta^{\beta/2}$, and $G^{A}_{U}(x,y)$ its Green function of $U$.
The second inequality in \eqref{n4.11} uses $j^{b}\le j^{A}$, where $j^{A}$ is the jump kernel of $X^{A}$, and $G^{b,q_{0}}_{U}\le \frac{3}{2}G_{U}\le 3G^{A}_{U}$.
The last two bounds follow from \eqref{n2.7} and from \eqref{them4.8.main} applied to $b\equiv A$. Now with $x_{1}=z_{x}+r\mathbf{n}(z_{x})/16$, by applying boundary Harnack inequality for $\Delta^{\alpha/2}+A\Delta^{\beta/2}$ (see, for example, \cite[Theorem 1.3]{CKS3}), we get
\begin{equation}
\P_{x}\big(X^{A}_{\tau^{A}_{U}}\in D\big)\le
c_{4}\P_{x_{1}}\big(X^{A}_{\tau^{A}_{U}}\in D\big)\frac{\delta_{D}(x)^{\alpha/2}}{\delta_{D}(x_{1})^{\alpha/2}}\le
c_{5}\frac{\delta_{D}(x)^{\alpha/2}}{\sqrt{t}}
,
\label{n4.12}
\end{equation}
where $c_{4}=c_{4}(d,\allowbreak \alpha,\allowbreak \beta,\allowbreak R_{0},\allowbreak \Lambda_{0},\allowbreak A)$
and $c_{5}=c_{5}(d,\allowbreak \alpha,\allowbreak \beta,\allowbreak R_{0},\allowbreak \Lambda_{0},\allowbreak A,\allowbreak r^{*},\allowbreak T)$ are positive constants.
Consequently, \eqref{lem7.2.2} follows from \eqref{n4.10}, \eqref{n4.11} and \eqref{n4.12}.
\end{proof}

For every $a\ge 0$, define
\[
J^{a}(x,y):=|x-y|^{-d-\alpha}+a |x-y|^{-d-\beta}.
\]
We write $J(x,y)$ for $J^{0}(x,y)$.
Recall that $p_{a}(t,x,y)$ denotes the transition density corresponding to the generator $\Delta^{\alpha/2}+a\Delta^{\beta/2}$.
For $\lambda>0$ and $b\in\mathcal{B}_{A}$, we write $a:=M_{b^{+},\lambda}$ in the rest of this subsection. It follows by \eqref{pa} that for every $\lambda,T>0$, $t\in (0,T]$ and
 $b\in\mathcal{B}_{A}$,
\begin{equation}\label{n4.13'}
p_{a}(t,x,y)\stackrel{c_{1}}{\asymp} t^{-d/\alpha}\wedge t J^{a}(x,y),
\end{equation}
where $c_{1}=c_{1}(d,\allowbreak \alpha,\allowbreak \beta,\allowbreak A,\allowbreak T)$ is a positive constant and $M_{b^{+},\lambda}=\esup_{x,z\in \R^{d},|z|>\lambda}b^{+}(x,z)$. Moreover, for a.e. $x,y\in \R^{d}$,
\begin{eqnarray}
j^{b}(x,y)&=&\frac{\mathcal{A}(d,-\alpha)}{|x-y|^{d+\alpha}}+\mathcal{A}(d,-\beta)\frac{b(x,y-x)}{|x-y|^{d+\beta}}\nonumber\\
&\le&\frac{\mathcal{A}(d,-\alpha)}{|x-y|^{d+\alpha}}+\mathcal{A}(d,-\beta)\frac{b^{+}(x,y-x)}{|x-y|^{d+\beta}}\nonumber\\
&=&\frac{\mathcal{A}(d,-\alpha)}{|x-y|^{d+\alpha}}+
\mathcal{A}(d,-\beta)|x-y|^{\alpha-\beta}\frac{b^{+}(x,y-x)}{|x-y|^{d+\alpha}}\1_{\{|x-y|\le \lambda\}}\nonumber\\
&&+\mathcal{A}(d,-\beta)\frac{b^{+}(x,y-x)}{|x-y|^{d+\beta}}\1_{\{|x-y|>\lambda\}}\nonumber\\
&\le&\big(\mathcal{A}(d,-\alpha)+A\lambda^{\alpha-\beta}\mathcal{A}(d,-\beta)\big)\frac{1}{|x-y|^{d+\alpha}}
+\mathcal{A}(d,-\beta)\frac{a}{|x-y|^{d+\beta}}\nonumber\\
&\le&c_{2}J^{a}(x,y).\label{n4.14'}
\end{eqnarray}
Here $c_{2}=c_{2}(d,\allowbreak \alpha,\allowbreak \beta,\allowbreak A,\allowbreak \lambda)$ is a positive constant.

\begin{lemma}\label{nlem4.4}
Suppose $D$ is a $C^{1,1}$ open set
with characteristics $(R_{0},\Lambda_{0})$.
For every $\lambda, T\in (0,{\infty})$, there is a positive constant
$ C_{40}= C_{40}(d,\allowbreak \alpha,\allowbreak \beta,\allowbreak R_{0},\allowbreak \Lambda_{0},\allowbreak A,\allowbreak M,\allowbreak T,\allowbreak \lambda)$
such that for every $t\in (0,T]$ and $x,y\in D$,
\begin{equation}
p^{b,q_{0}}_{D}(t/2,x,y)\le  C_{40}\big(1\wedge \frac{\delta_{D}(x)^{\alpha/2}}{\sqrt{t}}\big)\big(t^{-d/\alpha}\wedge tJ^{a}(x,y)\big),\label{n4.13}
\end{equation}
\begin{equation}
p^{b,q_{0}}_{D}(t/2,x,y)\le  C_{40}\big(1\wedge \frac{\delta_{D}(y)^{\alpha/2}}{\sqrt{t}}\big)\big(t^{-d/\alpha}\wedge tJ^{a}(x,y)\big).\label{n4.14}
\end{equation}
If, in addition,
\eqref{condi5} is satisfied for some $\eps\in (0,1]$, then there is a positive constant $ C_{41}=
C_{41}(d,\allowbreak \alpha,\allowbreak \beta,\allowbreak R_{0},\allowbreak \Lambda_{0},\allowbreak A,\allowbreak M,\allowbreak T,\allowbreak \eps)
$
such that for every $t\in (0,T]$ and $x,y\in D$,
\begin{equation}
p^{b,q_{0}}_{D}(t/2,x,y)\le  C_{41}\big(1\wedge \frac{\delta_{D}(x)^{\alpha/2}}{\sqrt{t}}\big)\big(t^{-d/\alpha}\wedge tJ(x,y)\big),\label{7.19}
\end{equation}
\begin{equation}
p^{b,q_{0}}_{D}(t/2,x,y)\le  C_{41}\big(1\wedge \frac{\delta_{D}(y)^{\alpha/2}}{\sqrt{t}}\big)\big(t^{-d/\alpha}\wedge tJ(x,y)\big).\label{7.20}
\end{equation}
\end{lemma}

\begin{proof}
We prove only \eqref{n4.14} and \eqref{7.20}. The proofs of \eqref{n4.13} and \eqref{7.19} are similar.
Fix $x,y\in D$, $T>0$ and $t\in (0,T]$. Let $r^{*}$ be the constant given in \eqref{e:rstar} and $r:=\big(\frac{t}{T}\big)^{1/\alpha}r^{*}$. If $|x-y|/8\le r$, then
\begin{equation}
t^{-d/\alpha}\wedge tJ^{a}(x,y)\ge t^{-d/\alpha}\wedge tJ(x,y)\ge t^{-d/\alpha}
\big(1\wedge \frac{T^{1/\alpha}}{8r^{*}}\big)^{d+\alpha} .  \label{n4.15}
\end{equation}
By the Chapman--Kolmogorov equation, \eqref{n3.4},
the upper bound of \eqref{ltdensity},
\eqref{n4.13'} and \eqref{lem7.2.2},
\begin{eqnarray}
p^{b,q_{0}}_{D}(t/2,x,y)&\le&\int_{D}p^{b,q_{0}}(t/4,x,z)p^{b,q_{0}}_{D}(t/4,z,y){\mathrm{d}z}\nonumber\\
&=&\int_{D}p^{b,q_{0}}(t/4,x,z)\wh{p}^{b,q_{0}}_{D}(t/4,y,z)\frac{h^{b,q_{0}}(y)}{h^{b,q_{0}}(z)}{\mathrm{d}z}\nonumber\\
&\le&c_{1}Me^{-q_{0}t/4}\int_{D}p_{a}(t/4,x,z)\wh{p}^{b,q_{0}}_{D}(t/4,y,z){\mathrm{d}z}\nonumber\\
&\le&c_{2}t^{-d/\alpha}\P_{y}\big(\wh{\tau}^{b,q_{0}}_{D}>t/4\big)\nonumber\\
&\le&c_{3}t^{-d/\alpha}\big(1\wedge \frac{\delta_{D}(y)^{\alpha/2}}{\sqrt{t}}\big).\label{n4.16}
\end{eqnarray}
Here $c_{1}=c_{1}(d,\allowbreak \alpha,\allowbreak \beta,\allowbreak A,\allowbreak \lambda)$, $c_{2}=c_{2}(d,\allowbreak \alpha,\allowbreak \beta,\allowbreak A,\allowbreak M,\allowbreak T,\allowbreak \lambda)$
and $c_{3}=c_{3}(d,\allowbreak \alpha,\allowbreak \beta,\allowbreak R_{0},\allowbreak \Lambda_{0},\allowbreak A,\allowbreak M,\allowbreak T,\allowbreak \lambda)$ are positive constants.
In this case, \eqref{n4.14} follows from \eqref{n4.15} and \eqref{n4.16}. We can prove \eqref{7.20} similarly by using \eqref{2.6} and \eqref{p0} instead of \eqref{ltdensity} and \eqref{n4.13'} respectively in the above argument.
If $\delta_{D}(y)\ge r/16$, then
\begin{equation}
1\wedge \frac{\delta_{D}(y)^{\alpha/2}}{\sqrt{t}}\ge 1\wedge
\frac{(r^{*}/16 )^{\alpha/2}}{\sqrt{T}} .\label{7.23}
\end{equation}
By
the upper bound of \eqref{ltdensity}
and \eqref{n4.13'}, for $t\in (0,T]$ and $x,y\in D$,
\begin{equation}\label{7.24}
p^{b,q_{0}}_{D}(t/2,x,y)\le p^{b}(t/2,x,y)\le c_{4}\big(t^{-d/\alpha}\wedge t J^{a}(x,y)\big)
\end{equation}
for some $c_{4}=c_{4}(d,\allowbreak \alpha,\allowbreak \beta,\allowbreak A,\allowbreak T,\allowbreak \lambda)>0$. Moreover, if \eqref{condi5} is satisfied for some $\eps\in (0,1]$, then by \eqref{2.6} and \eqref{p0}, for $t\in (0,T]$ and $x,y\in D$,
\begin{equation}\label{7.25}
p^{b,q_{0}}_{D}(t/2,x,y)\le p^{b}(t/2,x,y)\le c_{5}\big(t^{-d/\alpha}\wedge t J(x,y)\big)
\end{equation}
for some $c_{5}=c_{5}(d,\allowbreak \alpha,\allowbreak \beta,\allowbreak A,\allowbreak T,\allowbreak \eps)>0$. Thus \eqref{n4.14} follows from \eqref{7.23} and \eqref{7.24}, and \eqref{7.20} follows from \eqref{7.23} and \eqref{7.25}.
It remains to consider the case $|x-y|/8>r$ and $\delta_{D}(y)<r/16$. Let $z_{y}\in\partial D$ be such that $|y-z_{y}|=\delta_{D}(y)$. Let $U_{1}:=V(z_{y},r/8)$, $U_{3}:=\{z\in D:\ |z-y|>|x-y|/2\}$ and $U_{2}:=D\setminus (U_{1}\cup U_{3})$. Then by \eqref{lem7.3.2}, we have
\begin{equation}
p^{b,q_{0}}_{D}(t/2,x,y)\le M\P_{y}\big(\wh{X}^{b,q_{0}}_{\wh{\tau}^{b,q_{0}}_{U_{1}}}\in U_{2}\big)
\sup_{\substack{s<t/2\\ z\in U_{2}}}p^{b,q_{0}}_{D}(s,x,z)+M\big(\frac{t}{2}\wedge \E_{y}[\wh{\tau}^{b,q_{0}}_{U_{1}}]\big)\esup_{\substack{u\in U_{1}\\ z\in U_{3}}}j^{b}(z,u).\label{n4.19}
\end{equation}
For $s\in (0,t/2)$ and $z\in U_{2}$, we have $|x-z|\ge |x-y|-|y-z|\ge |x-y|/2$. Thus by \eqref{7.24}
\begin{equation}
p^{b,q_{0}}_{D}(s,x,z)\le c_{4} \big(s^{-d/\alpha}\wedge s J^{a}(x,z)\big)\le c_{6}t J^{a}(x,y),\label{n4.20}
\end{equation}
for some $c_{6}=c_{6}(d,\allowbreak \alpha,\allowbreak \beta,\allowbreak A,\allowbreak T,\allowbreak \lambda)>0$.
Moreover, if \eqref{condi5} is satisfied, then by \eqref{7.25},
\begin{equation}
p^{b,q_{0}}_{D}(s,x,z)\le c_{5} \big(s^{-d/\alpha}\wedge s J(x,z)\big)\le c_{7}t J(x,y),\label{n4.20'}
\end{equation}
for some $c_{7}=c_{7}(d,\allowbreak \alpha,\allowbreak \beta,\allowbreak A,\allowbreak T,\allowbreak \eps)>0$.
For $u\in U_{1}$ and $z\in U_{3}$, we have $|u-z|\ge |z-y|-|z_{y}-y|-|z_{y}-u|\ge |x-y|/2-
r/16-r/4\ge |x-y|/4$. Thus by \eqref{n4.14'},
\begin{equation}
\esup_{\substack{u\in U_{1}\\ z\in U_{3}}}j^{b}(z,u)\le
c_{8}(d,\allowbreak \alpha,\allowbreak \beta,\allowbreak A,\allowbreak \lambda)
\esup_{\substack{u\in U_{1}\\ z\in U_{3}}}J^{a}(z,u)\le
c_{9}(d,\allowbreak \alpha,\allowbreak \beta,\allowbreak A,\allowbreak \lambda)
J^{a}(x,y).\label{n4.21}
\end{equation}
If \eqref{condi5} is satisfied, then
\begin{equation}
\esup_{\substack{u\in U_{1}\\ z\in U_{3}}}j^{b}(z,u)\le \eps^{-1}\esup_{\substack{u\in U_{1}\\ z\in U_{3}}}J(z,u)\le
c_{10}(d,\allowbreak \alpha,\allowbreak \eps)
J(x,y).\label{n4.21'}
\end{equation}
It follows by \eqref{n4.19}, \eqref{n4.20} and \eqref{n4.21} that
\begin{eqnarray}
p^{b,q_{0}}_{D}(t/2,x,y)&\le&c_{11}J^{a}(x,y)\big(t\P_{y}\big(
\wh{X}^{b,q_{0}}_{\wh{\tau}^{b,q_{0}}_{U_{1}}}\in U_{2}\big)+\E_{y} \big[ \wh{\tau}^{b,q_{0}}_{U_{1}} \big] \big)\nonumber\\
&\le&c_{11}J^{a}(x,y)\big(t\P_{y}\big(
\wh{X}^{b,q_{0}}_{\wh{\tau}^{b,q_{0}}_{U_{1}}}\in D\big)+\E_{y} \big[ \wh{\tau}^{b,q_{0}}_{U_{1}} \big] \big)\label{n4.22}
\end{eqnarray}
for some $c_{11}=c_{11}(d,\allowbreak \alpha,\allowbreak \beta,\allowbreak A,\allowbreak M,\allowbreak T,\allowbreak \lambda)>0$.
Similarly, if \eqref{condi5} is satisfied, it follows by \eqref{n4.19}, \eqref{n4.20'} and \eqref{n4.21'} that
\begin{equation*}
p^{b,q_{0}}_{D}(t/2,x,y)\le c_{12}J(x,y)\big(t\P_{y}\big(
\wh{X}^{b,q_{0}}_{\wh{\tau}^{b,q_{0}}_{U_{1}}}\in D\big)+\E_{y} \big[ \wh{\tau}^{b,q_{0}}_{U_{1}} \big] \big)
\end{equation*}
for some $c_{12}=c_{12}(d,\allowbreak \alpha,\allowbreak \beta,\allowbreak A,\allowbreak M,\allowbreak T,\allowbreak \eps)>0$.
It follows from the argument in the proof of Lemma \ref{lem7.2} that
\[
\E_{y} \big[ \wh{\tau}^{b,q_{0}}_{U_{1}} \big] \lesssim t\P_{y}\big(
\wh{X}^{b,q_{0}}_{\wh{\tau}^{b,q_{0}}_{U_{1}}}\in D\big)\le c_{13}\sqrt{t}\delta_{D}(y)^{\alpha/2}
\]
 for some $c_{13}=
c_{13}(d,\allowbreak \alpha,\allowbreak \beta,\allowbreak R_{0},\allowbreak \Lambda_{0},\allowbreak A,\allowbreak M,\allowbreak T)>0
 $. By this and  \eqref{n4.22}, we have
\begin{equation} \label{n4.23}
p^{b,q_{0}}_{D}(t/2,x,y)\le
c_{14}(d,\allowbreak \alpha,\allowbreak \beta,\allowbreak R_{0},\allowbreak \Lambda_{0},\allowbreak A,\allowbreak M,\allowbreak T,\allowbreak \lambda)
\big(\frac{\delta_{D}(y)^{\alpha/2}}{\sqrt{t}}\big)\, t J^{a}(x,y) .
\end{equation}
Similarly if \eqref{condi5} is satisfied, then
\begin{equation}\label{n4.23'}
p^{b,q_{0}}_{D}(t/2,x,y)\le
c_{15}(d,\allowbreak \alpha,\allowbreak \beta,\allowbreak R_{0},\allowbreak \Lambda_{0},\allowbreak A,\allowbreak M,\allowbreak T,\allowbreak \eps)
\big(\frac{\delta_{D}(y)^{\alpha/2}}{\sqrt{t}}\big) \, t J(x,y) .
\end{equation}
Recall that $r=\big(\frac{t}{T}\big)^{1/\alpha}r^{*}$. For $x,y\in D$
with $\delta_{D}(y)<r/16$ and $|x-y|/8>r$, we have
\begin{equation}
1\wedge \frac{\delta_{D}(y)^{\alpha/2}}{\sqrt{t}}
\stackrel{c_{16}(r^{*},T)}{\asymp}
\frac{\delta_{D}(y)^{\alpha/2}}{\sqrt{t}},\label{n4.24}
\end{equation}
 and
 \begin{eqnarray}
tJ^{a}(x,y)&\le& t(1\vee A)\big(|x-y|^{-d-\alpha}+|x-y|^{-d-\beta}\big)\nonumber\\
&\stackrel{c_{17}(A,r^{*},T)}{\lesssim}&
t^{-d/\alpha}+t^{-d/\alpha+(1-\beta/\alpha)}\le t^{-d/\alpha}\big(1+T^{1-\beta/\alpha}\big).\nonumber
\end{eqnarray}
Thus
\begin{equation}
t^{-d/\alpha}\wedge tJ^{a}(x,y)
\stackrel{c_{18}(A,r^{*},T)}{\asymp}
tJ^{a}(x,y).\label{n4.25}
\end{equation}
Similarly, one proves that for $x,y\in D$ with $|x-y|/8>r$,
\begin{equation}
t^{-d/\alpha}\wedge tJ(x,y)
\stackrel{c_{19}(r^{*},T)}{\asymp}
tJ(x,y).\label{n4.25'}
\end{equation}
Therefore \eqref{n4.14} follows from \eqref{n4.23}, \eqref{n4.24} and \eqref{n4.25}, and \eqref{7.20} from \eqref{n4.23'}, \eqref{n4.24} and \eqref{n4.25'}.\end{proof}

\begin{lemma}\label{nlem4.5}
Suppose $D$ is a $C^{1,1}$ open set
with characteristics $(R_{0},\Lambda_{0})$.
For every $\lambda, T\in (0,{\infty})$, there is a positive constant
$ C_{42}= C_{42}(d,\allowbreak \alpha,\allowbreak \beta,\allowbreak R_{0},\allowbreak \Lambda_{0},\allowbreak A,\allowbreak M,\allowbreak T,\allowbreak \lambda)$
such that for every $t\in (0,T]$ and $x,y\in D$,
\begin{equation}\label{upforpbq}
p^{b,q_{0}}_{D}(t,x,y)\le  C_{42}\big(1\wedge \frac{\delta_{D}(x)^{\alpha/2}}{\sqrt{t}}\big)\big(1\wedge \frac{\delta_{D}(y)^{\alpha/2}}{\sqrt{t}}\big)\big(t^{-d/\alpha}\wedge tJ^{a}(x,y)\big).
\end{equation}
If, in addition,
\eqref{condi5} is satisfied for some $\eps\in (0,1]$, then there is a positive constant
$ C_{43}= C_{43}(d,\allowbreak \alpha,\allowbreak \beta,\allowbreak R_{0},\allowbreak \Lambda_{0},\allowbreak A,\allowbreak M,\allowbreak T,\allowbreak \eps)$
such that for every $t\in (0,T]$ and $x,y\in D$,
\begin{equation}\label{upforpbq'}
p^{b,q_{0}}_{D}(t,x,y)\le  C_{43}\big(1\wedge \frac{\delta_{D}(x)^{\alpha/2}}{\sqrt{t}}\big)\big(1\wedge \frac{\delta_{D}(y)^{\alpha/2}}{\sqrt{t}}\big)\big(t^{-d/\alpha}\wedge tJ(x,y)\big).
\end{equation}
\end{lemma}

\begin{proof}
We prove only \eqref{upforpbq}. The proof of \eqref{upforpbq'} is similar. By the Chapman--Kolmogorov equation, \eqref{n4.13}, \eqref{n4.14} and \eqref{n4.13'},
for every $t\in (0,T]$ and $x,y\in D$,
\begin{eqnarray}
&&p^{b,q_{0}}_{D}(t,x,y)\nonumber\\
&=&\int_{D}p^{b,q_{0}}_{D}(t/2,x,z)p^{b,q_{0}}_{D}(t/2,z,y){\mathrm{d}z}\nonumber\\
&\le&c_{1}\big(1\wedge \frac{\delta_{D}(x)^{\alpha/2}}{\sqrt{t}}\big)\big(1\wedge \frac{\delta_{D}(y)^{\alpha/2}}{\sqrt{t}}\big)\int_{D}\big(t^{-d/\alpha}\wedge tJ^{a}(x,z)\big)
\big(t^{-d/\alpha}\wedge tJ^{a}(z,y)\big){\mathrm{d}z}\nonumber\\
&\le&c_{2}\big(1\wedge \frac{\delta_{D}(x)^{\alpha/2}}{\sqrt{t}}\big)\big(1\wedge \frac{\delta_{D}(y)^{\alpha/2}}{\sqrt{t}}\big)\int_{D}p_{a}(t/2,x,z)p_{a}(t/2,z,y){\mathrm{d}z}\nonumber\\
&\le&c_{2}\big(1\wedge \frac{\delta_{D}(x)^{\alpha/2}}{\sqrt{t}}\big)\big(1\wedge \frac{\delta_{D}(y)^{\alpha/2}}{\sqrt{t}}\big)p_{a}(t,x,y)\nonumber\\
&\le&c_{3}\big(1\wedge \frac{\delta_{D}(x)^{\alpha/2}}{\sqrt{t}}\big)\big(1\wedge \frac{\delta_{D}(y)^{\alpha/2}}{\sqrt{t}}\big)\big(t^{-d/\alpha}\wedge tJ^{a}(x,y)\big).\nonumber
\end{eqnarray}
This completes the proof.\end{proof}

\begin{thm}\label{stupbound}
Suppose $D$ is a $C^{1,1}$ open set
with characteristics $(R_{0},\Lambda_{0})$.
For every $\lambda, T\in (0,{\infty})$, there is a positive constant
$ C_{44}= C_{44}(d,\allowbreak \alpha,\allowbreak \beta,\allowbreak R_{0},\allowbreak \Lambda_{0},\allowbreak A,\allowbreak M,\allowbreak T,\allowbreak \lambda)$
such that for every $t\in (0,T]$ and $x,y\in D$,
\begin{equation}\label{7.6.1}
p^{b}_{D}(t,x,y)\le  C_{44}\big(1\wedge \frac{\delta_{D}(x)^{\alpha/2}}{\sqrt{t}}\big)\big(1\wedge \frac{\delta_{D}(y)^{\alpha/2}}{\sqrt{t}}\big)\big(t^{-d/\alpha}\wedge tJ^{a}(x,y)\big).
\end{equation}
If, in addition,
\eqref{condi5} is satisfied for some $\eps\in (0,1]$, then there is a positive constant
$ C_{45}= C_{45}(d,\allowbreak \alpha,\allowbreak \beta,\allowbreak R_{0},\allowbreak \Lambda_{0},\allowbreak A,\allowbreak M,\allowbreak T,\allowbreak \eps)$
such that for every $t\in (0,T]$ and $x,y\in D$,
\begin{equation*}
p^{b}_{D}(t,x,y)\le  C_{45}\big(1\wedge \frac{\delta_{D}(x)^{\alpha/2}}{\sqrt{t}}\big)\big(1\wedge \frac{\delta_{D}(y)^{\alpha/2}}{\sqrt{t}}\big)\big(t^{-d/\alpha}\wedge \frac{t}{|x-y|^{d+\alpha}}\big).
\end{equation*}
\end{thm}
\begin{proof}
This theorem follows directly from \eqref{n1.2} and
Lemma \ref{nlem4.5}.
\end{proof}

\subsection{Small-time lower bound}\label{st.2}
Let $\eta(A)$, $M$ and $T_{0}$ be the constants given in \eqref{defiforetaA}, \eqref{defiforM} and \eqref{e:rstar}, respectively.
\begin{proposition}\label{nprop4.6}
Let $a_{1}\in (0,1)$, $a_{3}>a_{2}>0$ and $R\in (0,1/2]$.
There is a positive constant $ C_{46}= C_{46}(d,\allowbreak \alpha,\allowbreak \beta,\allowbreak a_{1},\allowbreak a_{2},\allowbreak a_{3},\allowbreak R,\allowbreak A)$ such that for any $x_{0}\in \R^{d}$, $r\in (0,R]$, $t\in [a_{2}r^{\alpha},a_{3}r^{\alpha}]$ and $x,y\in B(x_{0},a_{1}r)$
\begin{equation}\label{nprop4.6.1}
p^{b,q_{0}}_{B(x_{0},r)}(t,x,y)\ge  C_{46} r^{-d},
\end{equation}
\begin{equation}\label{nprop4.6.2}
\wh{p}^{b,q_{0}}_{B(x_{0},r)}(t,x,y)\ge  C_{46} M^{-1} r^{-d}.
\end{equation}
If, in addition, $b$ satisfies \eqref{condi3} for some $\eps>0$, then the same holds for every $R>0$, with $C_{46}$ depending also on $\eps$.
\end{proposition}

\begin{proof} This proposition follows directly from Proposition \ref{prop3.5}, \eqref{n1.2} and \eqref{n3.4}.\end{proof}

The following two-point estimate is used in the proof of Lemma \ref{nlem4.8}.

\begin{lemma}\label{nlem4.9}
There is a positive constant
$ C_{47}= C_{47}(d,\allowbreak \alpha,\allowbreak \beta,\allowbreak A,\allowbreak M)$
such that for every $t\in (0,T_{0}]$ and every $u,v\in \R^{d}$
with $|u-v|\le 153\eta(A)/160$,
\begin{equation}\label{nlem4.9.1}
p^{b,q_{0}}_{B(u,t^{1/\alpha}/2)\cup B(v,t^{1/\alpha}/2)}(t/3,u,v)\ge  C_{47}\big(t^{-d/\alpha}\wedge t J(u,v)\big).
\end{equation}
Moreover,
if \eqref{condi3} is satisfied for some $\eps>0$, then there is a positive constant
$ C_{48}= C_{48}(d,\allowbreak \alpha,\allowbreak \beta,\allowbreak A,\allowbreak M,\allowbreak \eps)$
such that for every $t\in (0,T_{0}]$ and every $u,v\in \R^{d}$,
\begin{equation}\label{nlem4.9.2}
p^{b,q_{0}}_{B(u,t^{1/\alpha}/2)\cup B(v,t^{1/\alpha}/2)}(t/3,u,v)\ge  C_{48}\big(t^{-d/\alpha}\wedge t J(u,v)\big).
\end{equation}
Furthermore,
if \eqref{condi4} is satisfied for some $\eps>0$, then there is a positive constant
$ C_{49}= C_{49}(d,\allowbreak \alpha,\allowbreak \beta,\allowbreak A,\allowbreak M,\allowbreak \eps)$
such that for every $t\in (0,T_{0}]$ and every $u,v\in \R^{d}$,
\begin{equation}\label{nlem4.9.3}
p^{b,q_{0}}_{B(u,t^{1/\alpha}/2)\cup B(v,t^{1/\alpha}/2)}(t/3,u,v)\ge  C_{49}\big(t^{-d/\alpha}\wedge t J^{1}(u,v)\big).
\end{equation}
\end{lemma}

\begin{proof} Fix $t\in (0,T_{0}]$ and $u,v\in \R^{d}$. If $|u-v|\le t^{1/\alpha}/4$, then by \eqref{nprop4.6.1},
\begin{equation}\nonumber
p^{b,q_{0}}_{B(u,t^{1/\alpha}/2)\cup B(v,t^{1/\alpha}/2)}(t/3,u,v)\ge \inf_{\substack{w\in\R^{d}\\ |u-w|\le t^{1/\alpha}/4}}
p^{b,q_{0}}_{B(u,t^{1/\alpha}/2)}(t/3,u,w)\ge
c_{0}(d,\allowbreak \alpha,\allowbreak \beta,\allowbreak A)
t^{-d/\alpha}.
\end{equation}
Thus \eqref{nlem4.9.1}, \eqref{nlem4.9.2} and \eqref{nlem4.9.3} hold in this case. It remains to consider the case $|u-v|>t^{1/\alpha}/4$.
Suppose $t^{1/\alpha}/4<|u-v|\le 153\eta(A)/160$. Let $U_{1}:=B(u,t^{1/\alpha}/80)$ and $U_{3}:=B(v,t^{1/\alpha}/80)$, so that $\mathrm{dist}(U_{1},U_{3})>0$. By \eqref{lem7.3.3}, \eqref{nprop4.6.1} and \eqref{nprop4.6.2}, we have
\begin{eqnarray}
&&p^{b,q_{0}}_{B(u,t^{1/\alpha}/2)\cup B(v,t^{1/\alpha}/2)}(t/3,u,v)\nonumber\\
&\ge&\frac{t}{3M}\P_{u}\big(\tau^{b,q_{0}}_{U_{1}}>\frac{t}{3}\big)
\P_{v}\big(\wh{\tau}^{b,q_{0}}_{U_{3}}>\frac{t}{3}\big)\einf_{w\in U_{1},z\in U_{3}}j^{b}(w,z)\nonumber\\
&\ge&\frac{t}{3M}\einf_{w\in U_{1},z\in U_{3}}j^{b}(w,z)\Big(\int_{B(u,t^{1/\alpha}/160)}p^{b,q_{0}}_{U_{1}}(t/3,u,x){\mathrm{d}x}\Big)
\Big(\int_{B(v,t^{1/\alpha}/160)}\wh{p}^{b,q_{0}}_{U_{3}}(t/3,v,y){\mathrm{d}y}\Big)\nonumber\\
&\ge&
c_{1}(d,\allowbreak \alpha,\allowbreak \beta,\allowbreak A,\allowbreak M)
t\einf_{w\in U_{1},z\in U_{3}}j^{b}(w,z)\big(t^{-d/\alpha}|B(u,t^{1/\alpha}/160)|\big)\big(t^{-d/\alpha}M^{-1}|B(v,t^{1/\alpha}/160)|\big)\nonumber\\
&\ge&
c_{2}(d,\allowbreak \alpha,\allowbreak \beta,\allowbreak A,\allowbreak M)
t\einf_{w\in U_{1},z\in U_{3}}j^{b}(w,z).\label{n4.50}
\end{eqnarray}
For $t^{1/\alpha}/4<|u-v|\le 153\eta(A)/160$, $w\in U_{1}$ and $z\in U_{3}$, we have $|w-z|\le |u-v|+t^{1/\alpha}/4 0\le 11|u-v|/10$ and, since $t^{1/\alpha}\le T_{0}^{1/\alpha}\le \eta(A)/64$, $|w-z|<\eta(A)$. Thus by \eqref{jb&j},
\begin{equation}
\einf_{w\in U_{1},z\in U_{3}}j^{b}(w,z)\ge \frac{1}{2}\einf_{w\in U_{1},z\in U_{3}}j(w,z)
\gtrsim
\einf_{w\in U_{1},z\in U_{3}}J(w,z)
\gtrsim
J(u,v). \label{n4.51}
\end{equation}
Consequently, \eqref{nlem4.9.1} follows from \eqref{n4.50} and \eqref{n4.51}.  For every $u,v\in \R^{d}$ with $|u-v|>t^{1/\alpha}/4$, $w\in U_{1}$ and $z\in U_{3}$, $|w-z|\le |u-v|+t^{1/\alpha}/40\le 11|u-v|/10$. If \eqref{condi3} is satisfied for some $\eps>0$, then
\begin{equation}
\einf_{w\in U_{1},z\in U_{3}}j^{b}(w,z)\ge \eps  \einf_{w\in U_{1},z\in U_{3}}J(w,z)
\gtrsim \eps
J(u,v). \label{n4.52}
\end{equation}
Consequently, \eqref{nlem4.9.2} follows from \eqref{n4.50} and \eqref{n4.52}. The proof of \eqref{nlem4.9.3} is similar.
\end{proof}

\begin{lemma}\label{nlem4.8}
Suppose $D$ is a $C^{1,1}$ open set
with characteristics $(R_{0},\Lambda_{0})$.
There is a positive constant
$ C_{50}= C_{50}(d,\allowbreak \alpha,\allowbreak \beta,\allowbreak R_{0},\allowbreak \Lambda_{0},\allowbreak A,\allowbreak M)$
 such that for every $t\in (0,T_{0}]$ and every $x,y\in D$ with $|x-y|<4\eta(A)/5$,
\begin{equation}\label{n4.30}
p^{b,q_{0}}_{D}(t,x,y)\ge  C_{50}\big(1\wedge \frac{\delta_{D}(x)^{\alpha/2}}{\sqrt{t}}\big)
\big(1\wedge \frac{\delta_{D}(y)^{\alpha/2}}{\sqrt{t}}\big)\big(t^{-d/\alpha}\wedge \frac{t}{|x-y|^{d+\alpha}}\big).
\end{equation}
Moreover,
if \eqref{condi3} is satisfied for some $\eps>0$, then there is a positive constant
$ C_{51}= C_{51}(d,\allowbreak \alpha,\allowbreak \beta,\allowbreak R_{0},\allowbreak \Lambda_{0},\allowbreak A,\allowbreak M,\allowbreak \eps)$
such that for every $t\in (0,T_{0}]$ and every $x,y\in D$,
\begin{equation}\label{n4.31}
p^{b,q_{0}}_{D}(t,x,y)\ge  C_{51}\big(1\wedge \frac{\delta_{D}(x)^{\alpha/2}}{\sqrt{t}}\big)
\big(1\wedge \frac{\delta_{D}(y)^{\alpha/2}}{\sqrt{t}}\big)\big(t^{-d/\alpha}\wedge \frac{t}{|x-y|^{d+\alpha}}\big).
\end{equation}
Furthermore, if
\eqref{condi4} is satisfied for some $\eps>0$, then there is a positive constant
$ C_{52}= C_{52}(d,\allowbreak \alpha,\allowbreak \beta,\allowbreak R_{0},\allowbreak \Lambda_{0},\allowbreak A,\allowbreak M,\allowbreak \eps)$
such that for every $t\in (0,T_{0}]$ and every $x,y\in D$,
\begin{equation}\label{n4.32}
p^{b,q_{0}}_{D}(t,x,y)\ge  C_{52}\big(1\wedge \frac{\delta_{D}(x)^{\alpha/2}}{\sqrt{t}}\big)
\big(1\wedge \frac{\delta_{D}(y)^{\alpha/2}}{\sqrt{t}}\big)\big(t^{-d/\alpha}\wedge \big(\frac{t}{|x-y|^{d+\alpha}}+\frac{t}{|x-y|^{d+\beta}}\big)\big).
\end{equation}
\end{lemma}

\begin{proof}
Fix $x,y\in D$ and $t\in (0,T_{0}]$, and note that $3t^{1/\alpha}<r_{0}$ by \eqref{e:rstar}. We choose $z^{t}_{x}$ as follows: If $\delta_{D}(x)\le 3 t^{1/\alpha}$, let $z_{x}\in\partial D$ be such that $|x-z_{x}|=\delta_{D}(x)$ and $z^{t}_{x}:=z_{x}+\frac{9}{2}t^{1/\alpha}\mathbf{n}(z_{x})$, where $\mathbf{n}(z_{x})$ is the unit inward normal to $\partial D$ at $z_{x}$, so that $B(z^{t}_{x},3t^{1/\alpha}/2)\subset
B(z_{x}+3t^{1/\alpha}\mathbf{n}(z_{x}),3t^{1/\alpha})\setminus \{x\}$ and $\delta_{D}(z)\ge t^{1/\alpha}$ for all $z\in B(z^{t}_{x},t^{1/\alpha}/2)$. Else if $\delta_{D}(x)>3 t^{1/\alpha}$, let
$z^{t}_{x}\in B(x,\delta_{D}(x))$ be such that $|x-z^{t}_{x}|=3t^{1/\alpha}/2$, so that $B(z^{t}_{x},3t^{1/\alpha}/2)\subset B(x,\delta_{D}(x))\setminus \{x\}$ and $\delta_{D}(z)\ge t^{1/\alpha}$ for all $z\in B(z^{t}_{x},t^{1/\alpha}/2)$. We also define $z^{t}_{y}$ and $z_{y}$ in the same way.

By the Chapman--Kolmogorov equation, we have
\begin{align}
&p^{b,q_{0}}_{D}(t,x,y)\nonumber\\
\ge&\int_{B(z^{t}_{y},t^{1/\alpha}/2)}\int_{B(z^{t}_{x},t^{1/\alpha}/2)}
p^{b,q_{0}}_{D}(t/3,x,u)p^{b,q_{0}}_{B(u,t^{1/\alpha}/2)\cup B(v,t^{1/\alpha}/2)}(t/3,u,v)p^{b,q_{0}}_{D}(t/3,v,y){\mathrm{d}u\mathrm{d}v}.\label{n4.33}
\end{align}
For $u\in B(z^{t}_{x},t^{1/\alpha}/2)$ and $v\in B(z^{t}_{y},t^{1/\alpha}/2)$,
\[|u-v|\le |u-z^{t}_{x}|+|z^{t}_{x}-x|+|x-y|+|y-z^{t}_{y}|+|z^{t}_{y}-v|\le |x-y|+10t^{1/\alpha}.\]
Thus for $x,y\in D$ with $|x-y|<4\eta(A)/5$, $|u-v|\le 4\eta(A)/5+5r^{*}/8\le 153\eta(A)/160$.
We now apply Lemma \ref{nlem4.9}.

Recall that $|u-v|\le |x-y|+10t^{1/\alpha}\le 2(|x-y|\vee 10t^{1/\alpha})$ for $u\in B(z^{t}_{x},t^{1/\alpha}/2)$ and $v\in B(z^{t}_{y},t^{1/\alpha}/2)$. If $|x-y|\ge 10t^{1/\alpha}$, then $|u-v|\le 2|x-y|$. Thus
\[t^{-d/\alpha}\wedge t J(u,v)
\gtrsim
\big(t^{-d/\alpha}\wedge t J(x,y)\big) \mbox{ and } t^{-d/\alpha}\wedge t J^{1}(u,v)
\gtrsim
\big(t^{-d/\alpha}\wedge t J^{1}(x,y)\big).\]
On the other hand, if $|x-y|<10t^{1/\alpha}$, then $|u-v|\le 20t^{1/\alpha}$. Thus
\[t^{-d/\alpha}\wedge t J(u,v)\asymp t^{-d/\alpha}\asymp t^{-d/\alpha}\wedge t J(x,y)\mbox{ and } t^{-d/\alpha}\wedge t J^{1}(u,v)\asymp t^{-d/\alpha} \asymp t^{-d/\alpha}\wedge t J^{1}(x,y).\]
These facts and Lemma \ref{nlem4.9} yield the following results:
there is a positive constant
$c_{1}=c_{1}(d,\allowbreak \alpha,\allowbreak \beta,\allowbreak A,\allowbreak M)$
such that for $t\in (0,T_{0}]$, $x,y\in D$ with $|x-y|<4\eta(A)/5$, $u\in B(z^{t}_{x},t^{1/\alpha}/2)$ and $v\in B(z^{t}_{y},t^{1/\alpha}/2)$,
\begin{equation}\label{7.38}
p^{b,q_{0}}_{B(u,t^{1/\alpha}/2)\cup B(v,t^{1/\alpha}/2)}(t/3,u,v)\ge
c_{1}
\big(t^{-d/\alpha}\wedge t J(x,y)\big).
\end{equation}
If, in addition,
\eqref{condi3} is satisfied for some $\eps>0$, there is a positive constant
$c_{2}=c_{2}(d,\allowbreak \alpha,\allowbreak \beta,\allowbreak A,\allowbreak M,\allowbreak \eps)$
such that for $t\in (0,T_{0}]$, $x,y\in D$, $u\in B(z^{t}_{x},t^{1/\alpha}/2)$ and $v\in B(z^{t}_{y},t^{1/\alpha}/2)$,
\begin{equation}\label{7.39}
p^{b,q_{0}}_{B(u,t^{1/\alpha}/2)\cup B(v,t^{1/\alpha}/2)}(t/3,u,v)\ge
c_{2}
\big(t^{-d/\alpha}\wedge t J(x,y)\big).
\end{equation}
Moreover, if \eqref{condi4} is satisfied for some $\eps>0$, then there is a positive constant
$c_{3}=c_{3}(d,\allowbreak \alpha,\allowbreak \beta,\allowbreak A,\allowbreak M,\allowbreak \eps)$
such that for $t\in (0,T_{0}]$, $x,y\in D$, $u\in B(z^{t}_{x},t^{1/\alpha}/2)$ and $v\in B(z^{t}_{y},t^{1/\alpha}/2)$,
\begin{equation}\label{7.40}
p^{b,q_{0}}_{B(u,t^{1/\alpha}/2)\cup B(v,t^{1/\alpha}/2)}(t/3,u,v)\ge
c_{3}
\big(t^{-d/\alpha}\wedge t J^{1}(x,y)\big).
\end{equation}
Hence by \eqref{n4.33} and \eqref{7.38}--\eqref{7.40}, to prove this lemma, it suffices to prove that
\begin{equation}
\int_{B(z^{t}_{x},t^{1/\alpha}/2)}p^{b,q_{0}}_{D}(t/3,x,u){\mathrm{d}u}\ge
c_{4}
\big(1\wedge\frac{\delta_{D}(x)^{\alpha/2}}{\sqrt{t}}\big),\label{n4.37}
\end{equation}
\begin{equation}
\int_{B(z^{t}_{y},t^{1/\alpha}/2)}p^{b,q_{0}}_{D}(t/3,v,y){\mathrm{d}v}\ge
c_{5}
\big(1\wedge\frac{\delta_{D}(y)^{\alpha/2}}{\sqrt{t}}\big)\label{n4.38}
\end{equation}
for some positive constants
$c_{i}=c_{i}(d,\allowbreak \alpha,\allowbreak \beta,\allowbreak R_{0},\allowbreak \Lambda_{0},\allowbreak A,\allowbreak M)$, $i=4,5$.
We prove only \eqref{n4.38}. The proof of \eqref{n4.37} is similar. If $\delta_{D}(y)\le 3 t^{1/\alpha}/2$, let $V_{1}:=D\cap B(z_{y},15t^{1/\alpha}/4)$; otherwise let $V_{1}:=B(y,3t^{1/\alpha}/4)$.
In both cases $V_{1}\subset D$ and $\mathrm{dist}(B(z^{t}_{y},t^{1/\alpha}/2),V_{1})\ge t^{1/\alpha}/4$. Indeed, in the first case $B(z^{t}_{y},t^{1/\alpha}/2)\subset B(z_{y},5t^{1/\alpha})\setminus B(z_{y},4t^{1/\alpha})$, while in the second case $|y-z^{t}_{y}|\ge 3t^{1/\alpha}/2$ by the construction of $z^{t}_{y}$.
By \eqref{lem7.3.3}, we have
\begin{align}
&\int_{B(z^{t}_{y},t^{1/\alpha}/2)}p^{b,q_{0}}_{D}(t/3,v,y){\mathrm{d}v}\nonumber\\
\ge&\frac{t}{3M}\einf_{\substack{z\in V_{1}\\ w\in B(z^{t}_{y},t^{1/\alpha}/2)}}j^{b}(w,z)
\P_{y}\big(\wh{\tau}^{b,q_{0}}_{V_{1}}>t/3\big)
\int_{B(z^{t}_{y},t^{1/\alpha}/4)}\P_{v}\big(\tau^{b,q_{0}}_{B(z^{t}_{y},t^{1/\alpha}/2)}>t/3\big){\mathrm{d}v}.\label{n4.39}
\end{align}
It follows from \eqref{nprop4.6.1} that
\begin{eqnarray}
&&\int_{B(z^{t}_{y},t^{1/\alpha}/4)}\P_{v}\big(\tau^{b,q_{0}}_{B(z^{t}_{y},t^{1/\alpha}/2)}>t/3\big){\mathrm{d}v}\nonumber\\
&\ge&\int_{B(z^{t}_{y},t^{1/\alpha}/4)}\int_{B(z^{t}_{y},t^{1/\alpha}/4)}p^{b,q_{0}}_{B(z^{t}_{y},t^{1/\alpha}/2)}(t/3,v,w)
\mathrm{d}w\mathrm{d}v\nonumber\\
&\ge&
c_{6}(d,\allowbreak \alpha,\allowbreak \beta,\allowbreak A)
t^{-d/\alpha}|B(z^{t}_{y},t^{1/\alpha}/4)|^{2}=
c_{7}(d,\allowbreak \alpha,\allowbreak \beta,\allowbreak A)
t^{d/\alpha}.\label{n4.40}
\end{eqnarray}
For $w\in B(z^{t}_{y},t^{1/\alpha}/2)$ and $z\in V_{1}$, $|w-z|\le 35 t^{1/\alpha}/4\le 35r^{*}/64<\eta(A)$.
Thus by \eqref{jb&j},
\begin{equation}
\einf_{\substack{z\in V_{1}\\ w\in B(z^{t}_{y},t^{1/\alpha}/2)}}j^{b}(w,z)\ge \frac{1}{2}\einf_{\substack{z\in V_{1}\\ w\in B(z^{t}_{y},t^{1/\alpha}/2)}}j(w,z)
\gtrsim
t^{-\frac{d+\alpha}{\alpha}}.\label{n4.41}
\end{equation}
By \eqref{n4.39}, \eqref{n4.40} and \eqref{n4.41}, to prove \eqref{n4.38}, it suffices to prove that
\begin{equation}
\P_{y}\big(\wh{\tau}^{b,q_{0}}_{V_{1}}>t/3\big)\ge
c_{8}
\big(1\wedge\frac{\delta_{D}(y)^{\alpha/2}}{\sqrt{t}}\big)
,
\label{n4.42}
\end{equation}
for some positive constant $c_{8}=c_{8}(d,\allowbreak \alpha,\allowbreak \beta,\allowbreak R_{0},\allowbreak \Lambda_{0},\allowbreak A,\allowbreak M)$.
If $\delta_{D}(y)>3t^{1/\alpha}/2$, it follows by \eqref{nprop4.6.2} that
\[\P_{y}\big(\wh{\tau}^{b,q_{0}}_{V_{1}}>t/3\big)=\int_{B(y,3t^{1/\alpha}/4)}\wh{p}^{b,q_{0}}_{B(y,3t^{1/\alpha}/4)}(t/3,y,z){\mathrm{d}z}
\ge \int_{B(y,3t^{1/\alpha}/8)}\wh{p}^{b,q_{0}}_{B(y,3t^{1/\alpha}/4)}(t/3,y,z){\mathrm{d}z}
\] is bounded from below by some positive constant
depending on $d,\alpha,\beta,A$ and $M$.
Thus \eqref{n4.42} holds in this case. Now we suppose $\delta_{D}(y)\le 3t^{1/\alpha}/2$ (in this case $V_{1}=D\cap B(z_{y},15t^{1/\alpha}/4)$).
Set $\theta_{0}:=1/8$ and $y_{0}:=z_{y}+\frac{29}{8}t^{1/\alpha}\mathbf{n}(z_{y})$. Since $15t^{1/\alpha}/4<r_{0}$, the ball $B(y_{0},\theta_{0}t^{1/\alpha})$ is contained in the interior ball $B(z_{y}+r_{0}\mathbf{n}(z_{y}),r_{0})\subset D$, and hence $B(y_{0},\theta_{0} t^{1/\alpha})\subset V_{1}\setminus B(z_{y},7t^{1/\alpha}/2)$.
Let $V_{2}:=V(z_{y},13t^{1/\alpha}/8)$ be the $C^{1,1}$ domain of Proposition \ref{P:2.2}(ii), so that $D\cap B(z_{y},13t^{1/\alpha}/8)\subset V_{2}\subset D\cap B(z_{y},13t^{1/\alpha}/4)$. Moreover, $y\in V_{2}$, since $|y-z_{y}|=\delta_{D}(y)\le 3t^{1/\alpha}/2$. By the Markov property of $\wh{X}^{b,q_{0}}$, we have
\begin{eqnarray}
&&\P_{y}\big(\wh{\tau}^{b,q_{0}}_{V_{1}}>t/3\big)\nonumber\\
&\ge&\P_{y}\big(\wh{\tau}^{b,q_{0}}_{V_{1}}>t/3;\ \wh{X}^{b,q_{0}}_{\wh{\tau}^{b,q_{0}}_{V_{2}}}\in B(y_{0},\theta_{0} t^{1/\alpha}/2)\big)\nonumber\\
&=&\E_{y}\big[\P_{\wh{X}^{b,q_{0}}_{\wh{\tau}^{b,q_{0}}_{V_{2}}}}\big(\wh{\tau}^{b,q_{0}}_{V_{1}}>t/3\big);\ \wh{X}^{b,q_{0}}_{\wh{\tau}^{b,q_{0}}_{V_{2}}}\in B(y_{0},\theta_{0} t^{1/\alpha}/2)\big]\nonumber\\
&\ge&\inf_{w\in B(y_{0},\theta_{0} t^{1/\alpha}/2)}\P_{w}\big(\wh{\tau}^{b,q_{0}}_{B(w,\theta_{0} t^{1/\alpha}/2)}>t/3\big)
\P_{y}\big(\wh{X}^{b,q_{0}}_{\wh{\tau}^{b,q_{0}}_{V_{2}}}\in B(y_{0},\theta_{0} t^{1/\alpha}/2)\big).\label{n4.43}
\end{eqnarray}
It follows by \eqref{nprop4.6.2} that $\inf_{w\in B(y_{0},\theta_{0} t^{1/\alpha}/2)}\P_{w}\big(\wh{\tau}^{b,q_{0}}_{B(w,\theta_{0} t^{1/\alpha}/2)}>t/3\big)$ is bounded from below by a positive constant
depending on $d,\alpha,\beta,A$ and $M$.
Moreover, $|z-w|\le 15t^{1/\alpha}/2<\eta(A)$ for $z\in B(y_{0},\theta_{0}t^{1/\alpha}/2)$ and $w\in V_{2}$.
By \eqref{poissonforwhxbq}, \eqref{n3.5}, \eqref{jb&j} and \eqref{n2.7},
we have
\begin{eqnarray}
\P_{y}\big(\wh{X}^{b,q_{0}}_{\wh{\tau}^{b,q_{0}}_{V_{2}}}\in B(y_{0},\theta_{0} t^{1/\alpha}/2)\big)
&=&\int_{B(y_{0},\theta_{0} t^{1/\alpha}/2)}\int_{V_{2}}G^{b,q_{0}}_{V_{2}}(w,y)j^{b}(z,w)\frac{h^{b,q_{0}}(z)}{h^{b,q_{0}}(y)}{\mathrm{d}w\mathrm{d}z}\nonumber\\
&\ge&\frac{1}{M}\int_{B(y_{0},\theta_{0} t^{1/\alpha}/2)}\int_{V_{2}}G^{b,q_{0}}_{V_{2}}(w,y)j^{b}(z,w){\mathrm{d}w\mathrm{d}z}\nonumber\\
&\ge&\frac{1}{8M}\int_{B(y_{0},\theta_{0} t^{1/\alpha}/2)}\int_{V_{2}}G_{V_{2}}(w,y)j(z,w){\mathrm{d}w\mathrm{d}z}\nonumber\\
&=&\frac{1}{8M}\P_{y}\big(X_{\tau_{V_{2}}}\in B(y_{0},\theta_{0} t^{1/\alpha}/2)\big).\nonumber
\end{eqnarray}
We choose a ball $B(y_{1},\theta_{1}t^{1/\alpha})\subset V_{2}$, where $\theta_{1}$ is a positive constant depending on $(R_{0},\Lambda_{0})$. By the boundary Harnack inequality for rotationally symmetric $\alpha$-stable process, we have
\begin{align}
\P_{y}\big(X_{\tau_{V_{2}}}\in B(y_{0},\theta_{0} t^{1/\alpha}/2)\big) \ge&
c_{9}(d,\allowbreak \alpha,\allowbreak \beta,\allowbreak R_{0},\allowbreak \Lambda_{0})
\P_{y_{1}}\big(X_{\tau_{V_{2}}}\in B(y_{0},\theta_{0} t^{1/\alpha}/2)\big)\frac{\delta_{D}(y)^{\alpha/2}}{\delta_{D}(y_{1})^{\alpha/2}}\nonumber\\
\ge&
c_{10}(d,\allowbreak \alpha,\allowbreak \beta,\allowbreak R_{0},\allowbreak \Lambda_{0})
\P_{y_{1}}\big(X_{\tau_{V_{2}}}\in B(y_{0},\theta_{0} t^{1/\alpha}/2)\big)\frac{\delta_{D}(y)^{\alpha/2}}{\sqrt{t}}.\label{n4.45}
\end{align}
For $w\in V_{2}$ and $z\in B(y_{0},\theta_{0}t^{1/\alpha}/2)$, $|w-z|\ge t^{1/\alpha}/4$. Thus
\begin{align}
\P_{y_{1}}\big(X_{\tau_{V_{2}}}\in B(y_{0},\theta_{0} t^{1/\alpha}/2)\big)
=&\int_{B(y_{0},\theta_{0}t^{1/\alpha}/2)}\int_{V_{2}}G_{V_{2}}(y_{1},w)j(w,z){\mathrm{d}w\mathrm{d}z}\nonumber\\
\ge&c_{11}(d,\allowbreak \alpha)
t^{-\frac{d+\alpha}{\alpha}}|B(y_{0},\theta_{0} t^{1/\alpha}/2)|\,\E_{y_{1}}[\tau_{V_{2}}]\nonumber\\
\ge&c_{12}(d,\allowbreak \alpha,\allowbreak R_{0},\allowbreak \Lambda_{0})
t^{-1}\E_{y_{1}}[\tau_{B(y_{1},\theta_{1}t^{1/\alpha})}]\nonumber\\
\ge&c_{13}(d,\allowbreak \alpha,\allowbreak R_{0},\allowbreak \Lambda_{0},\allowbreak T_{0}).\label{n4.46}
\end{align}
Consequently, \eqref{n4.42} follows from \eqref{n4.43}--\eqref{n4.46}. This completes the proof.\end{proof}

\begin{thm}\label{smalldist}
Suppose $D$ is a $C^{1,1}$ open set
with characteristics $(R_{0},\Lambda_{0})$.
For every $T\in (0,{\infty})$, there is a positive constant
$ C_{53}= C_{53}(d,\allowbreak \alpha,\allowbreak \beta,\allowbreak R_{0},\allowbreak \Lambda_{0},\allowbreak A,\allowbreak M,\allowbreak T)$
such that for every $t\in (0,T]$ and every $x,y\in D$ with $|x-y|<3\eta(A)/5$,
\begin{equation}\label{n4.47}
p^{b}_{D}(t,x,y)\ge  C_{53}\big(1\wedge \frac{\delta_{D}(x)^{\alpha/2}}{\sqrt{t}}\big)\big(1\wedge \frac{\delta_{D}(y)^{\alpha/2}}{\sqrt{t}}\big)\big(t^{-d/\alpha}\wedge \frac{t}{|x-y|^{d+\alpha}}\big).
\end{equation}
Moreover,
if \eqref{condi3} is satisfied for some $\eps>0$, then there is a positive constant
$ C_{54}= C_{54}(d,\allowbreak \alpha,\allowbreak \beta,\allowbreak R_{0},\allowbreak \Lambda_{0},\allowbreak A,\allowbreak M,\allowbreak T,\allowbreak \eps)$
such that for every $t\in (0,T]$ and every $x,y\in D$,
\begin{equation}\label{n4.48}
p^{b}_{D}(t,x,y)\ge  C_{54}\big(1\wedge \frac{\delta_{D}(x)^{\alpha/2}}{\sqrt{t}}\big)\big(1\wedge \frac{\delta_{D}(y)^{\alpha/2}}{\sqrt{t}}\big)\big(t^{-d/\alpha}\wedge \frac{t}{|x-y|^{d+\alpha}}\big).
\end{equation}
Furthermore, if
\eqref{condi4} is satisfied for some $\eps>0$, then there is a positive constant
$C_{55}= C_{55}(d,\allowbreak \alpha,\allowbreak \beta,\allowbreak R_{0},\allowbreak \Lambda_{0},\allowbreak A,\allowbreak M,\allowbreak T,\allowbreak \eps)$
such that for every $t\in (0,T]$ and every $x,y\in D$,
\begin{equation}\label{n4.49}
p^{b}_{D}(t,x,y)\ge  C_{55}\big(1\wedge \frac{\delta_{D}(x)^{\alpha/2}}{\sqrt{t}}\big)\big(1\wedge \frac{\delta_{D}(y)^{\alpha/2}}{\sqrt{t}}\big)\big(t^{-d/\alpha}\wedge\big( \frac{t}{|x-y|^{d+\alpha}}+\frac{t}{|x-y|^{d+\beta}}\big)\big).
\end{equation}
\end{thm}

\begin{proof}

Let $T_{0}$ and $r^{*}$ be the constants given in
\eqref{e:rstar},
and put $r:=T_{0}^{1/\alpha}=r^{*}/16$. Then $r\le \eta(A)/64$ and $r<R_{0}$.
Since $p^{b}_{D}\ge p^{b,q_{0}}_{D}$ by \eqref{n1.2}, estimate \eqref{n4.47} for $t\le T_{0}$ follows from \eqref{n4.30}.
Let $t\in (T_{0},T]$ and $x,y\in D$ with $|x-y|<3\eta(A)/5$.
Let $m=m(t)\in [2,T/T_{0}+1]$ be the integer such that $(m-1)T_{0}<t\le m T_{0}$, and put $s=s(t):=t/m\in [T_{0}/2,T_{0}]$.
By Proposition \ref{prop1}, there is a ball $B(\xi_{x},\theta_{0}r)\subset D\cap B(x,r)$. Set $K:=B(\xi_{x},\theta_{0}r/2)$.
For $z,z'\in K$ we have $\delta_{D}(z)\ge\theta_{0}r/2$, $|z-z'|<\theta_{0}r$, $|x-z|<r$ and $|z-y|<r+3\eta(A)/5<4\eta(A)/5$.
By the Chapman--Kolmogorov equation,
\[
p^{b}_{D}(t,x,y)\ge\int_{K^{m-1}}p^{b}_{D}(s,x,z_{1})\prod_{i=1}^{m-2}p^{b}_{D}(s,z_{i},z_{i+1})\,p^{b}_{D}(s,z_{m-1},y)\,{\mathrm{d}z}_{1}\cdots {\mathrm{d}z}_{m-1}.
\]
Consecutive points are at distance less than $4\eta(A)/5$, and $s\le T_{0}$. Hence \eqref{n4.30} and \eqref{n1.2} apply to each factor.
Since $s\in[T_{0}/2,T_{0}]$, $s\le t$, $\delta_{D}(z_{i})\ge\theta_{0}r/2$ and all distances involved are less than $\eta(A)$, there is $c_{1}=c_{1}(d,\allowbreak \alpha,\allowbreak \beta,\allowbreak R_{0},\allowbreak \Lambda_{0},\allowbreak A,\allowbreak M)>0$ such that
\[
p^{b}_{D}(s,x,z_{1})\ge c_{1}\big(1\wedge \frac{\delta_{D}(x)^{\alpha/2}}{\sqrt{t}}\big),\quad
p^{b}_{D}(s,z_{i},z_{i+1})\ge c_{1},\quad
p^{b}_{D}(s,z_{m-1},y)\ge c_{1}\big(1\wedge \frac{\delta_{D}(y)^{\alpha/2}}{\sqrt{t}}\big).
\]
Therefore
\[
p^{b}_{D}(t,x,y)\ge c_{1}^{m}|K|^{m-1}\big(1\wedge \frac{\delta_{D}(x)^{\alpha/2}}{\sqrt{t}}\big)\big(1\wedge \frac{\delta_{D}(y)^{\alpha/2}}{\sqrt{t}}\big).
\]
Since $m\le T/T_{0}+1$ and $t^{-d/\alpha}\wedge t|x-y|^{-d-\alpha}\le T_{0}^{-d/\alpha}$ for $t>T_{0}$, this proves \eqref{n4.47}.
The proofs of \eqref{n4.48} and \eqref{n4.49} are the same, with \eqref{n4.31} and \eqref{n4.32} in place of \eqref{n4.30}
, and
no restriction on $|x-y|$ is needed.
 For the last factor one uses
\begin{equation}\label{them6.10.n1}
s^{-d/\alpha}\wedge sJ(z_{m-1},y)\ge c_{2}\,(t^{-d/\alpha}\wedge tJ(x,y)),
\end{equation}
and similarly with $J^{1}$ in place of $J$, where $c_{2}$ depends on $T/T_{0}$.
Indeed, since $|z_{m-1}-y|\le |x-y|+r$, if $|x-y|\ge r$, then $|z_{m-1}-y|\le 2|x-y|$ and $J(z_{m-1},y)\ge 2^{-d-\alpha}J(x,y)$; if $|x-y|<r$, both sides of \eqref{them6.10.n1} are comparable to constants.
\end{proof}

\begin{thm}\label{largedist}

Suppose $D$ is a $C^{1,1}$ open set with characteristics $(R_{0},\Lambda_{0})$.
For every $\Lambda_{1}\ge 1$ and $T\in (0,{\infty})$, there exist positive constants
$C_{i}=C_{i}(d,\allowbreak \alpha,\allowbreak \beta,\allowbreak R_{0},\allowbreak \Lambda_{0},\allowbreak A,\allowbreak M,\allowbreak T,\allowbreak \Lambda_{1})$,
$i=56,57,58$ such that
the following holds. For every $t\in (0,T]$ and every $x,y$ in the same connected component of $D$ with $|x-y|\ge 3\eta(A)/5$
for which there is a rectifiable curve $l$ in $D$ connecting $x$ to $y$ with length no larger than $\Lambda_{1}|x-y|$, we have
\begin{equation}\label{n4.53}
p^{b}_{D}(t,x,y)\ge C_{56}\big(1\wedge \frac{\delta_{D}(x)^{\alpha/2}}{\sqrt{t}}\big)\big(1\wedge \frac{\delta_{D}(y)^{\alpha/2}}{\sqrt{t}}\big)\big(\frac{C_{57} t}{|x-y|}\big)^{C_{58}|x-y|}.
\end{equation}
In particular, if the path distance in each connected component of $D$ is comparable to the Euclidean distance with characteristic $\Lambda_{1}$, then \eqref{n4.53} holds for every $t\in(0,T]$ and every $x,y$ in the same connected component of $D$ with $|x-y|\ge 3\eta(A)/5$.
\end{thm}

\begin{proof}
Fix $\Lambda_{1}\ge 1$, $T>0$ and $x,y\in D$ with $|x-y|\ge 3\eta(A)/5$ such that there is a rectifiable curve $l$ in $D$ connecting $x$ to $y$
whose length $d_{D}(x,y)$ satisfies $d_{D}(x,y)\le \Lambda_{1}|x-y|$. Let $N:=3\eta(A)/5$ and $R:=\Lambda_{1}|x-y|$.
Let $n\ge 2$ be the integer such that $R/n\le N/2< R/(n-1)$, so that in this case, $R/n>N(1-1/n)/2\ge N/4$. Let $x_{0}=x,x_{1},\cdots,x_{n-1},x_{n}=y$ be points on $l$ that divide $l$ into $n$ arcs of equal length.
Then $|x_{i}-x_{i+1}|\le d_{D}(x,y)/n\le R/n$ for every $i\in \{0,\cdots,n-1\}$. Let $r:=(R_{0}\wedge N)/16$. By Proposition \ref{prop1}, for each $i\in \{1,\cdots,n-1\}$, there is a ball $B(\xi_{i},\theta_{0}r)\subset D\cap B(x_{i},r)$, where $\theta_{0}=\theta_{0}(\Lambda_{0})\in (0,1)$. Let $B_{i}:=B(\xi_{i},\theta_{0}r/2)$. By the Chapman--Kolmogorov equation, we have
\begin{equation}\label{n4.54}
p^{b}_{D}(t,x,y)\ge \int_{B_{1}\times B_{2}\times\cdots\times B_{n-1}}p^{b}_{D}(t/n,x,z_{1})p^{b}_{D}(t/n,z_{1},z_{2})\cdots
p^{b}_{D}(t/n,z_{n-1},y){\mathrm{d}z}_{1}{\mathrm{d}z}_{2}\cdots {\mathrm{d}z}_{n-1}.
\end{equation}
For $z_{1}\in B_{1}$ and $z_{n-1}\in B_{n-1}$, we have $|x-z_{1}|\vee |z_{n-1}-y|\le r+\frac{R}{n}\le 9N/16$.
For $z_{i}\in B_{i}$ and $z_{i+1}\in B_{i+1}$, $i\in \{1,\cdots,n-2\}$, we have
$\delta_{D}(z_{i})\ge \theta_{0}r/2$ and $|z_{i+1}-z_{i}|\le 2r+R/n\le 5N/8$. Thus by \eqref{n4.47}, for every $t\in (0,T]$ and $z_{1}\in B_{1}$,
\begin{eqnarray}
p^{b}_{D}(t/n,x,z_{1})&\ge&c_{1}\big(1\wedge\frac{\delta_{D}(x)^{\alpha/2}}{\sqrt{t/n}}
\big)\big(1\wedge\frac{\delta_{D}(z_{1})^{\alpha/2}}{\sqrt{t/n}}
\big)\big(\big(\frac{t}{n}\big)^{-d/\alpha}\wedge \frac{t/n}{|x-z_{1}|^{d+\alpha}}\big)\nonumber\\
&\ge&c_{1}\big(1\wedge \frac{\delta_{D}(x)^{\alpha/2}}{\sqrt{t/2}}\big)\big(1\wedge \frac{(\theta_{0}r/2)^{\alpha/2}}{\sqrt{T/2}}\big)\big(\big(\frac{T}{2}\big)^{-d/\alpha}\wedge \frac{t/n}{(9N/16)^{d+\alpha}}\big)\nonumber\\
&\ge&c_{2}\big(1\wedge \frac{\delta_{D}(x)^{\alpha/2}}{\sqrt{t}}\big)\big(1\wedge \frac{t}{n}\big).
 \nonumber
\end{eqnarray}
Similarly, \eqref{n4.47} and the above estimates give
\begin{equation}\label{n4.56}
p^{b}_{D}(t/n,z_{i},z_{i+1})\ge c_{3}\big(1\wedge \frac{t}{n}\big)\quad\mbox{ for }i\in \{1,\cdots,n-2\},
\end{equation}
\begin{equation}\label{n4.57}
p^{b}_{D}(t/n,z_{n-1},y)\ge c_{4}\big(1\wedge \frac{\delta_{D}(y)^{\alpha/2}}{\sqrt{t}}\big)\big(1\wedge\frac{t}{n}\big).
\end{equation}
Here
$c_{i}=c_{i}(d,\allowbreak \alpha,\allowbreak \beta,\allowbreak R_{0},\allowbreak \Lambda_{0},\allowbreak A,\allowbreak M,\allowbreak T)$,
$i=1,\cdots,4$ are positive constants, and $|B_{i}|=|B_{1}|$ for every $i$. By \eqref{n4.54}--\eqref{n4.57} and the facts that $2R/N\le n<(2R/N)+1<3R/N$ and $R=\Lambda_{1}|x-y|$, we have
\begin{eqnarray}
p^{b}_{D}(t,x,y)&\ge&c_{2}c_{3}^{n-2}c_{4}\prod_{i=1}^{n-1}|B_{i}|\big(1\wedge\frac{t}{n}\big)^{n}\big(1\wedge \frac{\delta_{D}(x)^{\alpha/2}}{\sqrt{t}}\big)\big(1\wedge \frac{\delta_{D}(y)^{\alpha/2}}{\sqrt{t}}\big)\nonumber\\
&\ge&c_{2}c_{3}^{-2}c_{4}|B_{1}|^{-1}\big(c_{3}|B_{1}|\big)^{n}\big(1\wedge\frac{t}{n}\big)^{n}\big(1\wedge \frac{\delta_{D}(x)^{\alpha/2}}{\sqrt{t}}\big)\big(1\wedge \frac{\delta_{D}(y)^{\alpha/2}}{\sqrt{t}}\big)\nonumber\\
&\ge&c_{6}\big(c_{7}\big(1\wedge\frac{t}{|x-y|}\big)\big)^{c_{8}|x-y|}\big(1\wedge \frac{\delta_{D}(x)^{\alpha/2}}{\sqrt{t}}\big)\big(1\wedge \frac{\delta_{D}(y)^{\alpha/2}}{\sqrt{t}}\big).\label{n4.58}
\end{eqnarray}
Since $|x-y|\ge N$ and $t\le T$, we have $1\wedge \frac{t}{|x-y|}\ge \big(1\wedge\frac{N}{T}\big)\frac{t}{|x-y|}$. This and \eqref{n4.58} imply \eqref{n4.53}. This completes the proof.\end{proof}

\section{
Small-time heat kernel estimates: the case $d=1\le \alpha$
}\label{sec9}

In this section we
prove Theorem \ref{T:1.2}(i)--(iii) for $d=1$. We deduce them from the case $d=2$, treated in Section \ref{sec7}, by viewing the one-dimensional process as the first coordinate of a process on $\R^{2}$. Only the range $\alpha\in[1,2)$ is new, since the case $d=1>\alpha$ is covered by Sections \ref{sec5}--\ref{sec7}. The argument, however, works for every $\alpha\in(0,2)$.

We write $\tilde{x}=(x_{1},x')\in\R\times\R$ for points of $\R^{2}$ and $\pi(\tilde{x}):=x_{1}$. Quantities in dimension $d$ carry the subscript or argument $d$ when necessary, e.g.\ $\mathcal{A}(d,-\alpha)$, $p^{(d)}_{a}$ and $\eta_{d}(A)$. Recall that
$\mathcal{A}(d,-\gamma)=\gamma 2^{\gamma-1}\pi^{-d/2}\Gamma((d+\gamma)/2)\Gamma(1-\gamma/2)^{-1}$ for $\gamma\in(0,2)$. The substitution $z'=|z_{1}|u$ and the Beta integral give $\int_{\R}(z_{1}^{2}+|z'|^{2})^{-\frac{2+\gamma}{2}}\mathrm{d}z'=\sqrt{\pi}\Gamma((1+\gamma)/2)\Gamma((2+\gamma)/2)^{-1}|z_{1}|^{-1-\gamma}$
for $z_{1}\not=0$.
Since $\mathcal{A}(2,-\gamma)/\mathcal{A}(1,-\gamma)=\pi^{-1/2}\Gamma((2+\gamma)/2)/\Gamma((1+\gamma)/2)$, we obtain
\begin{equation}\label{e:proj}
\int_{\R}\mathcal{A}(2,-\gamma)\big(z_{1}^{2}+|z'|^{2}\big)^{-\frac{2+\gamma}{2}}\mathrm{d}z'=\mathcal{A}(1,-\gamma)|z_{1}|^{-1-\gamma},\quad z_{1}\not=0.
\end{equation}
Let $a\ge 0$. For $d=2$ and $b\equiv a$, the process $X^{b}$ is a L\'{e}vy process on $\R^{2}$ with characteristic exponent $|\xi|^{\alpha}+a|\xi|^{\beta}$. Evaluating its characteristic function at $\xi=(\xi_{1},0)$ shows that the first coordinate of $X^{b}$ is a L\'evy process on $\R$ with characteristic exponent $|\xi_{1}|^{\alpha}+a|\xi_{1}|^{\beta}$.
Hence, for $t>0$ and $\tilde{x},\tilde{y}\in\R^{2}$,
\begin{equation}\label{e:projpa}
\int_{\R}p^{(2)}_{a}\big(t,\tilde{x},(y_{1},y')\big)\mathrm{d}y'=p^{(1)}_{a}(t,x_{1},y_{1}).
\end{equation}
For $d\ge 1$, $a\ge 0$ and $x\not= y$, set $J^{a}_{d}(x,y):=|x-y|^{-d-\alpha}+a|x-y|^{-d-\beta}$.

\begin{lemma}\label{L:9.1}
Let $d\ge 1$ and $A',T>0$. There is $c=c(d,\alpha,\beta,A',T)\ge 1$ such that
$p^{(d)}_{a}(t,x,y)\stackrel{c}{\asymp}t^{-d/\alpha}\wedge tJ^{a}_{d}(x,y)$ for $a\in[0,A']$, $t\in(0,T]$ and $x,y\in\R^{d}$.
\end{lemma}

\begin{proof}
This is \cite[(1.10)]{CW} with $k:=TA'^{\alpha/(\alpha-\beta)}$, since then $(0,T]\subset(0,ka^{-\alpha/(\alpha-\beta)}]$ for $a\in(0,A']$; for $a=0$ it is \eqref{p0}.
\end{proof}

Let $A>0$ and $b\in\mathcal{B}_{A}$ on $\R\times\R$. Let $\varphi(w):=\frac{1}{2}\1_{\{|w|<1\}}$ on $\R$, and for $\tilde{x},\tilde{z}\in\R^{2}$
define
\begin{align*}
\tilde{b}_{1}(\tilde{x},\tilde{z}):=&\frac{|\tilde{z}|^{2+\beta}}{\mathcal{A}(2,-\beta)}\Big(\mathcal{A}(1,-\beta)b^{+}(x_{1},z_{1})|z_{1}|^{-2-\beta}\varphi\big(z'/|z_{1}|\big)\\
&\qquad-\frac{\mathcal{A}(1,-\beta)\mathcal{A}(2,-\alpha)}{\mathcal{A}(1,-\alpha)}\,b^{-}(x_{1},z_{1})|z_{1}|^{\alpha-\beta}|\tilde{z}|^{-2-\alpha}\Big)
\1_{\{z_{1}\not=0\}},
\end{align*}
and
\[
\tilde{b}_{2}(\tilde{x},\tilde{z}):=\frac{|\tilde{z}|^{2+\beta}}{\mathcal{A}(2,-\beta)}\Big(\rho(x_{1},z_{1})\big(\mathcal{A}(2,-\alpha)|\tilde{z}|^{-2-\alpha}+\mathcal{A}(2,-\beta)
|\tilde{z}|^{-2-\beta}\big)-\mathcal{A}(2,-\alpha)|\tilde{z}|^{-2-\alpha}\Big)
\1_{\{z_{1}\not=0\}}
,
\]
where
$b^{-}:=(-b)\vee 0$ and
$\rho(x_{1},z_{1}):=j^{b}(x_{1},x_{1}+z_{1})\big/\big(\mathcal{A}(1,-\alpha)|z_{1}|^{-1-\alpha}+\mathcal{A}(1,-\beta)|z_{1}|^{-1-\beta}\big)$.
The lift $\tilde{b}_{1}$ spreads $b^{+}$ over the cone $\{|z'|<|z_{1}|\}$. This keeps the large jumps under control, which is needed for the upper bounds. The lift $\tilde{b}_{2}$ is used only under \eqref{condi4}; it preserves the $\beta$-stable tail of the jump kernel in every direction.

\begin{lemma}\label{L:9.2}
There is a constant $c_{*}=c_{*}(\alpha,\beta)\ge 1$ such that the following hold for every $b\in\mathcal{B}_{A}$ on $\R\times\R$.
\begin{enumerate}[label=\upshape(\roman*), ref=\upshape(\roman*)]
\item $\tilde{b}_{1},\tilde{b}_{2}\in\mathcal{B}_{c_{*}(1+A)}$ on $\R^{2}\times\R^{2}$.
\item For every $\tilde{x}\in\R^{2}$ and
every $z_{1}\in\R$ with $z_{1}\not=0$,
$$\int_{\R}j^{\tilde{b}_{1}}\big(\tilde{x},\tilde{x}+(z_{1},z')\big)\mathrm{d}z'=\int_{\R}j^{\tilde{b}_{2}}\big(\tilde{x},\tilde{x}+(z_{1},z')\big)\mathrm{d}z'=j^{b}(x_{1},x_{1}+z_{1}).$$
\item $M_{\tilde{b}_{1}^{+},\sqrt{2}\lambda}\le c_{*}M_{b^{+},\lambda}$ for every $\lambda>0$.
\item If $b$ satisfies \eqref{condi5} (resp.\ \eqref{condi3}) with $\eps\in(0,1]$, then $\tilde{b}_{1}$ satisfies \eqref{condi5} (resp.\ \eqref{condi3}) in
$\R^{2}\times\R^{2}$
with $\eps/c_{*}$ in place of $\eps$. If $b$ satisfies \eqref{condi4} with $\eps>0$, then $\tilde{b}_{2}$ satisfies \eqref{condi4} in
$\R^{2}\times\R^{2}$
with $\eps/c_{*}$ in place of $\eps$.
\end{enumerate}
\end{lemma}

\begin{proof}
Write $\kappa_{\gamma}(\tilde{z}):=\mathcal{A}(2,-\gamma)|\tilde{z}|^{-2-\gamma}$ and $k_{\gamma}(z_{1}):=\mathcal{A}(1,-\gamma)|z_{1}|^{-1-\gamma}$ for $\gamma\in\{\alpha,\beta\}$. By definition, for $z_{1}\not=0$,
\begin{align}
j^{\tilde{b}_{1}}(\tilde{x},\tilde{x}+\tilde{z})&=\kappa_{\alpha}(\tilde{z})\Big(1-\frac{\mathcal{A}(1,-\beta)}{\mathcal{A}(1,-\alpha)}b^{-}(x_{1},z_{1})|z_{1}|^{\alpha-\beta}\Big)+\mathcal{A}(1,-\beta)b^{+}(x_{1},z_{1})|z_{1}|^{-2-\beta}\varphi(z'/|z_{1}|),\label{e:jlift1}\\
j^{\tilde{b}_{2}}(\tilde{x},\tilde{x}+\tilde{z})&=\rho(x_{1},z_{1})\big(\kappa_{\alpha}(\tilde{z})+\kappa_{\beta}(\tilde{z})\big).\label{e:jlift2}
\end{align}
Since $b^{+}b^{-}=0$, the first expression equals $\kappa_{\alpha}(\tilde{z})\,j^{b}(x_{1},x_{1}+z_{1})/k_{\alpha}(z_{1})$ where $b(x_{1},z_{1})<0$, and is at least $\kappa_{\alpha}(\tilde{z})$ where $b(x_{1},z_{1})\ge 0$. Since $j^{b}\ge 0$, both kernels are non-negative, i.e.\ $\tilde{b}_{1}$ and $\tilde{b}_{2}$ satisfy \eqref{condi2}. Both are even in $\tilde{z}$, because $b(x_{1},\cdot)$ is even and $\varphi$ is radial.

\smallskip

\noindent(i) On the support of $\varphi(z'/|z_{1}|)$ we have $|z_{1}|\le|\tilde{z}|<\sqrt{2}|z_{1}|$, so the first term of $\tilde{b}_{1}$ is at most $2^{\beta/2}\frac{\mathcal{A}(1,-\beta)}{\mathcal{A}(2,-\beta)}A$. The second term is bounded by $\frac{\mathcal{A}(1,-\beta)\mathcal{A}(2,-\alpha)}{\mathcal{A}(1,-\alpha)\mathcal{A}(2,-\beta)}A$, since $|z_{1}|\le|\tilde{z}|$. For $\tilde{b}_{2}$ we have
$\tilde{b}_{2}(\tilde{x},\tilde{z})=\rho(x_{1},z_{1})+(\rho(x_{1},z_{1})-1)\frac{\mathcal{A}(2,-\alpha)}{\mathcal{A}(2,-\beta)}|\tilde{z}|^{\beta-\alpha}$, $0\le\rho(x_{1},z_{1})\le 1+A$, and
$|\rho(x_{1},z_{1})-1|\,|\tilde{z}|^{\beta-\alpha}\le\mathcal{A}(1,-\beta)|b(x_{1},z_{1})-1|\,|z_{1}|^{-1-\beta}k_{\alpha}(z_{1})^{-1}|\tilde{z}|^{\beta-\alpha}\le\frac{\mathcal{A}(1,-\beta)}{\mathcal{A}(1,-\alpha)}(1+A)$.

\smallskip

\noindent(ii) This follows from \eqref{e:jlift1}, \eqref{e:jlift2}, \eqref{e:proj} with $\gamma=\alpha,\beta$, and $\int_{\R}|z_{1}|^{-1}\varphi(z'/|z_{1}|)\mathrm{d}z'=1$.

\smallskip

\noindent(iii) When $\tilde{b}_{1}^{+}>0$ we have $b(x_{1},z_{1})>0$ and $|z'|<|z_{1}|$. If moreover $|\tilde{z}|>\sqrt{2}\lambda$, then $|z_{1}|>\lambda$, and $\tilde{b}_{1}^{+}(\tilde{x},\tilde{z})\le 2^{\beta/2}\frac{\mathcal{A}(1,-\beta)}{\mathcal{A}(2,-\beta)}b^{+}(x_{1},z_{1})$.

\smallskip

\noindent(iv) Assume \eqref{condi5}. When $b<0$, \eqref{e:jlift1} gives $j^{\tilde{b}_{1}}(\tilde{x},\tilde{x}+\tilde{z})\le\kappa_{\alpha}(\tilde{z})$. When $b\ge 0$,  we have $\mathcal{A}(1,-\beta)b^{+}(x_{1},z_{1})|z_{1}|^{-1-\beta}\le j^{b}(x_{1},x_{1}+z_{1})\le\eps^{-1}|z_{1}|^{-1-\alpha}$, so the last term in \eqref{e:jlift1} is at most $\frac{1}{2}\eps^{-1}|z_{1}|^{-2-\alpha}\le 2^{\alpha/2}\eps^{-1}|\tilde{z}|^{-2-\alpha}$ on its support.
Assume \eqref{condi3}. When $b\ge 0$, $j^{\tilde{b}_{1}}(\tilde{x},\tilde{x}+\tilde{z})\ge\kappa_{\alpha}(\tilde{z})$. When $b<0$,
 $j^{\tilde{b}_{1}}(\tilde{x},\tilde{x}+\tilde{z})=\kappa_{\alpha}(\tilde{z})j^{b}(x_{1},x_{1}+z_{1})/k_{\alpha}(z_{1})\ge\eps\,\kappa_{\alpha}(\tilde{z})/\mathcal{A}(1,-\alpha)$.
Assume \eqref{condi4}. Then $\rho\ge\eps/(\mathcal{A}(1,-\alpha)\vee\mathcal{A}(1,-\beta))$, and \eqref{e:jlift2} gives the assertion.
\end{proof}

\begin{lemma}\label{L:9.3}
Let $b\in\mathcal{B}_{A}$ on $\R\times\R$, and let $\hat{b}\in\{\tilde{b}_{1},\tilde{b}_{2}\}$. Let $\tilde{X}:=(X^{\hat{b}}_{t})_{t\ge 0}$ be the Feller process on $\R^{2}$ associated with
$\mathcal{L}^{\hat{b}}$
, and let $\tilde{x}=(x_{1},x')\in\R^{2}$. Then $
\pi(\tilde{X}):=
(\pi(\tilde{X}_{t}))_{t\ge 0}$ under $\P_{\tilde{x}}$ has the same law as
$X^{b}$
under $\P_{x_{1}}$.
Consequently, for every open set $D\subset\R$ and $\tilde{D}:=D\times\R$,
\begin{equation}\label{e:projpD}
p^{b}_{D}(t,x_{1},y_{1})=\int_{\R}p^{\hat{b}}_{\tilde{D}}\big(t,(x_{1},x'),(y_{1},y')\big)\mathrm{d}y'\quad\mbox{for }t>0\mbox{ and }x_{1},y_{1}\in D.
\end{equation}
\end{lemma}

\begin{proof}
Let $f\in\mathcal{S}(\R)$ and $g:=f\circ\pi\in C^{2}_{b}(\R^{2})$. By \cite[Theorem 1.1(iv)]{CW} and the Markov property, $g(\tilde{X}_{t})-g(\tilde{X}_{0})-\int_{0}^{t}\L^{\hat{b}}g(\tilde{X}_{s})\mathrm{d}s$ is a $\P_{\tilde{x}}$-martingale.
Since $j^{\hat{b}}(\tilde{x},\tilde{x}+\tilde{z})$ is even in $\tilde{z}$,
\[
\L^{\hat{b}}g(\tilde{x})=\frac{1}{2}\int_{\R^{2}}\big(g(\tilde{x}+\tilde{z})+g(\tilde{x}-\tilde{z})-2g(\tilde{x})\big)j^{\hat{b}}(\tilde{x},\tilde{x}+\tilde{z})\mathrm{d}\tilde{z},
\]
and the integral converges absolutely. The integrand equals $\big(f(x_{1}+z_{1})+f(x_{1}-z_{1})-2f(x_{1})\big)j^{\hat{b}}(\tilde{x},\tilde{x}+\tilde{z})$. By the Fubini theorem and Lemma \ref{L:9.2}(ii), $\L^{\hat{b}}g(\tilde{x})=\L^{b}f(x_{1})$.
Hence $\pi(\tilde{X})$ solves the martingale problem for $(\L^{b},\mathcal{S}(\R))$ with initial value $x_{1}$. Its paths are c\`adl\`ag. By the uniqueness in \cite[Theorem 1.3]{CW} (see also \cite[Theorem 1.6]{CW}), its law is that of $X^{b}$ under $\P_{x_{1}}$.

Since $\tau^{\hat{b}}_{\tilde{D}}=\inf\{t>0:\pi(\tilde{X}_{t})\notin D\}$, we get $\P_{x_{1}}(X^{b}_{t}\in E,\,t<\tau^{b}_{D})=\P_{\tilde{x}}(\tilde{X}_{t}\in E\times\R,\,t<\tau^{\hat{b}}_{\tilde{D}})$ for every Borel set $E\subset D$. Hence \eqref{e:projpD} holds for a.e.\ $y_{1}\in D$. Both sides are continuous in $y_{1}$. For the left side this is Theorem \ref{them3.4}. For the right side it follows from Theorem \ref{them3.4} in $\R^{2}$ and the generalized dominated convergence theorem, since
$p^{\hat{b}}_{\tilde{D}}(t,\cdot,\cdot)\le c\,\e^{ct}p^{(2)}_{A'}(t,\cdot,\cdot)$ for some $c=c(\alpha,\beta,A')$, where $A'=c_{*}(1+A)$ (see \eqref{ltdensity} and the remark preceding \eqref{n2.1}), and the integrals of the dominating functions, $c\,\e^{ct}p^{(1)}_{A'}(t,x_{1},y_{1})$ by \eqref{e:projpa}, are continuous in $y_{1}$. Therefore \eqref{e:projpD} holds for every $y_{1}\in D$.
\end{proof}

Let $D\subset\R$ be a $C^{1,1}$ open set with characteristics $(R_{0},\Lambda_{0})$. Then $\tilde{D}=D\times\R$ is a $C^{1,1}$ open set in $\R^{2}$ with characteristics $(R_{0},1)$. Indeed, its boundary is a union of lines
$\{(y_{1},y')\in\R^{2}:\ y_{1}=c\}$
at mutual distance at least $R_{0}$, so for $Q\in\partial\tilde{D}$ the set $\tilde{D}\cap B(Q,R_{0})$ is a half-disc. Moreover, $\delta_{\tilde{D}}(y_{1},y')=\delta_{D}(y_{1})$, and every connected component of $\tilde{D}$ has the form $I\times\R$ for an open interval $I$, hence is convex.

\begin{thm}\label{T:9.4}
Theorem \ref{T:1.2}{\rm(i)}--{\rm(iii)} hold for $d=1$ and every
$0<\beta<\alpha<2$,
with $R_{1}:=\frac{3}{10}\eta_{2}(c_{*}(1+A))$. In {\rm(ii)}, the path distance condition is automatically satisfied with $\Lambda_{1}=1$, since the components of $D$ are intervals.
\end{thm}

\begin{proof}
Fix $A,T>0$, $b\in\mathcal{B}_{A}$ on $\R\times\R$, $x_{1},y_{1}\in D$ and $t\in(0,T]$. Put $A':=c_{*}(1+A)$, $\tilde{x}:=(x_{1},0)$ and $\tilde{y}:=(y_{1},y')$, and note that $\Phi_{D}(t,x_{1},y_{1})=\Phi_{\tilde{D}}(t,\tilde{x},\tilde{y})$. Let $R_{1}^{(2)}:=3\eta_{2}(A')/5$ be the constant of Theorem \ref{T:1.2}(ii) for $d=2$ and $A'$, so that $R_{1}=R_{1}^{(2)}/2$. All constants below depend only on $\alpha,\beta,R_{0},\Lambda_{0},A,T$ and, where indicated, on $\lambda$ or $\eps$. Since \eqref{condi5}, \eqref{condi3} and \eqref{condi4} with $\eps$ imply the same conditions with $\eps\wedge 1$, we may assume $\eps\in(0,1]$.

\smallskip

\noindent(i) Let $\lambda>0$.
By Lemma \ref{L:9.2}(iii), $a:=M_{\tilde{b}_{1}^{+},\sqrt{2}\lambda}\le c_{*}M_{b^{+},\lambda}\le A'$.
Applying Theorem \ref{T:1.2}(i) for $d=2$ to $\tilde{b}_{1}\in\mathcal{B}_{A'}$, $\tilde{D}$ and $\sqrt{2}\lambda$, we get
$p^{\tilde{b}_{1}}_{\tilde{D}}(t,\tilde{x},\tilde{y})\le c_{0}\Phi_{D}(t,x_{1},y_{1})\big(t^{-2/\alpha}\wedge t J^{a}_{2}(\tilde{x},\tilde{y})\big)$.
By Lemma \ref{L:9.1}, \eqref{e:projpD} and \eqref{e:projpa},
\[
p^{b}_{D}(t,x_{1},y_{1})
=\int_{\R}p^{\tilde{b}_{1}}_{\tilde{D}}(t,\tilde{x},\tilde{y})\mathrm{d}y'
\le c_{1}\Phi_{D}(t,x_{1},y_{1})\int_{\R}p^{(2)}_{a}(t,\tilde{x},\tilde{y})\mathrm{d}y'=c_{1}\Phi_{D}(t,x_{1},y_{1})\,p^{(1)}_{a}(t,x_{1},y_{1}),
\]
and by Lemma \ref{L:9.1} again, $p^{(1)}_{a}(t,x_{1},y_{1})\le c_{2}\big(t^{-1/\alpha}\wedge tJ^{c_{*}M_{b^{+},\lambda}}_{1}(x_{1},y_{1})\big)\le c_{2}c_{*}\big(t^{-1/\alpha}\wedge tJ^{M_{b^{+},\lambda}}_{1}(x_{1},y_{1})\big)$. This is the first assertion of (i). Under \eqref{condi5}, Lemma \ref{L:9.2}(iv) and Theorem \ref{T:1.2}(i) for $d=2$ give the bound with $a=0$, and the same computation gives the second assertion.

\smallskip

\noindent(iii) Under \eqref{condi3}, Lemma \ref{L:9.2}(iv), Theorem \ref{T:1.2}(iii) for $d=2$, Lemma \ref{L:9.1}, \eqref{e:projpD} and \eqref{e:projpa} give $p^{b}_{D}(t,x_{1},y_{1})\ge c_{3}\Phi_{D}(t,x_{1},y_{1})p^{(1)}_{0}(t,x_{1},y_{1})\ge c_{4}f^{0}_{D}(t,x_{1},y_{1})$. Under \eqref{condi4} the same argument with $\tilde{b}_{2}$ in place of $\tilde{b}_{1}$ gives the bound with $f^{1}_{D}$.

\smallskip

\noindent(ii) First let $|x_{1}-y_{1}|<R_{1}$, and put $s:=t^{1/\alpha}\vee|x_{1}-y_{1}|$ and $\varrho:=R_{1}\wedge s$. For $|y'|<\varrho$ we have $|\tilde{x}-\tilde{y}|<2R_{1}=R^{(2)}_{1}$ and $|\tilde{x}-\tilde{y}|\le 2s$. By Theorem \ref{T:1.2}(ii) for $d=2$ and \eqref{e:projpD},
\[
p^{b}_{D}(t,x_{1},y_{1})\ge c_{5}\Phi_{D}(t,x_{1},y_{1})\int_{|y'|<\varrho}\big(t^{-2/\alpha}\wedge t(2s)^{-2-\alpha}\big)\mathrm{d}y'\ge c_{6}\Phi_{D}(t,x_{1},y_{1})\,\varrho\big(t^{-2/\alpha}\wedge ts^{-2-\alpha}\big).
\]
Since $t\le T$, we have $\varrho\ge c_{7}s$ with $c_{7}:=1\wedge (R_{1}T^{-1/\alpha})$. If $|x_{1}-y_{1}|\le t^{1/\alpha}$, then $s=t^{1/\alpha}$ and the right side is at least $
c_{6}c_{7}
\,\Phi_{D}(t,x_{1},y_{1})\, t^{-1/\alpha}$. Otherwise $s=|x_{1}-y_{1}|$ and it is at least $
c_{6}c_{7}
\,\Phi_{D}(t,x_{1},y_{1})\, t|x_{1}-y_{1}|^{-1-\alpha}$. In both cases $p^{b}_{D}(t,x_{1},y_{1})\ge c_{8}f^{0}_{D}(t,x_{1},y_{1})$.

Next let $x_{1},y_{1}$ lie in the same component of $D$ with $r:=|x_{1}-y_{1}|\ge R_{1}$. For $y'$ in the set $S:=\{R_{1}^{(2)}\le|y'|\le R_{1}^{(2)}+1\}$, the points $\tilde{x},\tilde{y}$ lie in the same convex component of $\tilde{D}$, $|\tilde{x}-\tilde{y}|\ge R^{(2)}_{1}$, and $r\le|\tilde{x}-\tilde{y}|\le\kappa r$ with $\kappa:=3+2/R^{(2)}_{1}$. By \eqref{them1.ii} for $d=2$ with $\Lambda_{1}=1$, $p^{\tilde{b}_{1}}_{\tilde{D}}(t,\tilde{x},\tilde{y})\ge c_{9}\Phi_{D}(t,x_{1},y_{1})(c_{10}t/\tilde{r})^{c_{11}\tilde{r}}$ with $\tilde{r}:=|\tilde{x}-\tilde{y}|$. We claim that $(c_{10}t/\tilde{r})^{c_{11}\tilde{r}}\ge c_{12}(c_{10}t/(\kappa r))^{c_{11}\kappa r}$. If $c_{10}t\le\tilde{r}$, this holds with $c_{12}=1$, because
$c_{10}t/(\kappa r)\le c_{10}t/\tilde{r}\le 1$
and $\tilde{r}\le\kappa r$. If $c_{10}t>\tilde{r}$, the left side is at least $1$, while
$r<c_{10}t\le c_{10}T$
and hence the right side without $c_{12}$ is at most $\big(1\vee\frac{c_{10}T}{\kappa R_{1}}\big)^{c_{11}\kappa c_{10}T}=:c_{12}^{-1}$. Integrating over $S$ and using \eqref{e:projpD}, we get \eqref{them1.ii} for $d=1$ with constants depending on $\alpha,\beta,R_{0},\Lambda_{0},A,T$.
\end{proof}

\section{Large-time heat kernel estimates}\label{sec8}

In this section $d\ge 1$ and $A>0$. Recall the definition of $\eta (A)$ from \eqref{defiforetaA}.
Set $R_{1}:=3\eta(A)/5$ if $d>\alpha$, and let $R_{1}$ be the constant of Theorem \ref{T:9.4} if $d=1\le\alpha$.
The proofs below use the small-time upper bound \eqref{7.6.1}, the small-time lower bounds \eqref{n4.47} and \eqref{n4.48}, and Theorems \ref{smalldist}
and \ref{largedist}. For $d=1\le\alpha$, the corresponding estimates are given by Theorem \ref{T:9.4}. To keep the notation simple, we cite only the results of Section \ref{sec7}; when $d=1\le \alpha$, they are to be replaced by the corresponding parts of Theorem \ref{T:9.4}.

For large time, the decay of $p^{b}_{D}(t,x,y)$ is governed by the principal eigenvalue of $\L^{b, D}$.
Since $\L^{b, D}$ is not self-adjoint, we use the spectral theory of positive compact operators.
Let $D$ be an arbitrary bounded open set in $\R^{d}$, and let $b\in\mathcal{B}_{A}$.
We first recall some facts from spectral theory.
For a linear operator $\A$ on $L^{2}(D;\mathrm{d}x)$ with domain $\mathrm{Dom}(\A)$, its resolvent set $\rho(\A)$ is the set of $\lambda\in\mathbb{C}$ such that $(\lambda I-\A)^{-1}$ exists as a bounded operator on $L^{2}(D;\mathrm{d}x)$, and its spectrum is $\sigma(\A):=\mathbb{C}\setminus\rho(\A)$. It is known that $\rho (\A)$ is an open set in $\mathbb{C}$.

 Define
\[
P^{b,D}_{t}f(x):=\int_{D}f(y)p^{b}_{D}(t,x,y)\mathrm{d}y,\quad   f\in L^{2}(D;{\mathrm{d}x}).
\]
Since for every $t>0$, $(x,y)\mapsto p^{b}_{D}(t,x,y)$ is bounded on
$D\times D$ (by
Proposition \ref{P:2.1}(i)) and $D$ is bounded,
it follows that
\[
\int_{D\times D} p^b_D(t, x, y)^2\, {\mathrm{d}x}\, {\mathrm{d}y} \le |D|^{2}\sup_{x,y\in D}p^{b}_{D}(t,x,y)^{2} <\infty.
\]
So for each $t>0$,
$P^{b,D}_{t}$ is a Hilbert--Schmidt operator, and hence compact.
Thus by the Riesz--Schauder theorem, $\sigma (P^{b, D}_t)$ is a countable set
 whose only possible accumulation point is $0$, and each non-zero
  $\lambda \in \sigma (P^{b, D}_{t})$
 is  an eigenvalue of finite multiplicity.
Let $\mathcal{L}^{b,D}$ be the infinitesimal generator of the semigroup $P^{b,D}_{t}$ on $L^{2}(D;\mathrm{d}x)$.
We use $(\cdot,\cdot)$ and $\|\cdot\|_{2}$ to denote the inner product
and norm in $L^2 (D; {\mathrm{d}x})$, respectively.

\begin{thm} \label{T:7.1}
Let $\lambda_{0}$ and $C_{30}$ be the constants in Lemma \ref{lem:exit}. The following hold.
\begin{enumerate}[label=\upshape(\roman*), ref=\upshape(\roman*)]
\item The spectrum $\sigma (\L^{b, D})$ is a non-empty discrete set consisting of eigenvalues of finite multiplicity, without finite accumulation points, and
 $\lambda^{b, D}_{1}:= - \sup {\rm Re} \, \sigma (\L^{b, D})\geq \lambda_0$. Moreover, $-\lambda^{b,D}_{1}$ is an eigenvalue of $\L^{b,D}$, and there is a non-negative $\phi\in L^{2}(D;\mathrm{d}x)$ with $\|\phi\|_{2}=1$ such that
\[P^{b,D}_{t} \phi=\e^{-t\lambda^{b,D}_{1}} \phi \quad\mbox{ for every }t>0.\]
 We call $-\lambda^{b,D}_{1}$ the principal eigenvalue of $\L^{b,D}$.
\item If, in addition,
\begin{equation}\label{e:irred}
p^{b}_{D}(t,x,y)>0 \quad\hbox{for every } t>0 \hbox{ and } x,y\in D,
\end{equation}
then
 the non-negative, normalized eigenfunction $\phi$ in (i) is unique, and is positive and continuous on $D$.
Moreover,
 $-\lambda^{b, D}_1   > {\rm Re}\,\mu$ for any other $\mu \in \sigma (\L^{b, D})$.
\end{enumerate}
\end{thm}

\begin{remark}\label{R:irred}
Condition \eqref{e:irred} holds if any two connected components of $D$ are at distance less than $\eta(A)$; in particular, if $D$ is connected. This follows from Lemma \ref{cor3.6}.
It also holds if \eqref{condi3} is satisfied, by \eqref{n4.48}, and if $\mathrm{diam}(D)<R_{1}$, by \eqref{n4.47}.
It cannot be dropped in
part (ii) of Theorem \ref{T:7.1}.
If $j^{b}$ vanishes for large jumps and $D$ has two components farther apart than the range of the jumps, then $X^{b,D}$ never moves between them.
An eigenfunction for the larger of the two principal eigenvalues then vanishes on one component.
\end{remark}

\begin{proof}

\medskip

\noindent(i)
Since each $P^{b, D}_t$ is compact, each resolvent operator $(\lambda I- \L^{b, D})^{-1}$ with $\lambda \in \rho (\L^{b, D})$
 is compact (cf.
 \cite[Theorem 2.3.3]{P}
 ).  Fix some $\lambda \in \rho (\L^{b, D})$. By the Riesz--Schauder theorem,
 $\sigma ((\lambda -\L^{b, D})^{-1})$ is a countable set
 whose only possible accumulation point is $0$, and
 each non-zero point in $\sigma( (\lambda -\L^{b, D})^{-1})$ is
 an eigenvalue of finite multiplicity.  It follows that $\sigma (\L^{b, D})$ is a discrete set
 consisting of eigenvalues without finite accumulation points, and each eigenvalue is of finite multiplicity.

 We also know by
 \cite[Theorem 2.2.4]{P}
 that
 \begin{equation}\label{e:7.3}
 \e^{t \sigma (\L^{b, D})} \subset \sigma (P^{b, D}_t) \subset
 \e^{t \sigma (\L^{b, D})} \cup \{0\} .
 \end{equation}
Let $r_{t}$ denote the spectral radius of $P^{b,D}_{t}$
 as an operator on $L^{2}(D;\mathrm{d}x)$, and let $\|\cdot\|_{L^{2}\to L^{2}}$ denote the corresponding operator norm.
By \eqref{e:7.3}, $r_{t}=\e^{-t\lambda^{b,D}_{1}}$ for every $t>0$, with $\e^{-\infty}:=0$.
We shall show that $r_{t}\le \e^{-\lambda_{0}t}$
where $\lambda_{0}$ is the constant given in
Lemma \ref{lem:exit}.
By
Proposition \ref{P:2.1}(i),
$c_{1}:=\sup_{x,y\in D}p^{b}_{D}(1,x,y)<\infty$.
For $f\in L^{2}(D;{\mathrm{d}x})$, $s>0$ and $x\in D$, the Markov property, the Cauchy--Schwarz inequality and
Lemma \ref{lem:exit}
give
\[
|P^{b,D}_{s+1}f(x)|\le \|P^{b,D}_{1}f\|_{\infty}\,\P_{x}(\tau^{b}_{D}>s)\le c_{1}|D|^{1/2}\|f\|_{2}\,C_{30}\e^{-\lambda_{0}s}.
\]
Hence $\|P^{b,D}_{s+1}\|_{L^{2}\to L^{2}}\le c_{1}C_{30}\e^{\lambda_{0}}|D|\,\e^{-\lambda_{0}(s+1)}$.
For $t>0$ and $n>1/t$, apply this bound with $s+1=nt$. Gelfand's formula then gives
\[
r_{t}=\lim_{n\to\infty}\|P^{b,D}_{nt}\|_{L^{2}\to L^{2}}^{1/n}\le \e^{-\lambda_{0}t}.
\]
Therefore $\lambda^{b,D}_{1}\ge\lambda_{0}$.

We next show that $\sigma(\L^{b,D})\neq\emptyset$. Let $K$ be a closed ball contained in $D$. Since $K$ lies in one connected component of $D$, Theorem \ref{them3.4} and Lemma \ref{cor3.6} give
$c_{K}:=\inf_{x,y\in K}p^{b}_{D}(1,x,y)>0$. Restricting all integrations in the Chapman--Kolmogorov equation to $K$ gives $(P^{b,D}_{n}\1_{K},\1_{K})\ge c_{K}^{n}|K|^{n+1}$ for $n\ge 1$. Hence $\|P^{b,D}_{n}\|_{L^{2}\to L^{2}}\ge (c_{K}|K|)^{n}$, and
$r_{1}\ge c_{K}|K|>0$. By \eqref{e:7.3}, $\sigma(\L^{b,D})\neq\emptyset$, so $\lambda^{b,D}_{1}<\infty$.

Since $r_{t}=\e^{-t\lambda^{b,D}_{1}}$ for every $t>0$, the growth bound of the semigroup $P^{b,D}_{t}$ equals $-\lambda^{b,D}_{1}$, which is also the spectral bound of $\L^{b,D}$. The semigroup is positive on the Banach lattice $L^{2}(D;\mathrm{d}x)$, so its spectral bound belongs to $\sigma(\L^{b,D})$; see \cite[Theorem C-III.1.1]{AGG86}. Hence $-\lambda^{b,D}_{1}$ is an eigenvalue of $\L^{b,D}$. Fix $\lambda>-\lambda^{b,D}_{1}$. Then $R_{\lambda}:=(\lambda-\L^{b,D})^{-1}=\int_{0}^{\infty}\e^{-\lambda t}P^{b,D}_{t}\,\mathrm{d}t$ is a positive compact operator. Its spectral radius is $1/\mathrm{dist}(\lambda,\sigma(\L^{b,D}))=1/(\lambda+\lambda^{b,D}_{1})$, because $-\lambda^{b,D}_{1}\in\sigma(\L^{b,D})$ and $|\lambda-\mu|\ge\lambda-\mathrm{Re}\,\mu\ge\lambda+\lambda^{b,D}_{1}$ for every $\mu\in\sigma(\L^{b,D})$. By the Krein--Rutman theorem (see \cite[Chapter V]{Schaefer}), there is a non-negative $\phi\in L^{2}(D;\mathrm{d}x)$ with
$\| \phi \|_{2}=1$ and $R_{\lambda} \phi =(\lambda+\lambda^{b,D}_{1})^{-1} \phi $. Thus $ \phi \in\mathrm{Dom}(\L^{b,D})$ and $\L^{b,D} \phi =-\lambda^{b,D}_{1} \phi $, which gives $P^{b,D}_{t} \phi =\e^{-t\lambda^{b,D}_{1}} \phi $ for every $t>0$.

\medskip

\noindent(ii)
Now assume \eqref{e:irred} and fix $t>0$.
 Since $P^{b,D}_{t}$ is compact and, by \eqref{e:irred}, has a strictly positive kernel,
 by \cite[Proposition V.6.6]{Schaefer},
 its spectral radius
 $r_{t}$
 is an eigenvalue of $P^{b, D}_t$,
  and has
 a unique eigenfunction $\phi^{(t)}$ with unit $L^2$-norm
 and $\phi^{(t)}>0$ a.e. on $D$.  Moreover,
     every other eigenvalue $\lambda$ of $P^{b,D}_{t}$ satisfies
 $| \lambda | < r_t$.

Let $\phi:=\phi^{(1)}$ and $\lambda_{1}:=\lambda^{b,D}_{1}$, so that $r_{1}=\e^{-\lambda_{1}}$ by \eqref{e:7.3}.
For $t>0$, $P^{b,D}_{1}P^{b,D}_{t}\phi=P^{b,D}_{t}P^{b,D}_{1}\phi=r_{1}P^{b,D}_{t}\phi$, and $P^{b,D}_{t}\phi$ is non-negative and not identically zero by \eqref{e:irred}.
By the uniqueness of the normalized eigenfunction for $r_{1}$, $P^{b,D}_{t}\phi=c(t)\phi$ for some $c(t)>0$.
The semigroup property gives $c(t+s)=c(t)c(s)$, and $t\mapsto c(t)=(P^{b,D}_{t}\phi,\phi)$ is continuous on $(0,\infty)$ by Theorem \ref{them3.4}, \eqref{ltdensity} and dominated convergence.
Hence $c(t)=\e^{-\lambda t}$ for some $\lambda\in\R$, and $c(1)=r_{1}$ gives $\lambda=\lambda_{1}$.
Thus $P^{b,D}_{t}\phi=\e^{-\lambda_{1}t}\phi$ for every $t>0$.
 Consequently, $ \phi \in\mathrm{Dom}(\L^{b,D})$ and $\L^{b,D} \phi =-\lambda_{1} \phi $.
Moreover, the non-negative normalized eigenfunction for $-\lambda_{1}$ is unique since the eigenspace of $P^{b,D}_{1}$ corresponding to $r_{1}$ is one-dimensional by \cite[Proposition V.6.6]{Schaefer}.
  Furthermore, Theorem \ref{them3.4} together with
Proposition \ref{P:2.1}(i) implies that $\phi = \e^{\lambda_1} P^{b, D}_1 \phi$   is a bounded positive continuous function on $D$.
 Suppose  $\mu\in\sigma(\L^{b,D})$ satisfies $\mathrm{Re}\,\mu=-\lambda_{1}$ and $\mu\not=-\lambda_{1}$.
Then ${\rm Im}\,\mu\neq 0$ and $\e^{t\mu}\neq r_{t}$ for some $t>0$.
By \eqref{e:7.3}, $\e^{t\mu}$ is a non-zero point of $\sigma(P^{b,D}_{t})$, hence an eigenvalue, and $|\e^{t\mu}|=r_{t}$. This contradicts the strict inequality $|\lambda|<r_{t}$ established above for the other eigenvalues of $P^{b,D}_{t}$.
Hence $-\lambda_{1}>{\rm Re}\,\mu$ for every other $\mu\in\sigma(\L^{b,D})$.
 \end{proof}

\begin{thm}\label{T:7.2}
Suppose $D$ is a bounded $C^{1,1}$ open set in
$\R^{d}$
with characteristics $(R_{0},\Lambda_{0})$
and $A,T\in (0,\infty)$.
Let $\lambda_0>0$ and $\lambda^{b, D}_1\geq \lambda_0 $ be the constants given in Theorem \ref{T:7.1}.
Then
there is a constant $C_{59}= C_{59}(d,\alpha,\beta,R_{0},\Lambda_{0},\mathrm{diam}(D),A,M,T)>0$
so that for every
$b\in\mathcal{B}_{A}$,
we have
\begin{equation}\label{e:7.5}
p^{b}_{D}(t,x,y)\le  C_{59} \e^{-t\lambda_0}
\delta_{D}(x)^{\alpha/2}\delta_{D}(y)^{\alpha/2}\quad \mbox{ for }(t,x,y)\in [T,\infty)\times D\times D.
\end{equation}
If, in addition, \eqref{e:irred} holds,
then there is a constant $ C_{60}= C_{60}(d,\allowbreak \alpha,\allowbreak \beta,\allowbreak D,\allowbreak A,\allowbreak M,\allowbreak b,\allowbreak T)>0$ such that
\begin{equation}\label{e:7.6}
p^{b}_{D}(t,x,y)\le  C_{60} \e^{-t \lambda^{b, D}_1 }
\delta_{D}(x)^{\alpha/2}\delta_{D}(y)^{\alpha/2}\quad \mbox{ for }(t,x,y)\in [T,\infty)\times D\times D.
\end{equation}
\end{thm}

\begin{proof}  Without loss of generality,
we assume $T=1$.
First, for $t>1$, we have by the Chapman--Kolmogorov equation,  \eqref{7.6.1}
 and
 Lemma \ref{lem:exit},
  \begin{eqnarray*}
 p^{b}_{D}(t, x, y)&=& \int_{D\times D} p^{b}_{D}(1/2, x, z) p^{b}_{D}(t-1, z, w) p^{b}_{D}(1/2, w, y) {\mathrm{d}z} {\mathrm{d}w} \\
 &\leq&   c_{1}
 (1\wedge \delta_D (x))^{\alpha /2}  (1\wedge \delta_D (y))^{\alpha /2}
 \int_{D\times D} p^{b}_{D}(t-1, z, w)  {\mathrm{d}z\mathrm{d}w} \\
   &=&c_{1}
 (1\wedge \delta_D (x))^{\alpha /2}  (1\wedge \delta_D (y))^{\alpha /2}
 \int_{D} \P_{z}\big(\tau^{b}_{D}>t-1\big) {\mathrm{d}z} \\
   &\leq&   c_{1} C_{30}
 (1\wedge \delta_D (x))^{\alpha /2}  (1\wedge \delta_D (y))^{\alpha /2}
 \e^{-\lambda_0 (t-1)} |D| ,
    \end{eqnarray*}
where
$c_{1}=c_{1}(d,\allowbreak \alpha,\allowbreak \beta,\allowbreak R_{0},\allowbreak \Lambda_{0},\allowbreak A,\allowbreak M )>0$.
This proves \eqref{e:7.5}.

For the proof of \eqref{e:7.6}, assume \eqref{e:irred}, and let $\phi$ be the positive eigenfunction in Theorem \ref{T:7.1}(ii).
Let $r:=R_{1}\wedge (R_{0}/2)$.
Since $r<R_{0}$,
by Proposition \ref{prop1} there is a constant $\theta_{0}=\theta_{0}(\Lambda_{0}) \in (0, 1)$ such that for every $x\in D$, there is a point $\xi(x)$ with $B(\xi(x), \theta_{0} r)
\subset B(x, r) \cap D$.
For notational simplicity, we write
$\lambda_1$ for $\lambda^{b, D}_1$ in the rest of this proof.
 Since $|x-y|<r\le R_{1}$ for all $y\in B(x,r)\cap D$, by
 \eqref{n4.47},
 \begin{eqnarray}
\phi (x) &=&  \e^{\lambda_1 } P^{b, D}_1 \phi (x)\nonumber\\
&\geq &
c_{2}\e^{\lambda_{1}}
(1\wedge \delta_D (x))^{\alpha /2} \int_{B(x, r)\cap D}
 (1\wedge \delta_D (y))^{\alpha /2}  \big( 1\wedge \frac{1}{|x-y|^{d+\alpha}} \big)
   \phi (y) {\mathrm{d}y}\nonumber\\
&\geq &
c_{2}\e^{\lambda_{1}}
(1\wedge \delta_D (x))^{\alpha /2} (1\wedge r^{-(d+\alpha)} )
 \int_{B(\xi(x), \theta_{0} r)}  (1\wedge \delta_D (y))^{\alpha /2} \phi (y) {\mathrm{d}y}\nonumber\\
 &\geq &
 c_{3}\e^{\lambda_{1}}
 (1\wedge \delta_D (x))^{\alpha /2}.\label{e:7.4}
\end{eqnarray}
  Here
  $c_{2}=c_{2}(d,\allowbreak \alpha,\allowbreak \beta,\allowbreak R_{0},\allowbreak \Lambda_{0},\allowbreak A,\allowbreak M)>0$
  and $c_{3}=c_{3}(d,\allowbreak \alpha,\allowbreak \beta,\allowbreak D,\allowbreak A,\allowbreak M,\allowbreak b)>0$.
 The last inequality is due to the fact that $v(z):=  \int_{B(z, \theta_{0} r/2)}  (1\wedge \delta_D (y))^{\alpha /2} \phi (y) {\mathrm{d}y}$
 is a positive continuous function on the compact set $\{z\in D: \delta_D(z)\geq \theta_{0} r/2 \}$ and its minimum there is strictly positive.
 For $t>1$, by the Chapman--Kolmogorov equation,  \eqref{7.6.1},
 \eqref{e:7.4} and Theorem \ref{T:7.1},
 \begin{eqnarray*}
 p^{b}_{D}(t, x, y)&=& \int_{D\times D} p^{b}_{D}(1/2, x, z) p^{b}_{D}(t-1, z, w) p^{b}_{D}(1/2, w, y) {\mathrm{d}z} {\mathrm{d}w} \\
 &\leq&   c_{4}
 (1\wedge \delta_D (x))^{\alpha /2}  (1\wedge \delta_D (y))^{\alpha /2}
 \int_{D\times D} p^{b}_{D}(t-1, z, w) (1\wedge \delta_D (w))^{\alpha /2} {\mathrm{d}z\mathrm{d}w} \\
 &\leq&   c_{4}c_{3}^{-1}
  \e^{-\lambda_{1}}
 (1\wedge \delta_D (x))^{\alpha /2}  (1\wedge \delta_D (y))^{\alpha /2}
 \int_{D\times D} p^{b}_{D}(t-1, z, w) \phi (w) {\mathrm{d}z\mathrm{d}w} \\
&=&
c_{4}c_{3}^{-1} \e^{-\lambda_{1}}
(1\wedge \delta_D (x))^{\alpha /2}  (1\wedge \delta_D (y))^{\alpha /2}
 \int_{D } \e^{-\lambda_1 (t-1)} \phi (z) {\mathrm{d}z} \\
 &\le&c_{4}c_{3}^{-1} \e^{-\lambda_{1}t}|D|^{1/2}\|\phi\|_{2}
(1\wedge \delta_D (x))^{\alpha /2}  (1\wedge \delta_D (y))^{\alpha /2}\nonumber\\
 &\leq &
 c_{5}c^{-1}_{3}\e^{-\lambda_1 t} (1\wedge \delta_D (x))^{\alpha /2}  (1\wedge \delta_D (y))^{\alpha /2}.
 \end{eqnarray*}
 Here
 $c_{4}=c_{4}(d,\allowbreak \alpha,\allowbreak \beta,\allowbreak R_{0},\allowbreak \Lambda_{0},\allowbreak A,\allowbreak M)>0$ and $c_{5}=c_{5}(d,\alpha,\beta,R_{0},\Lambda_{0},\mathrm{diam}(D),A,M)>0$.
 This establishes \eqref{e:7.6}.
 \end{proof}

\begin{thm}\label{ltthem1}
Suppose $D$ is a bounded $C^{1,1}$ open set
with characteristics $(R_{0},\Lambda_{0})$
and
$b\in\mathcal{B}_{A}$.
Assume one of the following conditions holds:
\begin{enumerate}[label=\upshape(\roman*), ref=\upshape(\roman*)]
\item $\mathrm{diam}(D)< R_{1}$, where $R_{1}=3\eta(A)/5$ if $d>\alpha$, and $R_{1}$ is the constant of Theorem \ref{T:9.4} if $d=1\le\alpha$;
\item $D$ is connected and the path distance in $D$ is comparable to the Euclidean distance with characteristic $\Lambda_{1}\ge 1$;
\item \eqref{condi3} is satisfied for some $\eps>0$.
\end{enumerate}
Then for every $T\in (0,\infty)$, there exists a constant
$ C_{61}\ge 1$ depending on
$d,\allowbreak\alpha,\allowbreak\beta,\allowbreak R_{0},\allowbreak\Lambda_{0},\allowbreak\mathrm{diam}(D),\allowbreak A,\allowbreak M,\allowbreak T$, (if {\rm(ii)} holds) $\Lambda_{1}$,
and (if {\rm(iii)} holds) $\eps$,
such that for every $t\in [T,\infty)$ and every $x,y\in D$,
\[
C_{61}^{-1} \e^{-t\lambda^{b,D}_{1}}\delta_{D}(x)^{\alpha/2}\delta_{D}(y)^{\alpha/2}
\le p^{b}_{D}(t,x,y)
 \le  C_{61} \e^{-t\lambda^{b,D}_{1}}\delta_{D}(x)^{\alpha/2}\delta_{D}(y)^{\alpha/2}.
 \]
Here $\lambda^{b,D}_{1}$ is as in Theorem \ref{T:7.1}.
  \end{thm}

\begin{proof}
 For notational simplicity, we write
$\lambda_1$ for $\lambda^{b, D}_1$ in this proof.
Under (i), condition \eqref{e:irred} follows from \eqref{n4.47}, since $|x-y|\le\mathrm{diam}(D)<R_{1}$ for $x,y\in D$ and $T>0$ is arbitrary. Under (ii) and (iii), it holds by Remark \ref{R:irred}.
 Let $\phi$ be the positive eigenfunction
 given in Theorem \ref{T:7.1}(ii).
By \eqref{7.6.1}, the Cauchy--Schwarz inequality and $\|\phi\|_{2}=1$, we have for every $x\in D$,
\begin{equation}\label{lteq2}
\phi(x)=\e^{\lambda_{1}}\int_{D}\phi(y)p^{b}_{D}(1,x,y){\mathrm{d}y}\le c_{1}\e^{\lambda_{1}}|D|^{1/2}(1\wedge \delta_{D}(x)^{\alpha/2})=:c_{2}\e^{\lambda_{1}}(1\wedge \delta_{D}(x)^{\alpha/2}),
\end{equation}
where $c_{1}=c_{1}(d,\allowbreak \alpha,\allowbreak \beta,\allowbreak R_{0},\allowbreak \Lambda_{0},\allowbreak A,\allowbreak M)>0$ and $c_{2}=c_{2}(d,\alpha,\beta,R_{0},\Lambda_{0},\mathrm{diam}(D),A,M)>0$.
By Theorem \ref{smalldist} (applied with $|x-y|\le\mathrm{diam}(D)$ under (i)), Theorem \ref{largedist} and the boundedness of $D$,
under each of the assumptions (i)--(iii),
\begin{equation}\label{8.9}
p^{b}_{D}(1,x,y)\ge c_{3}\big(1\wedge \delta_{D}(x)^{\alpha/2}\big)\big(1\wedge \delta_{D}(y)^{\alpha/2}\big)\quad\mbox{ for }x,y\in D,
\end{equation}
where $c_{3}>0$ depends on
$d,\allowbreak\alpha,\allowbreak\beta,\allowbreak R_{0},\allowbreak\Lambda_{0},\allowbreak\mathrm{diam}(D),\allowbreak A,\allowbreak M$, (if {\rm(ii)} holds) $\Lambda_{1}$
and (if {\rm(iii)} holds) $\eps$.
Set $I:=\int_{D}(1\wedge \delta_{D}(y)^{\alpha/2})\phi(y){\mathrm{d}y}$. By \eqref{lteq2} and $\|\phi\|_{2}=1$, $I\ge c_{2}^{-1}\e^{-\lambda_{1}}$. Hence, by \eqref{8.9},
\begin{equation}\label{7.8}
\phi(x)=\e^{\lambda_{1}}\int_{D}\phi(y)p^{b}_{D}(1,x,y){\mathrm{d}y}\ge c_{3}\e^{\lambda_{1}}\big(1\wedge \delta_{D}(x)^{\alpha/2}\big)I\ge c_{3}c_{2}^{-1}\big(1\wedge \delta_{D}(x)^{\alpha/2}\big).
\end{equation}
By \eqref{8.9} and \eqref{7.8},
\[
\e^{-\lambda_{1}}=\int_{D}\int_{D}\phi(x)p^{b}_{D}(1,x,y)\phi(y){\mathrm{d}x\mathrm{d}y}\ge c_{3}I^{2}\ge c_{3}^{3}c_{2}^{-2}\Big(\int_{D}\big(1\wedge \delta_{D}(y)^{\alpha/2}\big)^{2}{\mathrm{d}y}\Big)^{2}.
\]
By Proposition \ref{P:2.2}(iii), $D$ contains a ball
whose radius depends only on $R_{0}$ and $\Lambda_{0}$,
so that the last integral is bounded below by a positive constant depending on $d,\alpha,R_{0},\Lambda_{0}$. Since $\lambda_{1}\ge\lambda_{0}>0$, we obtain
\begin{equation}\label{7.7}
1\le \e^{\lambda_{1}}\le c_{4}
\end{equation}
for a constant $c_{4}\ge 1$ with the same dependence as $c_{3}$. By \eqref{lteq2}, \eqref{7.8} and \eqref{7.7}, $\phi(x)\stackrel{c_{5}}{\asymp}1\wedge \delta_{D}(x)^{\alpha/2}$ on $D$. Since
$\int_{D}\int_{D}\phi(x)p^{b}_{D}(t,x,y)\phi(y){\mathrm{d}x\mathrm{d}y}=\e^{-\lambda_{1}t}$ for every $t>0$, it follows that
\begin{equation}
c_{5}^{-2}\e^{-\lambda_{1}t}\le
\int_{D}\int_{D}(1\wedge
\delta_{D}(x)^{\alpha/2})p^{b}_{D}(t,x,y)(1\wedge
\delta_{D}(y)^{\alpha/2}){\mathrm{d}x\mathrm{d}y}
\le c_{5}^{2}\e^{-\lambda_{1}t}\quad\mbox{ for every }t>0.
\label{lteq3}
\end{equation}
 By the boundedness of $D$, Theorem \ref{stupbound}, Theorem \ref{smalldist} and Theorem
\ref{largedist},
under each of the assumptions (i)--(iii),
for every $T>0$ and every $x,y\in D$,
\begin{equation}\label{8.14}
p^{b}_{D}(T/4,x,y)\stackrel{c_{6}}{\asymp}
\big(1\wedge \delta_{D}(x)^{\alpha/2}\big)\big(1\wedge \delta_{D}(y)^{\alpha/2}\big)
\end{equation}
for some positive constant $c_{6}$
depending on $d,\allowbreak\alpha,\allowbreak\beta,\allowbreak R_{0},\allowbreak\Lambda_{0},\allowbreak\mathrm{diam}(D),\allowbreak A,\allowbreak M,\allowbreak T$, (if {\rm(ii)} holds) $\Lambda_{1}$, and (if {\rm(iii)} holds) $\eps$.
 By the Chapman--Kolmogorov equation, \eqref{8.14},
\eqref{lteq3} and \eqref{7.7},
we have
\begin{align}
&p^{b}_{D}(t,x,y)\nonumber\\
=&\int_{D}\int_{D}p^{b}_{D}(T/4,x,z)p^{b}_{D}(t-T/2,z,w)p^{b}_{D}(T/4,w,y)\mathrm{d}w\mathrm{d}z\nonumber\\
\stackrel{c_{7}}{\asymp}&
(1\wedge \delta_{D}(x)^{\alpha/2})(1\wedge
\delta_{D}(y)^{\alpha/2})\int_{D}\int_{D}(1\wedge
\delta_{D}(z)^{\alpha/2})p^{b}_{D}(t-T/2,z,w)(1\wedge
\delta_{D}(w)^{\alpha/2}){\mathrm{d}z\mathrm{d}w}\nonumber\\
\stackrel{c_{8}}{\asymp}&
\e^{-\lambda_{1}(t-\frac{T}{2})}(1\wedge
\delta_{D}(x)^{\alpha/2})(1\wedge \delta_{D}(y)^{\alpha/2})\nonumber\\
\stackrel{c_{9}}{\asymp}&
\e^{-\lambda_{1}t}(1\wedge
\delta_{D}(x)^{\alpha/2})(1\wedge \delta_{D}(y)^{\alpha/2})\nonumber.
\end{align}
 Here $c_{i},\ i=7,8,9$ are positive constants
with the same dependence as $c_{6}$.
This completes the proof. \end{proof}

\begin{proof}[Proof of Theorem \ref{T:1.2}] Let first $d>\alpha$, and set $R_{1}:=3\eta(A)/5$. The constant $M$ in \eqref{defiforM} depends only on $d,\alpha,\beta$ and $A$.
Assertion (i) is Theorem \ref{stupbound}.
Assertion (ii) follows from \eqref{n4.47} and Theorem \ref{largedist}, and assertion (iii) from \eqref{n4.48} and \eqref{n4.49}.
For $d=1\le\alpha$, assertions (i)--(iii) are Theorem \ref{T:9.4}.
Finally, for every $d\ge 1$, assertion (iv) follows from Theorems \ref{T:7.1}, \ref{T:7.2} and \ref{ltthem1}. Each of the conditions (a)--(c) implies \eqref{e:irred}, as shown in the proof of Theorem \ref{ltthem1}.
\end{proof}

\bigskip

 \textbf{Acknowledgment.}
  This work was carried out over a period of more than ten years. While it is just finalized now (November 23, 2026), a preliminary version was made available on arXiv in March 2015.
  After a complete draft of the manuscript has been prepared, we
used Claude and Kimi (K3.1 Pro)
to proofread the text, polish the presentation, and check our mathematical arguments.
  The original manuscript contained full proofs for the case $d>\alpha$, but treated the case $d=1\le \alpha$ only briefly. During this review, Claude identified that the extension to $d=1\le\alpha$ required additional justification.
  The final version now gives a full proof for this case
by lifting the one-dimensional process $X^{b}$   to  dimension 2 and viewing  it as the first coordinate of a suitable process on $\R^{2}$.
 The authors reviewed and verified all mathematical arguments and take full responsibility for the content of this paper.

We thank Renming Song for the discussion on spectral theory of non-symmetric operators.

\bigskip

\noindent {\bf Zhen-Qing Chen}

\noindent  Department of Mathematics, University of Washington, Seattle,
WA 98195, USA.

\noindent  E-mail: zqchen@uw.edu

\bigskip
\noindent  {\bf Ting Yang}

\noindent School of Mathematics and Statistics, Beijing Institute of Technology,  Beijing, 100081,
 P. R. China.

\noindent Beijing Key Laboratory on MCAACI, Beijing, 100081, P. R.  China.

\noindent  E-mail: yangt@bit.edu.cn


{\small

\begin{thebibliography}{99}

\bibitem{AGG86}  W. Arendt, A. Grabosch, G. Greiner, U. Groh, H. P. Lotz, U. Moustakas, R. Nagel, F. Neubrander and U. Schlotterbeck, \textit{One-Parameter Semigroups of Positive Operators}, Lecture Notes in Math. {\bf 1184}, Springer, Berlin, 1986.

\bibitem{AKSZ} H. Aikawa, T. Kilpel\"{a}inen, N. Shanmugalingam and X. Zhong, Boundary Harnack principle for $p$-harmonic functions in smooth Euclidean domains, \textit{Potential Anal.} {\bf 26} (2007), 281--301.

\bibitem{BBCK}  M. T. Barlow,  R. F. Bass, Z.-Q. Chen and M. Kassmann,
Non-local Dirichlet forms and symmetric jump processes, \textit{Trans. Amer. Math. Soc.} {\bf 361} (2009),
1963--1999.

\bibitem{BCDS06} B. Berkowitz, A. Cortis, M. Dentz and H. Scher, Modeling non-Fickian transport in geological formations as a continuous time random walk, \textit{Rev. Geophys.} {\bf 44} (2006), RG2003.

\bibitem{Bogdan1} K. Bogdan, The boundary Harnack principle for the fractional Laplacian, \textit{Studia Math.} {\bf 123} (1997), 43--80.

\bibitem{BB} K. Bogdan and T. Byczkowski, Potential theory for the $\alpha$-stable Schr\"{o}dinger operator on bounded Lipschitz domains, \textit{Studia Math.} {\bf 133} (1999), 53--92.




\bibitem{BGR1} K. Bogdan, T. Grzywny and M. Ryznar, Heat kernel estimates for the fractional Laplacian with Dirichlet conditions, \textit{Ann. Probab.} {\bf 38} (2010), 1901--1923.

\bibitem{BGR2} K. Bogdan, T. Grzywny and M. Ryznar, Dirichlet heat kernel for unimodal L\'{e}vy processes, \textit{Stochastic Process. Appl.} {\bf 124} (2014), 3612--3650.

\bibitem{BJ} K. Bogdan and T. Jakubowski, Estimates of heat kernel of fractional Laplacian perturbed by gradient operators, \textit{Comm. Math. Phys.} {\bf 271} (2007), 179--198.

\bibitem{BKN} K. Bogdan, T. Kulczycki and A. Nowak, Gradient estimates for harmonic and $q$-harmonic functions of symmetric stable processes, \textit{Illinois J. Math.} {\bf 46} (2002), 541--556.

\bibitem{BG90} J.-P. Bouchaud and A. Georges, Anomalous diffusion in disordered media: statistical mechanisms, models and physical applications, \textit{Phys. Rep.} {\bf 195} (1990), 127--293.

\bibitem{CSS} L. A. Caffarelli, S. Salsa and L. Silvestre, Regularity estimates for the solution and the free boundary of the obstacle problem for the fractional Laplacian, \textit{Invent. Math.} {\bf 171} (2008), 425--461.

\bibitem{Chen1} Z.-Q. Chen, Gaugeability and conditional gaugeability, \textit{Trans. Amer. Math. Soc.} {\bf 354} (2002), 4639--4679.

\bibitem{Chen} Z.-Q. Chen, Symmetric jump processes and their heat kernel estimates, \textit{Sci. China Ser. A} {\bf 52} (2009), 1423--1445.

\bibitem{Chen-Kim-Kumagai} Z.-Q. Chen, P. Kim and T. Kumagai, Weighted Poincar\'{e} inequality and heat kernel estimates for finite range jump processes, \textit{Math. Ann.} {\bf 342} (2008), 833--883.

\bibitem{CKK2} Z.-Q. Chen, P. Kim and T. Kumagai, Global heat kernel estimates for symmetric jump processes, \textit{Trans. Amer. Math. Soc.} {\bf 363} (2011), 5021--5055.

\bibitem{CKS1} Z.-Q. Chen, P. Kim and R. Song, Heat kernel estimates for the Dirichlet fractional Laplacian, \textit{J. Eur. Math. Soc.} {\bf 12} (2010), 1307--1329.

\bibitem{CKS3} Z.-Q. Chen, P. Kim and R. Song, Dirichlet heat kernel estimates for $\Delta^{\alpha/2}+\Delta^{\beta/2}$, \textit{Illinois J. Math.} {\bf 54} (2010), 1357--1392.

\bibitem{CKSrel} Z.-Q. Chen, P. Kim and R. Song, Sharp heat kernel estimates for relativistic stable processes in open sets, \textit{Ann. Probab.} {\bf 40} (2012), 213--244.

\bibitem{CKS2} Z.-Q. Chen, P. Kim and R. Song, Dirichlet heat kernel estimates for fractional Laplacian with gradient perturbation, \textit{Ann. Probab.} {\bf 40} (2012), 2483--2538.

\bibitem{CKSrot} Z.-Q. Chen, P. Kim and R. Song, Dirichlet heat kernel estimates for rotationally symmetric L\'{e}vy processes, \textit{Proc. Lond. Math. Soc.} {\bf 109} (2014), 90--120.


\bibitem{CKS15} Z.-Q. Chen, P. Kim and R. Song, Stability of Dirichlet heat kernel estimates for non-local operators under Feynman--Kac perturbation, \textit{Trans. Amer. Math. Soc.} {\bf 367} (2015), 5237--5270.

\bibitem{CK1} Z.-Q. Chen and T. Kumagai, Heat kernel estimates for stable-like processes on $d$-sets, \textit{Stochastic Process. Appl.} {\bf 108} (2003), 27--62.

\bibitem{Chen-Kumagai} Z.-Q. Chen and T. Kumagai, Heat kernel estimates for jump processes of mixed types on metric measure spaces, \textit{Probab. Theory Related Fields} {\bf 140} (2008),
    277--317.

\bibitem{CRY16} Z.-Q. Chen, Y.-X. Ren and T. Yang, Boundary Harnack principle and gradient estimates for fractional Laplacian perturbed by non-local operators, \textit{Potential Anal.} {\bf 45} (2016), 509--537.

\bibitem{CS2} Z.-Q. Chen and R. Song, Intrinsic ultracontractivity and conditional gauge for symmetric stable processes, \textit{J. Funct. Anal.} {\bf 150} (1997), 204--239.

\bibitem{CS} Z.-Q. Chen and R. Song, Estimates on Green functions and Poisson kernels of symmetric stable processes, \textit{Math. Ann.} {\bf 312} (1998), 465--501.

\bibitem{CW} Z.-Q. Chen and J.-M. Wang, Perturbation by non-local operators, \textit{Ann. Inst. Henri Poincar\'e Probab. Stat.} {\bf 54} (2018), 606--639.


\bibitem{CW2} Z.-Q. Chen and J.-M. Wang, Errata for Perturbation by non-local operators, \textit{Ann. Inst. Henri Poincar\'e Probab. Stat.} {\bf 56} (2020), 760--763.

\bibitem{CZ} Z.-Q. Chen and X. Zhang, Heat kernels and analyticity of non-symmetric jump diffusion semigroups, \textit{Probab. Theory Related Fields} {\bf 165} (2016), 267--312.

 \bibitem{CS25} S. Cho and R. Song, Approximate factorizations for non-symmetric jump processes, arXiv:2504.14763, 2025.

\bibitem{Chung} K. L. Chung and J. B. Walsh, \textit{Markov Processes, Brownian Motion, and Time Symmetry}, 2nd ed., Springer, New York, 2005.

\bibitem{CT04} R. Cont and P. Tankov, \textit{Financial Modelling with Jump Processes}, Chapman \& Hall/CRC,  2004.

\bibitem{Getoor} R. K. Getoor, Duality of L\'{e}vy systems, \textit{Z. Wahrscheinlichkeitstheorie und Verw. Gebiete} {\bf 19} (1971), 257--270.

\bibitem{KS} P. Kim and R. Song, Dirichlet heat kernel estimates for stable processes with singular drift in unbounded $C^{1,1}$ open sets, \textit{Potential Anal.} {\bf 41} (2014), 555--581.

\bibitem{KSZ96} J. Klafter, M. F. Shlesinger and G. Zumofen, Beyond Brownian motion, \textit{Physics Today} {\bf 49} (1996), 33--39.

\bibitem{Ku} T. Kulczycki, Properties of Green function of symmetric stable processes, \textit{Probab. Math. Statist.} {\bf 17} (1997), 381--406.


\bibitem{Liao2} M. Liao, Riesz representation and duality of Markov processes, in: \textit{S\'eminaire de Probabilit\'es XIX}, Lecture Notes in Math. {\bf 1123}, Springer, Berlin, 1985, 366--396.


\bibitem{MK00} R. Metzler and J. Klafter, The random walk's guide to anomalous diffusion: a fractional dynamics approach, \textit{Phys. Rep.} {\bf 339} (2000), 1--77.

\bibitem{Meyer} P.-A. Meyer, Renaissance, recollements, m\'elanges, ralentissement de processus de Markov, \textit{Ann. Inst. Fourier (Grenoble)} {\bf 25} (1975), 465--497.

\bibitem{P} A. Pazy, \textit{Semigroups of Linear Operators and Applications to Partial Differential Equations}, Springer, New York, 1983.

\bibitem{Schaefer} H. H. Schaefer, \textit{Banach Lattices and Positive Operators}, Springer, New York, 1974.

\bibitem{Silvestre} L. Silvestre, Regularity of the obstacle problem for a fractional power of the Laplace operator, \textit{Comm. Pure Appl. Math.} {\bf 60} (2007), 67--112.

\bibitem{SY22} Y. Shi and B. Yi, Intrinsic ultracontractivity of Laplacian and fractional Laplacian perturbed by non-local operators, \textit{Sci. Sin. Math.} {\bf 52} (2022), 1307--1322.

\bibitem{Song} R. Song, Estimates on the Dirichlet heat kernel of domains above the graphs of bounded $C^{1,1}$ functions, \textit{Glas. Mat. Ser. III} {\bf 39} (2004), 275--288.

\end{thebibliography}
}
\end{document}